\documentclass[10pt]{ijnam}
\usepackage{amsmath,amssymb,amsfonts}
\usepackage{xcolor}
\usepackage{algorithm} 
\usepackage{algorithmic} 
\usepackage{pifont}
\usepackage{hyperref}
\usepackage{cleveref}
\usepackage{graphicx}
\usepackage{geometry}
\def\currentvolume{1}
\def\currentissue{1}
\def\currentyear{2024}
\def\month{March}
\pagespan{1}{30}
\copyrightinfo{2026}{} 

\newcommand{\be}{\begin{equation}}
\newcommand{\ee}{\end{equation}}

\newtheorem{theorem}{Theorem}[section]

\newtheorem{lem}[theorem]{Lemma}

\theoremstyle{remark}
\newtheorem{remark}[theorem]{Remark}
\newtheorem{example}[theorem]{Example}

\begin{document}
	
	\title[IEQ-FEMs for CHNS equations]
	{Unconditional energy stable hybrid IEQ-FEMs for the Cahn-Hilliard-Navier-Stokes equations}
	
	\author[Y. Chen]{Yaoyao Chen}
	\address{
		School of Mathematics and Statistics, Anhui Normal University, Wuhu, Anhui 241000, PR China
	}
	\email{cyy1012xtu@126.com}
	
	\author[D. Li]{Dongqian Li}
	\address{
	Hunan Key Laboratory for Computation and Simulation in Science and Engineering, Hunan International Scientific and Technological Innovation Cooperation Base of Computational Science, Xiangtan University, Xiangtan 411105, Hunan,  China
	}
	\email{zuel\_ldq@outlook.com}

	\author[Y. Yang]{Yin Yang*}
	\address{
	Hunan Research Center of the Basic Discipline Fundamental Algorithmic Theory and Novel Computational Methods, Key Laboratory of Intelligent Computing and Information Processing of Ministry of Education, Xiangtan University, Xiangtan 411105, Hunan, China
	}
	\email{yangyinxtu@xtu.edu.cn}

	\author[P. Yin]{Peimeng Yin}
	\address{
	Department of Mathematical Sciences, The University of Texas at El Paso, El Paso, TX 79968, USA
	}
	\email{pyin@utep.edu}
	
	
	
	\date{December 9, 2024 and, in revised form, January 30, 2026.}
	
	\thanks{*Corresponding author.}
	
	
	\abstract{We investigate two unconditionally energy stable invariant energy quadratization (IEQ) finite element methods (FEMs) [Chen et al. Numerical Algorithms, 99:1161-1202, 2025] for solving the Cahn-Hilliard-Navier-Stokes (CHNS) equations. The time discretization of these IEQ-FEMs is based on the first- and second-order backward differentiation methods. {The auxiliary energy function introduced by the IEQ approach, modeling the square root of the nonlinear part of the energy, does not belong to the finite element space used for the spatial discretization.} These methods offer distinct advantages. Consequently, we propose a new hybrid IEQ-FEM that combines the strengths of both schemes, offering computational efficiency and unconditional energy stability in the finite element space. We provide rigorous proofs of mass conservation and energy dissipation for the proposed IEQ-FEMs. Several numerical experiments are presented to validate the accuracy, efficiency, and solution properties of the proposed method.}
	
	
	\maketitle
	
\section{Introduction}
This paper focuses on the development of unconditionally energy-stable finite element methods (FEMs) based on the invariant energy quadratization (IEQ) approach for solving the Cahn-Hilliard-Navier-Stokes (CHNS) problem
\cite{lowengrub_quasiincompressible_1998,shen_energy_2010}
\begin{subequations}\label{eq1.3}
	\begin{align}
		{\partial _t}\phi  + \nabla  \cdot \left( {{\mathbf{u}}\phi } \right) - \gamma \Delta w = 0\qquad &{\it \text{in}}\ \Omega  \times \textcolor{black}{(0,T]}, \label{eq1.3.1}\\
		w + \lambda \left( {\Delta \phi  - f\left( \phi  \right)} \right) = 0\qquad &{\it \text{in}}\ \Omega  \times \textcolor{black}{(0,T]},\label{eq1.3.2} \\
		{\partial _t}{\mathbf{u}} - \mu \Delta {\mathbf{u}} + \left( {{\mathbf{u}} \cdot \nabla } \right){\mathbf{u}} + \nabla p + \phi \nabla w = 0 \qquad &{\it \text{in}}\ \Omega  \times \textcolor{black}{(0,T]},\label{eq1.3.3}\\
		\nabla  \cdot {\mathbf{u}} = 0\qquad &{\it \text{in}}\ \Omega  \times \textcolor{black}{(0,T]},\label{eq1.3.4}\\
		{\mathbf{u}} = {\textcolor{black}{\mathbf{0}}},\ \dfrac{\partial \phi }{\partial {\textcolor{black}{\nu}}} = 0,\ {\dfrac{{\partial w}}{\partial {\textcolor{black}{\nu}}} = 0}
		\qquad &{\it \text{on}}\ \partial \Omega  \times \textcolor{black}{(0,T]},\label{eq1.3.5}\\
		{\mathbf{u}}\left( { \cdot ,0} \right) = {{\mathbf{u}}_0},\ \phi \left( { \cdot ,0} \right) = {\phi _0}\qquad &{\it \text{in}}\ \Omega,\label{eq1.3.6}
	\end{align}
\end{subequations} 
where $\Omega\subseteq{\mathbb R}^d(d=2,3)$ is a bounded domain, \textcolor{black}{and $\nu$ represents the unit outward normal vector on the boundary $\partial\Omega$.}  Here, $\gamma$ is a mobility constant related to the relaxation time scale, $\mu$ denotes the viscosity, and $\lambda$ represents the magnitude of the mixing energy. $\mathbf{u}$ and $p$ denote the velocity field and pressure, respectively. $\phi$ represents the phase field variable, where $\phi=\pm 1$ corresponds to two different fluids, and $w$ is the chemical potential. The nonlinear term $f(\phi)=F^{'}(\phi)$, where $F(\phi)$ represents the nonlinear bulk potential. The Dirichlet boundary conditions are imposed on the velocity field $\mathbf{u}$, while Neumann boundary conditions are applied to the phase field $\phi$ and potential $w$, as specified in \eqref{eq1.3.5}. One of well-known potentials is the double well potential 
\be\label{dwp}
F(\phi)=\frac{1}{4\epsilon^2}(\phi^2-1)^2,
\ee
where the parameter $\epsilon  > 0$. 
The CHNS equations \eqref{eq1.3} is endowed with energy dissipation law
\begin{equation}
	\frac{d}{{dt}}\mathcal{E}\left( {\phi ,{\mathbf{u}}} \right) =  - \int_\Omega  \left({\mu {\left| {\nabla {\mathbf{u}}} \right|^2} + \gamma {\left| {\nabla w} \right|^2}} \right) d\mathbf{x} \le 0,
\end{equation}
with the total energy $\mathcal{E}\left( {\phi ,{\mathbf{u}}} \right) =\int_\Omega \left(  { \frac{\lambda }{2}{{\left| {\nabla \phi } \right|}^2} + \lambda F\left( \phi  \right)}  +  {\frac{1}{2}{{\left| {\mathbf{u}} \right|}^2} } \right) d\textcolor{black}{\mathbf{x}}. $

The physical phenomena captured by the two-phase flow model frequently occur in nature and various industrial processes \cite{zhu_pore-scale_2017}. Over the past few decades, the CHNS equations and its modifications have been applied to numerous situations, such as the solidification of liquid metal alloys \cite{ENWR2013}, and the simulation of bubble dynamics \cite{ABHKWW2013}.

There are many challenges in developing efficient and easy-to-implement numerical schemes for solving the phase-field CHNS equations. 
A key difficulty arises from the higher-order derivatives in the Cahn-Hilliard (CH) equations. One approach to address this is the use of the mixed methods, which reformulates the fourth-order phase field problem into a lower-order system, making it more amenable to numerical approximation. 
In recent years, several numerical schemes have been actively designed and applied to solve the CH equation or CHNS equations, including the FEMs \cite{chen_fully-discrete_2021,chen_error_2022, chen2024recovery}, finite difference methods (FDMs) \cite{zhao_second-order_2021,zhu_efficient_2019}, and spectral methods \cite{chen_highly_2022,shen_energy_2010,yang_new_2021,yang_efficient_2018,ye_efficient_2022}.

Another significant challenge is handling nonlinear terms, including the convection term $ \left( {{\mathbf{u}} \cdot \nabla } \right){\mathbf{u}}$, the coupling term $\phi \nabla w$ in the NS equation, and the convection term $\nabla  \cdot \left( {{\mathbf{u}}\phi } \right)$ in the CH equation. For these terms, explicit-implicit methods are commonly used \cite{yang_efficient_2018}, and the projection method is applied to enforce a divergence-free velocity field \cite{guermond_overview_2006,JShenPROJ1992}.
Additionally, efficiently managing the nonlinear term $f(\phi)$ in the CH equation poses another challenge.
This term is often addressed using various techniques, such as \textcolor{black}{the convex-splitting method \cite{cai2023optimal, han2015second}}, the stabilization method \cite{chen_decoupled_2021, shen_energy_2010,xu_error_2020}, the IEQ approach \cite{chen_IEQ_FEM, liu_unconditionally_2021,yang_efficient_2018,yang_convergence_2020}, the scalar auxiliary variable (SAV) approach \cite{li_unconditional_2022,shen_new_2019,zhu_efficient_2019}, and other related methods based on IEQ or SAV apparoach \cite{chen_fully-discrete_2021,chen_highly_2022,yang_family_2019,ye_efficient_2022}. 

\textcolor{black}{FEMs offer flexibility, accuracy, and the ability to handle complex geometries, nonlinearities, and multi-physics systems, making them widely used for  solving the CHNS equations \cite{chen_fully-discrete_2021,chen_error_2022,YANGJT2024115577}. 
	A straightforward idea for combining the IEQ approach with FEMs is to discretize the auxiliary variable introduced by IEQ directly within the finite element space. However, this leads to a scheme that is computationally expensive due to the need to invert the mass matrix \cite{chen_IEQ_FEM}.
	The choice of function space for the auxiliary variable plays a crucial role in both computational efficiency and, in some cases, the physical fidelity of the numerical solution.
	Comparative studies on various IEQ-FEM schemes for the Cahn–Hilliard and Allen–Cahn equations were carried out in \cite{chen_IEQ_FEM}, identifying two efficient approaches. In the first approach, the intermediate function is constrained to the continuous function space $C^0(\Omega)$, which enables efficient computation. However, because the computed energy cannot be directly evaluated in $C^0(\Omega)$, it must be interpolated into the finite element space. This may result in non-monotonic energy decay and potentially nonphysical behavior, despite the method's efficiency.
	The second approach computes the intermediate function in $C^0(\Omega)$ and then projects it into the finite element space. This projection ensures that the discrete energy law is unconditionally satisfied, albeit at the expense of additional computational cost.
	Motivated by these two methods, we extend them to the CHNS system and propose a novel hybrid IEQ-FEM scheme that combines their strengths—achieving both computational efficiency and unconditional energy stability.}

Solving the CHNS equations using the IEQ approach within the finite element method framework has not been thoroughly explored in the literature. In this work, we study several IEQ-FEMs for solving the CHNS equations. Specifically, we introduce two initial types of IEQ-FEMs, where the intermediate function is either in the combination of the polynomial space and the continuous function space, or only continuous function space.
In the method where the intermediate function is in the combination of the polynomial space and the continuous function space, the function is first computed in the continuous function space $C^0(\Omega)$ and then projected onto the finite element space. This approach ensures that the numerical solutions are unconditionally energy-stable. The second method situates the intermediate function directly in the continuous function space $C^0(\Omega)$ rather than the finite element space. {\color{black}Since the intermediate function does not lie in the finite element space, the computed energy cannot
be exactly evaluated in $C^0(\Omega)$ and must instead be interpolated within the finite element space. However, due to interpolation errors, the computed energy may not consistently exhibit decay, potentially leading to nonphysical solutions \cite{chen_IEQ_FEM}. Nonetheless, conditional energy decay in finite element space can still be achieved. }

Building on these, we propose a hybrid IEQ-FEM that begins with the IEQ-FEM with the intermediate function in the continuous function space. If the computed energy fails to exhibit dissipation, the scheme switches to the IEQ-FEM with projection onto the finite element space. This hybrid IEQ-FEM combines the advantages of both earlier methods: computational efficiency and unconditional energy stability in finite element space.
Additionally, the phase field and velocity variables are approximated using piecewise quadratic finite elements, while the pressure variable is approximated by using piecewise linear finite elements. For time discretization, we apply first-order and second-order backward differentiation methods. 

The rest of the paper is organized as follows. In Section \ref{se:2}, we first introduce the semi-discrete FEM for the CHNS equations, followed by an equivalent PDE system reformulated using the IEQ approach. In Section \ref{se:3}, we construct first- and second-order fully discrete IEQ-FEM schemes for the CHNS equations with an intermediate function in different function spaces. We rigorously prove the existence and uniqueness of the numerical solutions, along with the energy dissipation property. In Section \ref{se:4}, we present several numerical examples to demonstrate the effectiveness of the proposed schemes in solving the CHNS equations.

\section{IEQ-FEMs for the CHNS Equations} 
\label{se:2}

In this section, we present the IEQ-FEMs for the CHNS equations (\ref{eq1.3}). We begin by introducing the notations that will be used throughout the paper.
Let $\alpha=(\alpha_1, \ldots, \alpha_d)\in\mathbb Z^d_{\geq0}$ be a multi-index with $\partial^\alpha:=\partial_1^{\alpha_1}\cdots\partial^{\alpha_d}_{x_d}$ and $|\alpha|:=\sum_{i=1}^d\alpha_i$.
The Sobolev space $H^m(\Omega)$, $m\geq 0$, consists of functions whose derivatives, corresponding to the multi-index $\alpha$, are square integrable. Denote by $H^1_0(\Omega)\subset H^1(\Omega)$ the subspace consisting of functions with zero trace on the boundary  $\partial\Omega$. Let $L^2(\Omega):=H^0(\Omega)$, and let $L^2_0(\Omega)$ be the subset of $L^2(\Omega)$ consisting of functions with zero average. For $s\geq 0$, $(H^s(\Omega))^d$ represents the vector space, such that $\mathbf{v}=(v_1,\ldots, v_d)^\top \in (H^s(\Omega))^d$ represents $v_i \in H^s(\Omega)$, $i=1,\ldots, d$, where $\top$ is the transposition of a matrix or a vector. 

Let $\mathcal{T}_{h}$ be a regular shape partition of $\Omega$, with each $K \in \mathcal{T}_{h}$ being a quasi-uniform element in the finite element mesh $\mathcal{T}_{h}$. The diameter of each element $K$ is denoted by $h:=\max\{h_K|h_K=\text{diam}(K), K\in\mathcal{T}_{h}\}$. 
Define $Y_{h}$ as a finite dimensional subspace of $H^{1}(\Omega)$, and 
${\mathbf{X}}_{h}\times M_{h}\subset (H_0^1(\Omega))^d\times L_0^2(\Omega)$ as a pair of mixed finite element spaces, which satisfies the following $inf$-$sup$ condition:
\begin{equation} \label{infsup}
\inf_{q\in M_{h}}\sup_{\boldsymbol{v}\in {\mathbf{X}}_{h}}\frac{\left(\mathbf{B}^\top q,\boldsymbol{v}\right)}{\|\boldsymbol{v}\|_{H^{1}(\Omega)}\|q\|_{L^{2}(\Omega)}}\geq \beta_{\ast},
\end{equation}
where $\beta_{\ast}>0$ is a constant, and the operator $\mathbf{B}^\top: M_{h}\rightarrow {\mathbf{X}}_{h}$ is the transpose of the discrete divergence operator $\mathbf{B}: {\mathbf{X}}_{h}\rightarrow M_{h}$. For every pair $(\boldsymbol{v},\Phi)\in {\mathbf{X}}_{h}\times M_{h}$, it holds
\[(\mathbf{B}\boldsymbol{v},\Phi)=(\boldsymbol{v},\mathbf{B}^\top\Phi)=-(\nabla\cdot\boldsymbol{v},\Phi).\]
Let ${\mathbf{V}}_{h}$ be a finite-dimensional subspace of $(L^{2}(\Omega))^{d}$. We assume either that $\mathbf{V}_{h}$ is conformal in $H_{0}^{\mathrm{div}}(\Omega)$
or that $M_{h}$ is conformal in $H^{1}(\Omega)$. 
In this work, following \cite{cai2018error,guermond1998approximation}, we choose the finite element spaces as
\begin{equation*}
\begin{aligned}
	&{{\mathbf{X}}_h} \textcolor{black}{= {{\mathbf{V}}_h}} = \left\{ {{{\mathbf{v}}_h} \in {{\left( {{C^0}\left( {\bar \Omega } \right)} \right)}^d} \cap {{\left( {H_0^1(\Omega)} \right)}^d};{{\mathbf{v}}_h}\left| {_K} \right. \in {{\left( {{P_2}\left( K \right)} \right)}^d}} \right\},\\
	&{M_h} = \left\{ {{\phi _h} \in C^0(\bar \Omega)\cap L_0^2(\Omega);{\phi _h}\left| {_K} \right. \in {P_1}\left( K \right)} \right\},\\
	&{Y_h} = \left\{ {\psi  \in {C^0}\left( {\bar \Omega } \right);{\psi _h}\left| {_K} \right. \in {P_2}\left( K \right)} \right\} . 
\end{aligned} 
\end{equation*}

\subsection{The semi-discrete FEM for the CHNS equations.} 
The semi-discrete finite element scheme for the
CHNS equations (\ref{eq1.3}) is to find $(\phi_h, w_h,{\mathbf{u}}_h,p_h)\in Y_h\times Y_h\times {\mathbf{X}}_h\times M_h $ such that
\begin{subequations}\label{eq2.1}
\begin{align}
	\left({\partial _t}\phi_h,\varphi\right)  \textcolor{black}{- \left( {{\mathbf{u}_h}\phi_h } ,\nabla \varphi\right)} - \gamma\left( \nabla w_h,\nabla \varphi\right) = 0,\quad & {\forall\varphi\in Y_h}, \label{eq2.1.1}\\
	\left(w_h,\psi\right)- \lambda \left( \nabla  \phi_h,\nabla \psi\right)  - \lambda \left(f\left( \phi_h  \right),\psi \right) = 0,\quad & {\forall\psi \in Y_h},\label{eq2.1.2} \\
	\left({\partial _t}{\mathbf{u}_h},{\mathbf{v}}\right) + \mu \left(\nabla {\mathbf{u}_h},\nabla {\mathbf{v}}\right) + \left(\left( {{\mathbf{u}_h} \cdot \nabla } \right){\mathbf{u}_h},{\mathbf{v}}\right) + \left(\nabla p_h,{\mathbf{v}}\right) + \left(\phi_h \nabla w_h,{\mathbf{v}}\right) = 0, \quad &{\forall{\mathbf{v}} \in {\mathbf{X}_h}},\label{eq2.1.3}\\
	\left(\nabla  \cdot {\mathbf{u}}_h,q\right) = 0,\quad &{\forall q \in M_h}.\label{eq2.1.4}
\end{align}
\end{subequations}
\textcolor{black}{Here and in what follows, we denote by $\Pi$ and $\mathbf{\Pi}$ the $L^2$ projection operators for scalar-valued and vector-valued functions, respectively. For any $\phi \in L^2(\Omega)$ and $\mathbf{u} \in (L^2(\Omega))^d$, the projections $\Pi \phi \in Y_h$ and $\mathbf{\Pi} \mathbf{u} \in \mathbf{X}_h$ are defined by
\begin{subequations}
	\begin{align}
		\int_\Omega (\Pi \phi(x) - \phi(x)) \, \varphi \, d\mathbf{x} &= 0, \quad \forall \varphi \in Y_h, \label{PJ1}\\
		\int_\Omega (\mathbf{\Pi} \mathbf{u}(x) - \mathbf{u}(x)) \cdot \mathbf{v} \, d\mathbf{x} &= 0, \quad \forall \mathbf{v} \in \mathbf{X}_h. \label{PJ2}
	\end{align}
\end{subequations}
The initial conditions for the semi-discrete finite element scheme \eqref{eq2.1} are given by
\begin{equation} \label{phiuInitial}
	{\mathbf{u}_h}\left( { \cdot ,0} \right) = \mathbf{\Pi}{{\mathbf{u}}_0},\ \phi_h \left( { \cdot ,0} \right) = \Pi{\phi _0}.
	\end{equation}}
	We introduce the discrete free energy
	\be\label{freeenergy}
	\mathcal{E}\left( {\phi_h ,{\mathbf{u}_h}} \right)  = \int_\Omega  {\left( {\frac{1}{2}{{\left| {\mathbf{u}_h} \right|}^2} + \frac{\lambda }{2}{{\left| {\nabla \phi_h } \right|}^2} + \lambda F\left( \phi_h  \right)} \right)} d\mathbf{x} .
	\ee
	Then the following results hold.
	
	\begin{lem}\label{CHsemiLem}
The semi-discrete finite element scheme \eqref{eq2.1} conserves the total mass
\be\label{mass0}
\frac{d}{dt}\int_{\Omega }\phi_{h} d\mathbf{x} = 0,
\ee
and the solution satisfies the energy dissipation law
\be\label{energy0}
\frac{d}{{dt}}\mathcal{E}\left( {\phi_h ,{\mathbf{u}_h}} \right) =  - \int_\Omega  \left({\mu {\left| {\nabla {\mathbf{u}_h}} \right|^2} + \gamma {\left| {\nabla w_h} \right|^2}} \right) d\mathbf{x} \leq 0.
\ee
\end{lem}

\subsection{The IEQ reformulation} \label{section2.2IEQ}
Note that the potential $F(\phi_h)$ is uniformly bounded from below. The IEQ approach \cite{X.Yang2016} transforms $F(\phi_h)$ into a quadratic form through an intermediate function
$$
U = \sqrt {F\left( \phi_h \right) + B},
$$
where $B$ is a constant to ensure $F(\phi_h)+B> 0$. Consequently, the free energy \eqref{freeenergy} becomes
\begin{equation}\label{se2-1}
\mathcal{E}\left( {\phi_h ,{\mathbf{u}_h},U} \right) = \int_\Omega  {\left( {\frac{1}{2}{{\left| {\mathbf{u}_h} \right|}^2} + \frac{\lambda }{2}{{\left| {\nabla \phi_h } \right|}^2} + \lambda U^2} \right)} d\mathbf{x} - \lambda B \vert \Omega \vert .
\end{equation}
By using the intermediate function $U$, the nonlinear term $F^{'}(\phi_h)$ is substituted by $$F^{'}(\phi_h)= H\left( \phi_h \right) U,$$ where $$H\left( \nu\right) =\frac{F^{'}(\nu) }   {\sqrt{F \left(\nu \right) +B}}.$$ 
The update of $U$ is governed by 
\be\label{uode}
\partial_t U = \frac{1}{2}H(\phi_h) \partial_t \phi_{h},
\ee
\textcolor{black}{which is an ordinary differential equation for $U$,} subject to the initial data
\be\label{CHinit+}
U(x, 0)=\sqrt{F(\phi_0(x))+B}.
\ee
With the IEQ approach, the semi-discrete finite element scheme \eqref{eq2.1} can be reformulated as the following semi-discrete IEQ-FEM scheme
\begin{subequations}\label{eq2.1+}
\begin{align}
	\left({\partial _t}\phi_h,\varphi\right)  \textcolor{black}{- \left( {{\mathbf{u}_h}\phi_h } ,\nabla \varphi\right)} - \gamma\left( \nabla w_h,\nabla \varphi\right) = & 0,\quad  {\forall\varphi\in Y_h}, \label{eq2.1.1+}\\
	\left(w_h,\psi\right)- \lambda \left( \nabla  \phi_h,\nabla \psi\right)  - \lambda \left(H\left( \phi_h \right) U,\psi \right) = & 0,\quad  {\forall\psi \in Y_h},\label{eq2.1.2+} \\
	\left({\partial _t}{\mathbf{u}_h},{\mathbf{v}}\right) + \mu \left(\nabla {\mathbf{u}_h},\nabla {\mathbf{v}}\right) + \left(\left( {{\mathbf{u}_h} \cdot \nabla } \right){\mathbf{u}_h},{\mathbf{v}}\right) + \left(\nabla p_h,{\mathbf{v}}\right) + \left(\phi_h \nabla w_h,{\mathbf{v}}\right) = & 0, \quad {\forall{\mathbf{v}} \in {\mathbf{X}_h}},\label{eq2.1.3+}\\
	\left(\nabla  \cdot {\mathbf{u}}_h,q\right) = & 0,\quad {\forall q \in M_h},\label{eq2.1.4+}\\
	\partial_t U = & \frac{1}{2}H(\phi_h) \partial_t \phi_{h}.\label{eq2.1.5+}
\end{align}
\end{subequations}
Then the following results hold.
\begin{lem}\label{CHsemiLem+}
The semi-discrete IEQ-FEM scheme \eqref{eq2.1+} conserves the total mass
\be\label{mass0+}
\frac{d}{dt}\int_{\Omega }\phi_{h} d\mathbf{x} = 0,
\ee
and the solution satisfies the energy dissipation law
\be\label{energy0+}
\frac{d}{{dt}}\mathcal{E}\left( {\phi_h ,{\mathbf{u}_h}},U \right) =  - \int_\Omega  \left({\mu {\left| {\nabla {\mathbf{u}_h}} \right|^2} + \gamma {\left| {\nabla w_h} \right|^2}} \right) d\mathbf{x} \leq 0.
\ee
\end{lem}
\begin{proof}
Choosing  $\varphi=1$ in (\ref{eq2.1.1+}) gives the total mass conservation (\ref{mass0+}).
From (\ref{eq2.1.5+}), it follows 
\be\label{HUtoU}
(H(\phi_h)U, \partial_t\phi_h) = (U, H(\phi_h) \partial_t\phi_h) = 2(U, \partial_tU) = \partial_t \|U\|^2.
\ee
By setting $\varphi=-w_h$ in (\ref{eq2.1.1+}), $\psi=\partial_t \phi_{h}$ in (\ref{eq2.1.2+}), and ${\mathbf{v}}={\mathbf{u}_h}$ in (\ref{eq2.1.3+}), the summation of (\ref{eq2.1.1+}), (\ref{eq2.1.2+}), and (\ref{eq2.1.3+}) together with \eqref{HUtoU} yields (\ref{energy0+}).
\end{proof}

\section{Fully discrete schemes}\label{se:3}
The main purpose of this section is to present the IEQ-FEM fully discrete schemes, where the spatial discretization is based on the finite element method, and the first- and second-order backward differentiation (BDF1 and BDF2 for short) are applied for the temporal discretization.


We consider a partition $0=t_{0}<t_{1}<\cdots<t_{N}=T$ of $[0,T]$, the $(n+1)$-th subinterval is defined as $I_{n+1}:=(t_n, t_{n+1}]$, and the corresponding time step $\tau_{n+1}:=t_{n+1}-t_{n},  n=0, \ldots, N-1$. For simplicity, we assume the time step $\tau_n = \tau>0$ is a constant.
For any given (vector) function $v(x,t)$ and $n\geq 0$, we denote $v^n:= v(x, t_n)$ or $v_h^n$ the approximation of $v(x,t)$ at $t_n$.

\subsection{first-order fully discrete schemes} In this section, we first introduce two types of first-order fully discrete IEQ-FEM schemes for the CHNS equations, where the intermediate function $U$ in the semi-discrete IEQ-FEM scheme \eqref{eq2.1+} is positioned in either the polynomial space or the continuous function space. We then propose a new method that combines both of these schemes.

\subsubsection{\textcolor{black}{P-BDF1} scheme}
To design a first-order, fully decoupled, and fully discrete IEQ-FEM scheme, we apply a splitting strategy to the CHNS problem \eqref{eq1.3}. Following the idea in \cite{GuermondShen2003splitting}, the pressure $p$ in \eqref{eq1.3.3} can be computed separately from the velocity field $\mathbf{u}$. 
Assume that $w$, and $\phi$ has been computed. 
First, if an approximation $p(=\tilde p)$ of pressure is available, we can directly solve for $\mathbf{u}$ from 
\begin{equation}\label{ptou}
{\partial _t}{{\mathbf{u}}} - \mu \Delta {{\mathbf{u}}} + \left( {{\mathbf{u}} \cdot \nabla } \right){\mathbf{u}} + \nabla \tilde p + \phi \nabla w = 0.
\end{equation}
Note that $\mathbf{u}$ from \eqref{ptou} is not necessarily divergence-free.
Second, if $\mathbf{u}$ is known, the pressure $p$ can be recovered by testing \eqref{eq1.3.3} with $\nabla q $. Using $(\partial_t\mathbf{u}, \nabla q) = - (\nabla \cdot \partial_t \mathbf{u}, q)=0$ yields
\be\label{sp}
\left(\nabla  p,\nabla q\right)=\left(\mu \Delta {\mathbf{u}} - \left( {{\mathbf{u}} \cdot \nabla } \right){\mathbf{u}} - \phi \nabla w,\nabla q\right), \quad \forall q \in H^1(\Omega),
\ee
which corresponds to the relation
\be\label{sp+}
\nabla \bar p:=\mu \Delta {\mathbf{u}} - \left( {{\mathbf{u}} \cdot \nabla } \right){\mathbf{u}} - \phi \nabla w.
\ee
Observe that the true pressure $p$ in \eqref{eq1.3.3} satisfies  
\begin{equation}\label{sp+=}
\nabla p= -{\partial _t}{\mathbf{u}} + \mu \Delta {\mathbf{u}} - \left( {{\mathbf{u}} \cdot \nabla } \right){\mathbf{u}} - \phi \nabla w.
\end{equation}
Subtracting \eqref{sp+} from \eqref{sp+=} yields
\be\label{spd}
{\partial _t}{\mathbf{u}} + \nabla (p -\bar p) = 0,
\ee
so enforcing the incompressible condition \eqref{eq1.3.4} on \eqref{spd} provides a corrected (divergence-free) pair $(\mathbf{u},p)$ that serves as an alternative solution to \eqref{eq1.3.3}. 

Hence, the solution $(\mathbf{u},p)$ at $t_{n+1}$ to the NS equation \eqref{eq1.3.3} can be approximated from the following first-order splitting scheme:
\begin{subequations}\label{splitting}
\begin{align}
	& {\partial _t}{\tilde{\mathbf{u}}} - \mu \Delta {\tilde{\mathbf{u}}} + \left( {\tilde{\mathbf{u}} \cdot \nabla } \right)\tilde{\mathbf{u}} + \nabla p^n + \phi \nabla w = 0, \quad \tilde{\mathbf{u}}(x, t_n) =  {\mathbf{u}}^n, \label{splitting-1}\\
	& {\partial _t}{\mathbf{u}} + \nabla (p -p^n) = 0, \quad \quad \mathbf{u}(x, t_n) = \tilde{\mathbf{u}}(x, t_{n+1}).\label{splitting-2}
\end{align}
\end{subequations}
Here, we slightly abuse the notation ${\mathbf{u}}^n$ in \eqref{splitting-1} and $\mathbf{u}(x, t_n)$ in \eqref{splitting-2}. Although they should be identical in meaning, ${\mathbf{u}}^n$ in \eqref{splitting-1} denotes the velocity obtained from the previous time step, while $\mathbf{u}(x, t_n)$ in \eqref{splitting-2} represents the initial condition of the current subproblem at time $t_n$.
Here, the first substep \eqref{splitting-1} advances an intermediate velocity \(\tilde{\mathbf{u}}\) using the known pressure \(p^n\), 
while the second substep \eqref{splitting-2} updates the pressure and corrects the velocity field to enforce the incompressible condition \eqref{eq1.3.4}. 
The intermediate velocity \(\tilde{\mathbf{u}}\) is generally not divergence-free, whereas the corrected velocity \(\mathbf{u}\) at \(t_{n+1}\) satisfies \(\nabla\cdot\mathbf{u}=0\).

The first-order fully discrete IEQ-FEM for \eqref{eq1.3} with \eqref{eq1.3.3} split by \eqref{splitting} is now ready to present. We first project the intermediate function $U^n$ onto the polynomial space $Y_h$.
Given $(\phi^{n}_h,{\mathbf{u}}_h^{n}) \in Y_h\times {\mathbf{V}}_h$ and $U^n \in C^0(\Omega)$, 
the \textcolor{black}{P-BDF1 scheme} is to find $(\phi_h^{n+1}, w_h^{n+1}, {{\mathbf{\tilde u}}_h^{n + 1}, {\mathbf{u}}_h^{n+1}},p_h^{n+1}) \in Y_h  \times Y_h\times  {\mathbf{X}}_h\times {\mathbf{V}}_h \times M_h$ and $U^{n+1} \in C^0(\Omega)$, \textcolor{black}{$\mathbf{\hat u}^{n + 1} \in (L^2(\Omega))^d $} such that
\begin{subequations}\label{se3-1}
\begin{align}
\left( {\dfrac{{\phi _h^{n + 1} - \phi _h^n}}{\tau },{\varphi _h}} \right) \textcolor{black}{-\left( { {{\mathbf{\hat u}}^{n + 1}\phi _h^n} ,\nabla {\varphi _h}} \right)} + \gamma a\left( {w_h^{n + 1},{\varphi _h}} \right) = 0, \qquad &\forall {\varphi _h} \in {Y_h}, \label{se3-1a} \\
\left( {w_h^{n + 1},{\psi _h}} \right) - \lambda a\left( {\phi _h^{n + 1},{\psi _h}} \right)- \lambda \left( {H\left( {\phi _h^n} \right)U^{n+1},{\psi _h}} \right) = 0, \qquad & {\forall {\psi _h} \in {Y_h}},\label{se3-1b}\\
{\mathbf{\hat u}^{n + 1} =  {{\mathbf{u}}_h^n}  - \tau  {\phi _h^n\nabla w_h^{n + 1}} }, \qquad & \label{se3-1c}\\
\left( {U_h^n,{\mu _h}} \right) = \left( {U^n,{\mu _h}} \right),\qquad & {\forall {\mu _h} \in {Y_h}}, \label{se3-1 projection} \\
U^{n + 1} = U_h^n + \frac{1}{2}H\left( {\phi _h^n} \right)\left( {\phi _h^{n + 1} - \phi _h^n} \right),\qquad& \label{se3-1e}
\end{align}
\end{subequations}
\begin{equation}\label{se3-2}
\begin{aligned}
\left( {\dfrac{{{\mathbf{\tilde u}}_h^{n + 1} - {\mathbf{u}}_h^n}}{\tau },{{\mathbf{v}}_h}} \right) + \mu \tilde a\left( {{\mathbf{\tilde u}}_h^{n + 1},{{\mathbf{v}}_h}} \right) + b\left( {{\mathbf{u}}_h^n,{\mathbf{\tilde u}}_h^{n + 1},{{\mathbf{v}}_h}} \right) 
- \left( {p_h^n,\nabla  \cdot {{\mathbf{v}}_h}} \right)& \\+ \left( {\phi _h^n\nabla w_h^{n + 1},{{\mathbf{v}}_h}} \right) = 0,& \qquad\forall {{\mathbf{v}}_h} \in {{\mathbf{X}}_h},
\end{aligned}
\end{equation}
and
\begin{equation}\label{se3-3}
\begin{aligned}
\left( {\dfrac{{{\mathbf{u}}_h^{n + 1} - {\mathbf{\tilde u}}_h^{n + 1}}}{\tau },{{\mathbf{\chi }}_h}} \right) - \left( { \left( {p_h^{n + 1} - p_h^n} \right),{\nabla \cdot {\mathbf{\chi }}_h}} \right) = 0,\qquad&\forall {{\mathbf{\chi }}_h} \in {{\mathbf{V}}_h},\\
\left( {\nabla  \cdot {\mathbf{u}}_h^{n + 1},{q_h}} \right) = 0,\qquad&\forall {q_h} \in {M_h},
\end{aligned}
\end{equation}
where
\begin{equation}
\begin{aligned}
a(w,\phi)=(\nabla w, \nabla \phi),\qquad &\forall w,\phi \in H^1(\Omega), \\ 
\tilde a(\mathbf{u},\mathbf{v})=(\nabla \mathbf{u}, \nabla \mathbf{v}),\qquad &\forall \mathbf{u},\mathbf{v} \in (H^1(\Omega))^d, \\
b\left({\mathbf{u}},{\mathbf{v}},{\mathbf{w}} \right) = ((\mathbf{u}\cdot \nabla) \mathbf{v}, \mathbf{w}) + \frac{1}{2} ((\nabla \cdot \mathbf{u})\mathbf{v},\mathbf{w}), \qquad &\forall \mathbf{u},\mathbf{v} ,\mathbf{w}\in (H^1(\Omega))^d.
\end{aligned}
\end{equation}
\textcolor{black}{The initial data $(\phi^0_h,\mathbf{u}_h^0,U^0)$ is given by \eqref{phiuInitial} and \eqref{CHinit+}.} 
\begin{remark} 
In the P-BDF1 scheme, the projection step \eqref{se3-1 projection} represents an $L^2$ projection, while \eqref{se3-3} corresponds to the projection step used in the Navier-Stokes equations. The latter realizes the identity $\mathbf{u}_h^{n+1} = P_H \mathbf{\tilde u}_h^{n+1}$, where the operator $P_H: \mathbf{X}_h \rightarrow \text{ker}(\mathbf{B})$ \cite{guermond1998approximation}. To avoid ambiguity, we will distinguish these two types of projection by explicitly referencing their equation numbers throughout the remainder of this article.
\end{remark}
\begin{remark}\label{rem:P-BDF1}
\textcolor{black}{For the P-BDF1 scheme \eqref{se3-1}-\eqref{se3-3}, given \( (\phi_h^{n}, w_h^{n}, U^n, \mathbf{u}_h^n, p_h^n) \), we first compute \( U_h^n \) from the projection step in equation \eqref{se3-1 projection}. 
$\mathbf{\hat u}^{n + 1} \in (L^2(\Omega))^d$ in \eqref{se3-1a} is an intermediate function containing a stabilization term, and itself is not computed in the algorithm.
\( (\phi_h^{n+1}, w_h^{n+1}) \) is obtained by solving the linear system formed by equations \eqref{se3-1a} and \eqref{se3-1b}, with the substitutions from equations \eqref{se3-1c} and \eqref{se3-1e}, respectively
\begin{align}
	&\left( {\dfrac{{\phi _h^{n + 1} - \phi _h^n}}{\tau },{\varphi _h}} \right) -\left( \mathbf{u}_h^n \phi_h^n, \nabla \varphi_h \right) + \tau \left( (\phi_h^n)^2 \nabla w_h^{n+1}, \nabla \varphi_h \right) + \gamma a\left( {w_h^{n + 1},{\varphi _h}} \right) = 0, \label{equation 3.5}\\
	&\left( w_h^{n+1}, \psi_h \right) - \lambda a( \phi_h^{n+1}, \psi_h ) 
	= \lambda \left( H(\phi_h^n) U_h^n, \psi_h \right)
	+ \frac{1}{2} \lambda \left( H(\phi_h^n)^2 (\phi_h^{n+1} - \phi_h^n), \psi_h \right),\label{equation 3.6}
\end{align}
where each term in \eqref{equation 3.5}-\eqref{equation 3.6} can be computed by using Gaussian quadrature, which is equivalent to the process of $L^2$ projection.
Then, we update $U^{n+1}$ from \eqref{se3-1e}.
For the Navier-Stokes part, the intermediate velocity \( \mathbf{\tilde u}_h^{n+1} \) is computed from \eqref{se3-2}, followed by the projection step \eqref{se3-3} to obtain \( (\mathbf{u}_h^{n+1}, p_h^{n+1}) \).
All resulting subproblems are linear systems, which can be solved by appropriate linear solvers.}
\end{remark}
\begin{remark}
\textcolor{black}{Although the original Cahn--Hilliard equation contains a fourth-order differential operator, we reduce its order by introducing an auxiliary variable \( w \). This reformulation converts the fourth-order equation into a system of two coupled second-order equations, thereby avoiding the challenges of directly discretizing fourth-order terms and allowing the use of \( C^0 \) finite elements.}
\end{remark}
\begin{remark}
\textcolor{black}{ 
	The proposed methods and their analysis are applicable to other boundary conditions, such as nonhomogeneous mixed Dirichlet-Neumann boundary conditions and the periodic boundary conditions \cite{HINTERMULLER2013810,ye_efficient_2022}. 
	For the nonhomogeneous Neumann boundary condition, assume that the phase field variables $\phi$ and $w$ satisfy $\frac{\partial \phi}{\partial \nu}=\phi_N,\;\frac{\partial w}{\partial \nu}=w_N$. Applying the integration-by-parts formula, we reformulate Eqs. \eqref{se3-1a}-\eqref{se3-1b} in the P-BDF1 scheme as follows
	\begin{align*}
		&\left( {\dfrac{{\phi _h^{n + 1} - \phi _h^n}}{\tau },{\varphi _h}} \right) -\left( \mathbf{u}_h^n \phi_h^n, \nabla \varphi_h \right) + \tau \left( (\phi_h^n)^2 \nabla w_h^{n+1}, \nabla \varphi_h \right) + \gamma a\left( {w_h^{n + 1},{\varphi _h}} \right) -\gamma (w_N^{n+1},\varphi_h )_{\partial\Omega}= 0, \\
		&\left( w_h^{n+1}, \psi_h \right) - \lambda a( \phi_h^{n+1}, \psi_h ) + \lambda (\phi_N^{n+1},\psi_h )_{\partial\Omega}
		= \lambda \left( H(\phi_h^n) U_h^n, \psi_h \right)
		+ \frac{1}{2} \lambda \left( H(\phi_h^n)^2 (\phi_h^{n+1} - \phi_h^n), \psi_h \right),
	\end{align*}
	where $(p,q )_{\partial\Omega} = \int_{\partial\Omega}pq ds$. 
	Similarly, the proposed methods can be adapted to handle nonhomogeneous Dirichlet boundary conditions and periodic boundary conditions by modifying the finite element spaces and enforcing the prescribed boundary values appropriately. It is worth noting that the types of the boundary condition does not affect the ODE \eqref{uode} or its discrete forms, as they depend only on the initial data.
}
\end{remark}

\begin{lem}\label{lem3e1}\cite{chen_error_2022}
The bilinear form $b\left({\mathbf{u}}_h^{n},\cdot,\cdot\right)$ in the fully discrete scheme \eqref{se3-1}-\eqref{se3-3} is skew symmetric. Especially, it holds that
\begin{equation}\label{se3-4}
b\left( {\mathbf{u}}_h^n,{\mathbf{\tilde u}}_h^{n + 1},{\mathbf{\tilde u}}_h^{n + 1}\right)=0,\qquad \forall{\mathbf{\tilde u}}_h^{n + 1}\in {{\mathbf{X}}_h}.
\end{equation}
\end{lem}

\begin{lem}\label{lem3e2}
The solution of \eqref{se3-1}-\eqref{se3-3} conserves the total mass
\begin{equation}
\int_{\Omega} \phi_h^{n+1} d\mathbf{x} = \int_{\Omega} \phi_h^{0} d\mathbf{x}.
\end{equation}
\end{lem}
\begin{proof}
Taking $\varphi_h = 1$ in \eqref{se3-1a} yields the conservation of total mass.
\end{proof}

For scheme \eqref{se3-1}-\eqref{se3-3}, we establish the following result.
\begin{theorem}
The \textcolor{black}{P-BDF1 scheme} \eqref{se3-1}-\eqref{se3-3} admits a unique solution set 
$$(\phi_h^{n+1}, w_h^{n+1}, \mathbf{u}_h^{n+1}, p_h^{n+1}) \in Y_h \times Y_h \times \mathbf{V}_h\times M_h,$$ 
and $U^{n+1} \in C^0(\Omega)$.
\end{theorem}
\begin{proof}
\textcolor{black}{Since equations \eqref{se3-1}–\eqref{se3-3} constitute a finite-dimensional linear system, the existence of a solution is equivalent to its uniqueness. Now, assume that the linear system \eqref{se3-1}–\eqref{se3-3} has two solutions, and denote their difference
at $t
^{n+1}$ by 
$(\delta\phi_{h}^{n+1}, \delta w_{h}^{n+1}, \delta{\mathbf{\hat u}}_{h}^{n + 1}, \delta{{\mathbf{\tilde u}}_{h}^{n + 1}, \delta{\mathbf{u}}_{h}^{n+1}},\delta p_{h}^{n+1})$}, then it follows 
\begin{subequations}\label{dse3-1}
\begin{align}
\frac{1}{\tau}\left( {\delta \phi _h^{n + 1},{\varphi _h}} \right) - \left( { {{\delta 
	\mathbf{\hat u}}_h^{n + 1}\phi _h^n} ,\nabla{\varphi _h}} \right) + \gamma a\left( {\delta w_h^{n + 1},{\varphi _h}} \right) = 0, \qquad &\forall {\varphi _h} \in {Y_h}, \label{dse3-1a}\\
\left( {\delta w_h^{n + 1},{\psi _h}} \right) - \lambda a\left( {\delta \phi _h^{n + 1},{\psi _h}} \right)- \frac{\lambda}{2} \left( {H^2\left( {\phi _h^n} \right)\delta \phi_h^{n+1},{\psi _h}} \right) = 0, \qquad & {\forall {\psi _h} \in {Y_h}},\label{dse3-1b}\\
{\delta \mathbf{\hat u}}^{n + 1} = - \tau  {\phi _h^n\nabla \delta w_h^{n + 1}}.& \label{dse3-1c}
\end{align}
\end{subequations}
Taking $\varphi_h = \tau \delta w_h^{n+1},\; \psi_h = \delta \phi_h^{n+1}$ in \eqref{dse3-1a} and \eqref{dse3-1b}, respectively, multiplying both sides of equation \eqref{dse3-1c} by 
$\delta \mathbf{\hat u}^{n + 1}$, and summing them up yields
\begin{equation}
\tau\gamma \vert \delta w_h^{n+1} \vert^2_{H^1(\Omega)} + \lambda \vert \delta \phi_h^{n+1} \vert^2_{H^1(\Omega)} + \frac{\lambda}{2} \left\| H( {\phi _h^n} )\delta \phi_h^{n+1} \right\| ^2 + \left\| \delta \mathbf{\hat u}^{n + 1} \right\|^2 = 0,
\end{equation}
which implies $\delta \mathbf{\hat u}^{n + 1} = 0$ and $ \delta w_h^{n+1} = \delta \phi_h^{n+1} = \text{Constant}$. Then, plugging them into \eqref{dse3-1a} and \eqref{dse3-1b} leads to $\left( {\delta \phi _h^{n + 1},{\varphi _h}} \right) = \left( {\delta w_h^{n + 1},{\psi _h}} \right) = 0,\; \forall \varphi_h, \; \psi_h \in Y_h$, which implies $ \delta w_h^{n+1} = \delta \phi_h^{n+1} = 0$. 

We note that $U^{n+1}$ is determined by the known variables $\phi_h^n,\;U_h^n$ and the unique solution $\phi_h^{n+1}$, and is therefore unique. $\delta{\mathbf{\tilde u}}_{h}^{n + 1}$ satisfies
\begin{equation}\label{dse3-2}
\begin{aligned}
\frac{1}{\tau}\left( \delta{\mathbf{\tilde u}}_{h}^{n + 1},{{\mathbf{v}}_h} \right) + \mu \tilde a\left( \delta{\mathbf{\tilde u}}_{h}^{n + 1},{{\mathbf{v}}_h} \right) + b\left( {\mathbf{u}}_h^n,{\delta{\mathbf{\tilde u}}_{h}^{n + 1},{{\mathbf{v}}_h}} \right) 
= 0,& \qquad\forall {{\mathbf{v}}_h} \in {{\mathbf{X}}_h}.
\end{aligned}
\end{equation}
Taking $\mathbf{v}_h = \delta{\mathbf{\tilde u}}_{h}^{n + 1}$ and using Lemma \ref{lem3e1}, we obtain $$\left\| \delta{\mathbf{\tilde u}}_{h}^{n + 1} \right\|^2 = 0,$$ namely, $$\delta{\mathbf{\tilde u}}_{h}^{n + 1} = 0.$$

Finally, the variable $\delta{\mathbf{u}}_{h}^{n + 1}$ satisfies the following equations:
\begin{equation}\label{dse3-3}
\begin{aligned}
\left( \delta{\mathbf{u}}_{h}^{n + 1},{{\mathbf{\chi }}_h} \right) - \left(   (\delta p_{h}^{n + 1})  ,{\nabla \cdot {\mathbf{\chi }}_h} \right) = 0,\qquad&\forall {{\mathbf{\chi }}_h} \in {{\mathbf{V}}_h},\\
\left( {\nabla  \cdot (\delta{\mathbf{u}}_{h}^{n + 1})}, q_h \right) = 0,\qquad&\forall {q_h} \in {M_h}.
\end{aligned}
\end{equation}
Taking $\chi_h = \delta{\mathbf{u}}_{h}^{n + 1}$ and $q_h = \delta p_{h}^{n + 1} $ in \eqref{dse3-3} yield $\left\| \delta{\mathbf{\tilde u}}_{h}^{n + 1} \right\|^2 = 0$, namely, $\delta{\mathbf{\tilde u}}_{h}^{n + 1} = 0$, 
which together with the first equation in \eqref{dse3-3} and inf-sup condition \eqref{infsup} further implies $\nabla ( \delta p_{h}^{n + 1}) = 0$. 
\end{proof}
\begin{theorem}\label{thm1e1}
For the CHNS equations \eqref{eq1.3}, the \textcolor{black}{P-BDF1} scheme \eqref{se3-1}-\eqref{se3-3} is unconditionally energy stable and satisfies the following modified discrete energy law:
\begin{equation}\label{se3-6-1}
\begin{aligned}
E\left( \phi_h^{n+1},{\mathbf{u}}_h^{n+1},U_h^{n+1},p_h^{n+1} \right)&\leq E\left( \phi_h^{n+1},{\mathbf{u}}_h^{n+1},U^{n+1} ,p_h^{n+1}\right) \\
&=E\left( \phi_h^{n},{\mathbf{u}}_h^{n},U_h^{n},p_h^{n} \right)-\tau \gamma {\left\| {\nabla w_h^{n + 1}} \right\|^2} \\
&\quad- \frac{\lambda }{2} {{\left\| {\nabla \left( {\phi _h^{n + 1} - \phi _h^n} \right)} \right\|}^2}  - \lambda   {{\left\| {U^{n + 1} - U_h^n} \right\|}^2} \\
&\quad-
\frac{1}{2} {{\left\| {{\mathbf{\hat u}}^{n + 1} - {\mathbf{u}}_h^n} \right\|}^2} - \frac{1}{2}{{\left\| {{\mathbf{\tilde u}}_h^{n + 1} - {\mathbf{\hat u}}^{n + 1}} \right\|}^2}  - \tau \mu \left\| {\nabla {\mathbf{\tilde u}}_h^{n + 1}} \right\|^2,
\end{aligned}
\end{equation}
where
$$
E\left( \phi_h^{n+1},{\mathbf{u}}_h^{n+1},U_h^{n+1},p_h^{n+1} \right)=\mathcal{E}\left( \phi_h^{n+1},{\mathbf{u}}_h^{n+1},U_h^{n+1} \right) + \frac{\tau^2}{2}\left\| \nabla p _h^{n+1} \right\|,
$$
\textcolor{black}{and $\mathcal{E}\left( \phi_h^{n+1},{\mathbf{u}}_h^{n+1},U^{n+1} \right)$ is defined in equation \eqref{se2-1}}.
\end{theorem}

\begin{proof}
Step 1: Taking ${\varphi _h} = \tau w_h^{n + 1}\in {Y_h}$ in (\ref{se3-1}a),  \; ${\psi _h} =  - \left( {\phi _h^{n + 1} - \phi _h^n} \right) \in {Y_h}$ in (\ref{se3-1}b),\; ${{\mathbf{v}}_h} = \tau {\mathbf{\tilde u}}_h^{n + 1} \in {{\mathbf{X}}_h}$ in (\ref{se3-2}) give
\begin{equation}\label{se3-7}
\begin{aligned}
\left( {\phi _h^{n + 1} - \phi _h^n,w_h^{n + 1}} \right) + \tau \left( {\nabla  \cdot \left( {{\mathbf{\hat u}}^{n + 1}\phi _h^n} \right),w_h^{n + 1}} \right) + \tau \gamma a\left( {w_h^{n + 1},w_h^{n + 1}} \right) = 0 ,&\\
- \left( {w_h^{n + 1},\phi _h^{n + 1} - \phi _h^n} \right) + \lambda a\left( {\phi _h^{n + 1},\phi _h^{n + 1} - \phi _h^n} \right) + \lambda \left( {H\left( {\phi _h^n} \right)U^{n + 1},\phi _h^{n + 1} - \phi _h^n} \right) = 0 ,&\\
\left( {{\mathbf{\tilde u}}_h^{n + 1} - {\mathbf{u}}_h^n,{\mathbf{\tilde u}}_h^{n + 1}} \right) + \tau \mu \tilde a\left( {{\mathbf{\tilde u}}_h^{n + 1},{\mathbf{\tilde u}}_h^{n + 1}} \right) - \tau \left( {p_h^n,\nabla  \cdot {\mathbf{\tilde u}}_h^{n + 1}} \right) + \tau \left( {\phi _h^n\nabla w_h^{n + 1},{\mathbf{\tilde u}}_h^{n + 1}} \right) = 0.&
\end{aligned}
\end{equation}
The summation of \eqref{se3-7} yields
\begin{equation}\label{se3-8}
\begin{aligned}
\tau \left( {\nabla  \cdot \left( {{\mathbf{\hat u}}^{n + 1}\phi _h^n} \right),w_h^{n + 1}} \right) + \tau \gamma {\left\| {\nabla w_h^{n + 1}} \right\|^2}
+ \lambda a\left( {\phi _h^{n + 1},\phi _h^{n + 1} - \phi _h^n} \right) \\+ \lambda \left( {H\left( {\phi _h^n} \right)U^{n + 1},\phi _h^{n + 1} - \phi _h^n} \right)
+ \left( {{\mathbf{\tilde u}}_h^{n + 1} - {\mathbf{u}}_h^n,{\mathbf{\tilde u}}_h^{n + 1}} \right) \\+ \tau \mu \left\| {\nabla {\mathbf{\tilde u}}_h^{n + 1}} \right\|^2 - \tau \left( {p_h^n,\nabla  \cdot {\mathbf{\tilde u}}_h^{n + 1}} \right) + \tau \left( {\phi _h^n\nabla w_h^{n + 1},{\mathbf{\tilde u}}_h^{n + 1}} \right) = 0.
\end{aligned}
\end{equation}
Upon regrouping, it follows
\begin{equation}
\begin{array}{l}
- \tau \left( {{\mathbf{\hat u}}^{n + 1}\phi _h^n,\nabla w_h^{n + 1}} \right) + \tau \gamma {\left\| {\nabla w_h^{n + 1}} \right\|^2}
+ \frac{\lambda }{2}\left( {{{\left\| {\nabla \phi _h^{n + 1}} \right\|}^2} + {{\left\| {\nabla \left( {\phi _h^{n + 1} - \phi _h^n} \right)} \right\|}^2} - {{\left\| {\nabla \phi _h^n} \right\|}^2}} \right) \\
\qquad\qquad\qquad \qquad\quad\ + \lambda \left( {{{\left\| {U^{n + 1}} \right\|}^2} + {{\left\| {U^{n + 1} - U_h^n} \right\|}^2} - {{\left\| {U_h^n} \right\|}^2}} \right)
+ \left( {{\mathbf{\tilde u}}_h^{n + 1} - {\mathbf{u}}_h^n,{\mathbf{\tilde u}}_h^{n + 1}} \right) \\ \qquad\qquad\qquad \qquad\quad\ + \tau \mu \left\| {\nabla {\mathbf{\tilde u}}_h^{n + 1}} \right\|^2 - \tau \left( {p_h^n,\nabla  \cdot {\mathbf{\tilde u}}_h^{n + 1}} \right) + \tau \left( {\phi _h^n\nabla w_h^{n + 1},{\mathbf{\tilde u}}_h^{n + 1}} \right) = 0 .
\end{array} \label{se3-9}
\end{equation}

Step 2: Multiplying (\ref{se3-1}c) by ${\mathbf{\tilde u}}_h^{n + 1}$ and ${\mathbf{\hat u}}^{n + 1}$, respectively, gives
\begin{equation}
\left( {{\mathbf{\hat u}}^{n + 1},{\mathbf{\tilde u}}_h^{n + 1}} \right) = \left( {{\mathbf{u}}_h^n,{\mathbf{\tilde u}}_h^{n + 1}} \right) - \tau \left( {\phi _h^n\nabla w_h^{n + 1},{\mathbf{\tilde u}}_h^{n + 1}} \right) \label{se3-10},
\end{equation}
\begin{equation}
\left( {{\mathbf{\hat u}}^{n + 1},{\mathbf{\hat u}}^{n + 1}} \right) = \left( {{\mathbf{u}}_h^n,{\mathbf{\hat u}}^{n + 1}} \right) - \tau \left( {\phi _h^n\nabla w_h^{n + 1},{\mathbf{\hat u}}^{n + 1}} \right) \label{se3-11}.
\end{equation}

Plugging \eqref{se3-10} and \eqref{se3-11} into \eqref{se3-9} gives
\begin{equation}
\begin{aligned}
\left( {{\mathbf{\hat u}}^{n + 1} - {\mathbf{u}}_h^n,{\mathbf{\hat u}}^{n + 1}} \right) + \tau \gamma {\left\| {\nabla w_h^{n + 1}} \right\|^2}
+ \frac{\lambda }{2}\left( {{{\left\| {\nabla \phi _h^{n + 1}} \right\|}^2} + {{\left\| {\nabla \left( {\phi _h^{n + 1} - \phi _h^n} \right)} \right\|}^2} - {{\left\| {\nabla \phi _h^n} \right\|}^2}} \right) \\
+ \lambda \left( {{{\left\| {U^{n + 1}} \right\|}^2} + {{\left\| {U^{n + 1} - U_h^n} \right\|}^2} - {{\left\| {U_h^n} \right\|}^2}} \right)
+ \left( {{\mathbf{\tilde u}}_h^{n + 1} - {\mathbf{u}}_h^n,{\mathbf{\tilde u}}_h^{n + 1}} \right)\\ + \tau \mu \left\| {\nabla {\mathbf{\tilde u}}_h^{n + 1}} \right\|^2 - \tau \left( {p_h^n,\nabla  \cdot {\mathbf{\tilde u}}_h^{n + 1}} \right) - \left( {{\mathbf{\hat u}}^{n + 1} - {\mathbf{u}}_h^n,{\mathbf{\tilde u}}_h^{n + 1}} \right) = 0.
\end{aligned} \label{se3-12}
\end{equation}
Upon simplification, it follows
\begin{equation}
\begin{aligned}
\frac{1}{2}\left( {{{\left\| {{\mathbf{\hat u}}^{n + 1}} \right\|}^2} + {{\left\| {{\mathbf{\hat u}}^{n + 1} - {\mathbf{u}}_h^n} \right\|}^2} - {{\left\| {{\mathbf{u}}_h^n} \right\|}^2}} \right)
+ \frac{\lambda }{2}\left( {{{\left\| {\nabla \phi _h^{n + 1}} \right\|}^2} + {{\left\| {\nabla \left( {\phi _h^{n + 1} - \phi _h^n} \right)} \right\|}^2} - {{\left\| {\nabla \phi _h^n} \right\|}^2}} \right) \\
+ \tau \gamma {\left\| {\nabla w_h^{n + 1}} \right\|^2} + \lambda \left( {{{\left\| {U^{n + 1}} \right\|}^2} + {{\left\| {U^{n + 1} - U_h^n} \right\|}^2} - {{\left\| {U_h^n} \right\|}^2}} \right)\\
+ \frac{1}{2}\left( {{{\left\| {{\mathbf{\tilde u}}_h^{n + 1}} \right\|}^2} + {{\left\| {{\mathbf{\tilde u}}_h^{n + 1} - {\mathbf{\hat u}}^{n + 1}} \right\|}^2} - {{\left\| {{\mathbf{\hat u}}^{n + 1}} \right\|}^2}} \right) + \tau \mu \left\| {\nabla {\mathbf{\tilde u}}_h^{n + 1}} \right\|^2 - \tau \left( {p_h^n,\nabla  \cdot {\mathbf{\tilde u}}_h^{n + 1}} \right) = 0.
\end{aligned}\label{se3-13}
\end{equation}

Step 3: Let ${{\mathbf{\chi }}_h} = \tau \nabla p_h^n,\ {q_h} = p_h^n$ in \eqref{se3-3}. Then the summation gives
\begin{equation}
- \left( {{\mathbf{\tilde u}}_h^{n + 1},\nabla p_h^n} \right) + \frac{\tau }{2}\left( {{{\left\| {\nabla p_h^{n + 1}} \right\|}^2} - {{\left\| {\nabla \left( {p_h^{n + 1} - p_h^n} \right)} \right\|}^2} - {{\left\| {\nabla p_h^n} \right\|}^2}} \right) = 0.
\label{se3-14}
\end{equation}
By plugging \eqref{se3-14} into \eqref{se3-13} to replace $- \tau \left( {p_h^n,\nabla  \cdot {\mathbf{\tilde u}}_h^{n + 1}} \right)$, it holds
\begin{equation}
\begin{aligned}
\frac{1}{2}\left( {{{\left\| {{\mathbf{\hat u}}^{n + 1}} \right\|}^2} + {{\left\| {{\mathbf{\hat u}}^{n + 1} - {\mathbf{u}}_h^n} \right\|}^2} - {{\left\| {{\mathbf{u}}_h^n} \right\|}^2}} \right) 
+ \frac{\lambda }{2}\left( {{{\left\| {\nabla \phi _h^{n + 1}} \right\|}^2} + {{\left\| {\nabla \left( {\phi _h^{n + 1} - \phi _h^n} \right)} \right\|}^2} - {{\left\| {\nabla \phi _h^n} \right\|}^2}} \right) \\ + \tau \gamma {\left\| {\nabla w_h^{n + 1}} \right\|^2}
+ \lambda \left( {{{\left\| {U^{n + 1}} \right\|}^2}  + {{\left\| {U^{n + 1} - U_h^n} \right\|}^2} - {{\left\| {U_h^n} \right\|}^2}} \right)\\ 
+ \frac{1}{2}\left( {{{\left\| {{\mathbf{\tilde u}}_h^{n + 1}} \right\|}^2} + {{\left\| {{\mathbf{\tilde u}}_h^{n + 1} - {\mathbf{\hat u}}^{n + 1}} \right\|}^2} - {{\left\| {{\mathbf{\hat u}}^{n + 1}} \right\|}^2}} \right) + \tau \mu \left\| {\nabla {\mathbf{\tilde u}}_h^{n + 1}} \right\|^2\\
+ \frac{{{\tau ^2}}}{2}\left( {{{\left\| {\nabla p_h^{n + 1}} \right\|}^2} - {{\left\| {\nabla \left( {p_h^{n + 1} - p_h^n} \right)} \right\|}^2} - {{\left\| {\nabla p_h^n} \right\|}^2}} \right) = 0.
\end{aligned} \label{se3-15}
\end{equation}

Step 4: Further, taking ${{\mathbf{\chi }}_h}$ as $\tau \left( {{\mathbf{u}}_h^{n + 1} + {\mathbf{\tilde u}}_h^{n + 1}} \right)$ and $\tau \nabla \left( {p_h^{n + 1} - p_h^n} \right)$, respectively, and ${q_h}=p_h^{n + 1} - p_h^n$ in \eqref{se3-3}, it gives 
\begin{equation}\label{se3-16}
\begin{aligned}
\left( {{\mathbf{u}}_h^{n + 1} - {\mathbf{\tilde u}}_h^{n + 1},{\mathbf{u}}_h^{n + 1} + {\mathbf{\tilde u}}_h^{n + 1}} \right) + \tau \left( {\nabla \left( {p_h^{n + 1} - p_h^n} \right),{\mathbf{u}}_h^{n + 1} + {\mathbf{\tilde u}}_h^{n + 1}} \right) = 0,\\
\left( {{\mathbf{u}}_h^{n + 1} - {\mathbf{\tilde u}}_h^{n + 1},\nabla \left( {p_h^{n + 1} - p_h^n} \right)} \right) + \tau \left( {\nabla \left( {p_h^{n + 1} - p_h^n} \right),\nabla \left( {p_h^{n + 1} - p_h^n} \right)} \right) = 0 ,\\
\left( {{\mathbf{u}}_h^{n + 1},\nabla \left( {p_h^{n + 1} - p_h^n} \right)} \right) = 0.
\end{aligned}
\end{equation}
Summing up \eqref{se3-16} gives 
\begin{equation}
{\left\| {{\mathbf{\tilde u}}_h^{n + 1}} \right\|^2} = {\left\| {{\mathbf{u}}_h^{n + 1}} \right\|^2} + {\tau ^2}\left\| {\nabla \left( {p_h^{n + 1} - p_h^n} \right)}  \right\| ^2 .\label{se3-17}
\end{equation}
Finally, plugging \eqref{se3-17} into \eqref{se3-15} yields
\begin{equation}
\begin{aligned}
\frac{1}{2} {\left\| {{\mathbf{u}}_h^{n + 1}} \right\|}^2+ \frac{\lambda }{2} {\left\| {\nabla \phi _h^{n + 1}} \right\|}^2 +\lambda  {\left\| {U^{n + 1}} \right\|}^2 + \frac{{{\tau ^2}}}{2}{\left\| {\nabla p_h^{n + 1}} \right\|^2}\\  +\tau \gamma {\left\| {\nabla w_h^{n + 1}} \right\|^2} + \frac{\lambda }{2} {{\left\| {\nabla \left( {\phi _h^{n + 1} - \phi _h^n} \right)} \right\|}^2}  + \lambda  {{\left\| {U^{n + 1} - U_h^n} \right\|}^2}\\
+
\frac{1}{2} {{\left\| {{\mathbf{\hat u}}^{n + 1} - {\mathbf{u}}_h^n} \right\|}^2} + \frac{1}{2}  {{\left\| {{\mathbf{\tilde u}}_h^{n + 1} - {\mathbf{\hat u}}^{n + 1}} \right\|}^2} + \tau \mu \left\| {\nabla {\mathbf{\tilde u}}_h^{n + 1}} \right\|^2 \\
= \frac{1}{2}{\left\| {{\mathbf{u}}_h^n} \right\|^2} + \frac{\lambda }{2}{\left\| {\nabla \phi _h^n} \right\|^2} + \lambda {\left\| {U_h^n} \right\|^2} + \frac{{{\tau ^2}}}{2}{\left\| {\nabla p_h^n} \right\|^2} ,
\end{aligned}
\label{fnse3-18}
\end{equation}
which establishes the energy stability \eqref{se3-6-1}. 
\end{proof}

\subsubsection{\textcolor{black}{C-BDF1} scheme} 

An alternative first-order, fully discrete IEQ-FEM scheme—referred to as the C-BDF1 scheme—is formulated by retaining the intermediate function $U^n$ from the P-BDF1 scheme \eqref{se3-1}-\eqref{se3-3} in continuous function space $C^0(\Omega)$. Specifically, the scheme seeks functions $U^{n+1} \in C^0(\Omega)$, \textcolor{black}{$\mathbf{ \hat u^{n+1}} \in (L^2(\Omega))^d$}, and  $(\phi^{n+1}_h, w_h^{n+1}, {\mathbf{\tilde u}}_h^{n + 1}, p_h^{n+1}) \in Y_h  \times Y_h \times {\mathbf{X}}_h\times M_h$ such that 
\begin{subequations} \label{apdxA2.1}
\begin{align}
{\left( {\dfrac{{\phi _h^{n + 1} - \phi _h^n}}{\tau },{\varphi _h}} \right) - \left( { \left( {{\mathbf{\hat u}}^{n + 1}\phi _h^n} \right),{\nabla \varphi _h}} \right) + \gamma a\left( {w_h^{n + 1},{\varphi _h}} \right) = 0},&\qquad {\forall {\varphi _h} \in {Y_h}},\\
\left( {w_h^{n + 1},{\psi _h}} \right) - \lambda a\left( {\phi _h^{n + 1},{\psi _h}} \right) - \lambda \left( {H\left( {\phi _h^n} \right)U^{n+1},{\psi _h}} \right) = 0 ,&\qquad {\forall {\psi _h} \in {Y_h}},\\
{{\mathbf{\hat u}}^{n + 1} =  {{\mathbf{u}}_h^n} - \tau  {\phi _h^n\nabla w_h^{n + 1} } }, & \\
{U^{n + 1} = U^n + \frac{1}{2}H\left( {\phi _h^n} \right)\left( {\phi _h^{n + 1} - \phi _h^n} \right)}, &
\end{align}
\end{subequations}
\begin{equation} \label{apdxA2.2}
\begin{aligned}
\left( {\dfrac{{{\mathbf{\tilde u}}_h^{n + 1} - {\mathbf{u}}_h^n}}{\tau },{{\mathbf{v}}_h}} \right) + \mu \tilde a\left( {{\mathbf{\tilde u}}_h^{n + 1},{{\mathbf{v}}_h}} \right) + b\left( {{\mathbf{u}}_h^n,{\mathbf{\tilde u}}_h^{n + 1},{{\mathbf{v}}_h}} \right) 
- \left( {p_h^n,\nabla  \cdot {{\mathbf{v}}_h}} \right)& \\+ \left( {\phi _h^n\nabla w_h^{n + 1},{{\mathbf{v}}_h}} \right) = 0,& \qquad\forall {{\mathbf{v}}_h} \in {{\mathbf{X}}_h},
\end{aligned}
\end{equation}
and
\begin{equation} \label{apdxA2.3}
\begin{aligned}
\left( {\dfrac{{{\mathbf{u}}_h^{n + 1} - {\mathbf{\tilde u}}_h^{n + 1}}}{\tau },{{\mathbf{\chi }}_h}} \right) - \left( { \left( {p_h^{n + 1} - p_h^n} \right),{\nabla \cdot {\mathbf{\chi }}_h}} \right) = 0,\qquad&\forall {{\mathbf{\chi }}_h} \in {{\mathbf{V}}_h},\\
\left( {\nabla  \cdot {\mathbf{u}}_h^{n + 1},{q_h}} \right) = 0,\qquad&\forall {q_h} \in {M_h}.
\end{aligned}
\end{equation}
Similar to the \textcolor{black}{P-BDF1} scheme \eqref{se3-1}-\eqref{se3-3}, \textcolor{black}{the initial data set $(\phi^0_h,\mathbf{u}_h^0,U^0)$ is given by \eqref{phiuInitial} and \eqref{CHinit+}.} Additionally, we have the following result.

\begin{remark} \color{black}{The main difference between the C-BDF1 scheme \eqref{apdxA2.1}--\eqref{apdxA2.3} and the P-BDF1 scheme \eqref{se3-1}--\eqref{se3-3} lies in equations (\ref{apdxA2.1}d) and (\ref{se3-1}e). Specifically, the right-hand side of (\ref{apdxA2.1}d) does not enforce the previous-step auxiliary variable \( U^n \) to lie in the finite element space.} 

\begin{theorem}\label{lem-appendix1}
For the CHNS equations \eqref{eq1.3}, the \textcolor{black}{C-BDF1} scheme \eqref{apdxA2.1}-\eqref{apdxA2.3} is unconditionally energy stable and satisfies the following modified discrete energy law:
\begin{equation}\label{seappendix-1}
\begin{aligned}
E\left( \phi_h^{n+1},{\mathbf{u}}_h^{n+1},U^{n+1},p_h^{n+1} \right) 
&=E\left( \phi_h^{n},{\mathbf{u}}_h^{n},U^{n} ,p_h^{n}\right)-\tau \gamma {\left\| {\nabla w_h^{n + 1}} \right\|^2} - \frac{\lambda }{2} {{\left\| {\nabla \left( {\phi _h^{n + 1} - \phi _h^n} \right)} \right\|}^2}  \\
&\quad- \lambda   {{\left\| {U^{n + 1} - U^n} \right\|}^2} -
\frac{1}{2} {{\left\| {{\mathbf{\hat u}}^{n + 1} - {\mathbf{u}}_h^n} \right\|}^2} \\
&\quad- \frac{1}{2}{{\left\| {{\mathbf{\tilde u}}_h^{n + 1} - {\mathbf{\hat u}}^{n + 1}} \right\|}^2}  - \tau \mu \left\| {\nabla {\mathbf{\tilde u}}_h^{n + 1}} \right\|^2,
\end{aligned}
\end{equation}
where
$$E\left( \phi_h^{n+1},{\mathbf{u}}_h^{n+1},U^{n+1} ,p_h^{n+1}\right)=\mathcal{E}\left( \phi_h^{n+1},{\mathbf{u}}_h^{n+1},U^{n+1} \right)+\frac{\tau^2}{2}\left\| \nabla p _h^{n+1} \right\|,$$
\textcolor{black}{and $\mathcal{E}\left( \phi_h^{n+1},{\mathbf{u}}_h^{n+1},U^{n+1} \right)$ is defined in equation \eqref{se2-1}.}
\end{theorem}
The proof of \Cref{lem-appendix1} is similar to that of \Cref{thm1e1}.

\end{remark}
The \textcolor{black}{C-BDF1} scheme offers a computationally efficient alternative while still preserving the energy dissipation law~\eqref{seappendix-1}. However, the exact computation of the discrete energy 
\(
E(\phi_h^{n}, \mathbf{u}_h^{n}, U^{n}, p_h^n)
\)
is no longer feasible, since $U^{n} \notin Y_h$. Instead, the energy must be approximated by
\(
E(\phi_h^{n}, \mathbf{u}_h^{n}, U_{I}^{n}, p_h^n),
\)
where $U_{I}^{n} \in Y_h$ denotes the interpolation of $U^n$ evaluated at Gaussian quadrature points.

\begin{remark}
\color{black}{
At each time step, the following identity holds:
\begin{equation}
E\left( \phi_h^{n},{\mathbf{u}}_h^{n},U^{n} ,p_h^{n}\right) =  E\left( \phi_h^{n},{\mathbf{u}}_h^{n},U^{n}_I,p_h^{n} \right) + \lambda(\|U^{n}\|^2-\|U_I^{n}\|^2). \label{DES EU}
\end{equation}
Based on this relation \eqref{DES EU}, the discrete energy dissipation law \eqref{seappendix-1} can be equivalently expressed in terms of the approximated energy as:
\begin{equation}\label{seappendix-1+}
\begin{aligned}
& E\left( \phi_h^{n+1},{\mathbf{u}}_h^{n+1},U^{n+1}_I,p_h^{n+1} \right) 
-E\left( \phi_h^{n},{\mathbf{u}}_h^{n},U^{n}_I ,p_h^{n}\right) \\ 
& = -\tau \gamma {\left\| {\nabla w_h^{n + 1}} \right\|^2} - \frac{\lambda }{2} {{\left\| {\nabla \left( {\phi _h^{n + 1} - \phi _h^n} \right)} \right\|}^2} 
- \lambda   {{\left\| {U^{n + 1} - U^n} \right\|}^2} -
\frac{1}{2} {{\left\| {{\mathbf{\hat u}}^{n + 1} - {\mathbf{u}}_h^n} \right\|}^2} \\
& - \frac{1}{2}{{\left\| {{\mathbf{\tilde u}}_h^{n + 1} - {\mathbf{\hat u}}^{n + 1}} \right\|}^2}  - \tau \mu \left\| {\nabla {\mathbf{\tilde u}}_h^{n + 1}} \right\|^2 +\lambda(\|U_I^{n+1}\|^2-\|U_I^{n}\|^2-\|U^{n+1}\|^2+\|U^{n}\|^2).
\end{aligned}
\end{equation}
However, the decay of the approximated energy may not be guaranteed if the last term in \eqref{seappendix-1+} becomes positive and dominant.
}
\end{remark}

\begin{remark}\label{Enerdis} 
The error between the discrete energy $E(\phi_h^{n},{\mathbf{u}}_h^{n}, U^{n}, p_h^n)$ in \Cref{lem-appendix1} and the approximated energy $E(\phi_h^{n},{\mathbf{u}}_h^{n}, U_{I}^{n}, p_h^n)$ satisfies
\begin{equation}
\left| E\left( \phi_h^{n},{\mathbf{u}}_h^{n},U^{n} ,p_h^{n}\right) -  E\left( \phi_h^{n},{\mathbf{u}}_h^{n},U^{n}_I,p_h^{n} \right) \right| \leq  \lambda(\|U^{n}\|+\|U_I^{n}\|) (\|U^{n}-U_I^{n}\|) . \label{DES EU2}
\end{equation}
where the interpolation error $\|U^n-U_{I}^{n}\|$ may be reduced by refining meshes or adjusting the constant $B$ to an appropriate larger number \cite{LY22}.
\end{remark}
\subsubsection{\textcolor{black}{CP-BDF1} scheme}

The \textcolor{black}{C-BDF1} scheme is straightforward to implement and stands out as the most computationally efficient. However, as noted in \Cref{Enerdis}, its approximated energy in $V_h$ is only conditionally dissipative. On the other hand, the \textcolor{black}{P-BDF1} scheme includes an additional projection step compared to the \textcolor{black}{C-BDF1} scheme, but it ensures unconditional energy dissipation, which can be directly computed in $V_h$.
 
To maintain the advantages of both methods, we design a hybrid \textcolor{black}{CP-BDF1} scheme, which starts with the \textcolor{black}{C-BDF1} scheme and switches to the \textcolor{black}{P-BDF1} scheme at time $t_n=t _*$ if the approximated energy increases, i.e., $E(\phi_h^{n+1},{\mathbf{u}}_h^{n+1}, U_{I}^{n+1}, p_h^{n+1})>E(\phi_h^{n},{\mathbf{u}}_h^{n}, U_{I}^{n}, p_h^n)$. The detailed algorithm is presented in \Cref{CP algorithm}.

\begin{algorithm}[!ht]
\caption{The \textcolor{black}{CP-BDF1} scheme}
\begin{algorithmic}[1]\label{CP algorithm}
\STATE Start with time step $\tau$, total time $T$   and initial solution $(\phi_h^{0},\mathbf{u}_h^{0},p_h^0,U^{0})$;\;

\STATE Set $n=0$ and $t_0=0$;\;

\STATE Compute the energy $E_h^0 = E\left( \phi_h^{0},{\mathbf{u}}_h^{0},U_I^{0} ,p_h^{0}\right)$;\;
\FOR {$0\leq n \leq \frac{T}{\tau}$}

\STATE Set $t_{n+1}=t_n+\tau$;\;
\STATE Solve $(\phi_h^{n+1},w_h^{n+1},\mathbf{u}_h^{n+1},p_h^{n+1})$ and $U_I^{n+1}$ by using the \textcolor{black}{C-BDF1} scheme;\; \label{STATE Cscheme}
\STATE Compute the energy $E_h^{n+1} = E\left( \phi_h^{n+1},{\mathbf{u}}_h^{n+1},U_I^{n+1} ,p_h^{n+1}\right)$;\; \label{STATE Cscheme energy}
\IF{$E_h^{n+1}>E_h^n$}
\STATE Store $n^* = n$, $(\phi_h^{n},w_h^{n},\mathbf{u}_h^{n},p_h^{n})$ and $U_I^{n}$;
\STATE Break;
\ENDIF



\STATE $n=n+1$;\;
\ENDFOR
\FOR {$n^*\leq n\leq\frac{T}{\tau}$}
\STATE Set $t_{n+1}=t_n+\tau$;\;
\STATE Compute the $L^2$ projection $U_h^n$ of $U_I^n$ in the finite element space $Y_h$;\;
\STATE Solve $(\phi_h^{n+1},w_h^{n+1},\mathbf{u}_h^{n+1},p_h^{n+1})$ and $U_h^{n+1}$ by using the \textcolor{black}{P-BDF1} scheme;\; 
\STATE Compute the energy $E_h^{n+1} = E\left( \phi_h^{n+1},{\mathbf{u}}_h^{n+1},U_h^{n+1} ,p_h^{n+1}\right)$;\; 
\STATE $n=n+1$;\;
\ENDFOR

\end{algorithmic}
\end{algorithm}
\begin{theorem} \label{Thcpenergylaw}
For $n\geq 0$, let $E_h^{n}$ be the discrete energy of the solution from the hybrid \textcolor{black}{CP-BDF1} scheme at $t_n$. Then it holds
\begin{equation}\label{cpenergylaw}
E_h^{n+1} \leq E_h^{n},
\end{equation}
where the approximated energy
$$
E_h^{n}=E(\phi_h^{n},{\mathbf{u}}_h^{n}, U_{I}^{n}, p_h^{n})
$$ for $n \leq n^*$, and $$
E_h^{n}=E(\phi_h^{n},{\mathbf{u}}_h^{n}, U_h^{n}, p_h^{n})
$$ for $n > n^*$, where $n^*$ given in \Cref{CP algorithm}.
\end{theorem}
\begin{proof}
\Cref{CP algorithm} together with \Cref{lem-appendix1} and \Cref{Enerdis} implies that \eqref{cpenergylaw} holds for $n < n^*$, namely
$$ E(\phi_h^{n+1},{\mathbf{u}}_h^{n+1}, U_{I}^{n+1}, p_h^{n+1}) \leq E(\phi_h^{n},{\mathbf{u}}_h^{n}, U_{I}^{n}, p_h^n).$$
\Cref{thm1e1} implies that \eqref{cpenergylaw} also holds for $n>n^*$. To this end, we only need to prove that \eqref{cpenergylaw} holds for $n=n^*$.
Recall that 
\begin{align}
E_h^{n^*}= & E(\phi_h^{n^*},{\mathbf{u}}_h^{n^*}, U_{I}^{n^*}, p_h^{n^*}),\\
E_h^{n^*+1}= & E(\phi_h^{n^*+1},{\mathbf{u}}_h^{n^*+1}, U_h^{n^*+1}, p_h^{n^*+1}).
\end{align}
Also recall that $U_h^{n^*+1}$ is obtained in the following steps
\begin{align}
& \left( {U_h^{n^*},{\mu _h}} \right) = \left( {U^{n^*}_I,{\mu _h}} \right),\qquad {\forall {\mu _h} \in {Y_h}},\label{uiproj1}\\
& U^{n^* + 1} =  U_h^{n^*} + \frac{1}{2}H\left( {\phi _h^{n^*}} \right)\left( {\phi _h^{n^* + 1} - \phi _h^{n^*}} \right),\label{uiproj2}\\
& \left( {U_h^{n^*+1},{\mu _h}} \right) = \left( {U^{n^*+1},{\mu _h}} \right),\qquad {\forall {\mu _h} \in {Y_h}}.\label{uiproj3}
\end{align}
Here, \eqref{uiproj1} gives
\begin{align}\label{Enstar1}
E(\phi_h^{n^*},{\mathbf{u}}_h^{n^*}, U_{h}^{n^*}, p_h^{n^*}) \leq E(\phi_h^{n^*},{\mathbf{u}}_h^{n^*}, U_{I}^{n^*}, p_h^{n^*}).
\end{align}
Similar to \Cref{thm1e1}, \eqref{uiproj2} and \eqref{uiproj3} together with scheme \eqref{se3-1}-\eqref{se3-3} yield
\begin{align}\label{Enstar2}
E\left( \phi_h^{n^*+1},{\mathbf{u}}_h^{n^*+1},U_h^{n^*+1},p_h^{n^*+1} \right)&\leq E\left( \phi_h^{n^*+1},{\mathbf{u}}_h^{n^*+1},U^{n^*+1} ,p_h^{n^*+1}\right) \leq E\left( \phi_h^{n^*},{\mathbf{u}}_h^{n^*},U_h^{n^*},p_h^{n^*} \right).
\end{align}
Finally, \eqref{Enstar1} and \eqref{Enstar2} imply that \eqref{cpenergylaw} holds $n=n^*$.
\end{proof}

\begin{remark}
The \textcolor{black}{CP-BDF1} scheme in \Cref{CP algorithm} combines the computational efficiency of the \textcolor{black}{C-BDF1} scheme with the unconditional energy stability of the \textcolor{black}{P-BDF1} scheme.
\end{remark}


\subsection{Second order fully discrete schemes}
\subsubsection{\textcolor{black}{P-BDF2} scheme} In this subsection, we present the second order fully discrete IEQ-FEM scheme with intermediate function in polynomial space (\textcolor{black}{P-BDF2} scheme) for the CHNS equations.

Motivated by the techniques for the \textcolor{black}{P-BDF1} scheme, we consider the following second order fully discrete backward time-differentiation formula (BDF2) for the CHNS equations by approximating the intermediate function $U^{n+1}\in C^{0}(\Omega)$ and functions $(\phi^{n+1}_h, w_h^{n+1}, {{\mathbf{\tilde u}}_h^{n + 1}, {\mathbf{u}}_h^{n+1}},p^{n+1}_h) \in Y_h  \times Y_h\times {\mathbf{X}}_h\times {\mathbf{V}}_h\times M_h $, such that
\begin{subequations}\label{se3-18}
\begin{align}
\left( {\dfrac{{3\phi _h^{n + 1} - 4\phi _h^n + \phi _h^{n - 1}}}{{2\tau }},{\varphi _h}} \right) - \left( { \left( {{\mathbf{\tilde u}}_h^{n + 1}\phi _h^{n,\star}} \right),{\nabla\varphi _h}} \right)&\nonumber\\ + \gamma a\left( {w_h^{n + 1},{\varphi _h}} \right) = 0,\qquad\forall {\varphi _h} \in {Y_h},\\
\left( {w_h^{n + 1},{\psi _h}} \right) - \lambda a\left( {\phi _h^{n + 1},{\psi _h}} \right) - \lambda\left( H{{\left( {\phi _h^{n,\star}} \right)}U ^{n + 1}},{\psi _h}\right)=0,&\qquad\forall {\psi _h} \in {Y_h},\\
\left( {\dfrac{{3{\mathbf{\tilde u}}_h^{n + 1} - 4{\mathbf{u}}_h^n + {\mathbf{u}}_h^{n - 1}}}{{2\tau }},{{\mathbf{v}}_h}} \right) + \mu \tilde a\left( {{\mathbf{\tilde u}}_h^{n + 1},{{\mathbf{v}}_h}} \right) + b\left( {{\mathbf{u}}_h^{n,\star},{\mathbf{\tilde u}}_h^{n + 1},{{\mathbf{v}}_h}} \right)\nonumber&\\
- \left( {p_h^n,\nabla  \cdot {{\mathbf{v}}_h}} \right) + \left( {\phi _h^{n,\star}\nabla w_h^{n + 1},{{\mathbf{v}}_h}} \right) = 0
,&\qquad\forall {{\mathbf{v}}_h} \in {{\mathbf{X}}_h},\\
\left( {U_h^n,{\mu _h}} \right) = \left( {U^n,{\mu _h}} \right),&\qquad  {\forall {\mu _h} \in {Y_h}}, \label{se3-18 projection}\\
U^{n + 1} = \dfrac{{4U_h^n - U_h^{n - 1}}}{3} + \dfrac{1}{2}H\left( {\phi _h^{n,\star}} \right)\left( {\dfrac{{3\phi _h^{n + 1} - 4\phi _h^n + \phi _h^{n - 1}}}{3}} \right),&
\end{align}
\end{subequations}
and
\begin{equation}\label{se3-19}
\begin{aligned}
\left( {\dfrac{{3{\mathbf{u}}_h^{n + 1} - 3{\mathbf{\tilde u}}_h^{n + 1}}}{{2\tau }},{{\mathbf{\chi }}_h}} \right) - \left( { \left( {p_h^{n + 1} - p_h^n} \right),{\nabla \cdot {\mathbf{\chi }}_h}} \right) = 0,\qquad&\forall {{\mathbf{\chi }}_h} \in {{\mathbf{V}}_h},\\
\left( {\nabla  \cdot {\mathbf{u}}_h^{n + 1},{q_h}} \right) = 0,\qquad&\forall {q_h} \in {M_h},
\end{aligned}
\end{equation}
where $v_h^{n,\star}$ is defined as 
\[v_h^{n,\star}=2v_h^{n}-v_h^{n-1}.\]
\textcolor{black}{The initial data set $(\phi^0_h,\mathbf{u}_h^0,U^0)$ is given by \eqref{phiuInitial} and \eqref{CHinit+}, and the first-step values $(\phi^1_h,\mathbf{u}_h^1,U^1)$ are obtained by using the P-BDF1 scheme \eqref{se3-1}-\eqref{se3-3}.}
\begin{lem}\label{lem2e2} \cite{liu_unconditionally_2021}
For any symmetric bilinear form $A(\cdot,\cdot)$, it satisfies 
\begin{equation}\label{se3-20}
\begin{aligned}
A(\phi+\psi, \phi-\psi) = & A(\phi,\phi) - A(\psi,\psi),\\
2A\left( {3{\phi _1} - 2{\phi _2} - {\phi _3},{\phi _1}} \right) = &A\left( {{\phi _1},{\phi _1}} \right) + A\left( {2{\phi _1} - {\phi _2},2{\phi _1} - {\phi _2}} \right) - A\left( {{\phi _2},{\phi _2}} \right)\\
&+ A\left( {{\phi _1} - {\phi _3},{\phi _1} - {\phi _3}} \right) - A\left( {{\phi _3},{\phi _3}} \right) ,
\end{aligned} 
\end{equation}
where $\phi, \ \psi, \ \phi_1,\ \phi_2,\ \phi_3 \in \mathbf V_h$. 
\end{lem}

Next, similar to the \textcolor{black}{P-BDF1} scheme \eqref{se3-1}-\eqref{se3-3}, the following result holds for the \textcolor{black}{P-BDF2} scheme \eqref{se3-18}-\eqref{se3-19}. 

\begin{theorem}\label{thm1e2}
The \textcolor{black}{P-BDF2} scheme \eqref{se3-18}-\eqref{se3-19} satisfies the following energy dissipation law
\begin{equation}\label{se3-6-2}
\begin{aligned}
&\bar E\left( \phi_h^{n+1},\phi_h^{n+1,\star},{\mathbf{u}}_h^{n+1},{\mathbf{u}}_h^{n+1,\star},U_h^{n+1},U_h^{n+1,\star},p_h^{n+1} \right)\\ &\leq \bar E\left(  \phi_h^{n+1},\phi_h^{n+1,\star},{\mathbf{u}}_h^{n+1},{\mathbf{u}}_h^{n+1,\star},{2U^{n+1} - U_h^{n}},p_h^{n+1} \right)  \\
&=\bar E\left( \phi_h^{n},\phi_h^{n,\star},{\mathbf{u}}_h^{n},{\mathbf{u}}_h^{n,\star},U_h^{n},U_h^{n,\star},p_h^n \right)- 2\tau \mu \left\| {\nabla {\mathbf{\tilde u}}_h^{n + 1}} \right\|^2 - 2\tau \gamma {\left\| {\nabla w_h^{n + 1}} \right\|^2} \\&\quad- \frac{\lambda }{2}{\left\| {\nabla \left( {\phi _h^{n + 1} - \phi _h^{n,\star}} \right)} \right\|^2} - \lambda {\left\| {U^{n + 1} - U_h^{n,\star}} \right\|^2}
- \frac{1}{2}{\left\| {{\mathbf{u}}_h^{n + 1} - {\mathbf{u}}_h^{n,\star}} \right\|^2} \\&\quad- \frac{{2{\tau ^2}}}{3}{\left\| {\nabla \left( {p_h^{n + 1} - p_h^n} \right)} \right\|^2},
\end{aligned}
\end{equation}
where 
\begin{equation}
\begin{aligned}
&\bar E\left( \phi_h^{n},\phi_h^{n,\star},{\mathbf{u}}_h^{n},{\mathbf{u}}_h^{n,\star},U_h^{n},U_h^{n,\star},p_h^n \right)=\frac{\lambda }{2}\left( {{{\left\| {\nabla \phi _h^n} \right\|}^2} + {{\left\| {\nabla \phi _h^{n,\star}} \right\|}^2}} \right) \\&\qquad\qquad\qquad+ \lambda \left( {{{\left\| {U_h^n} \right\|}^2} + {{\left\| {U_h^{n,\star}} \right\|}^2}} \right) + \frac{1}{2}\left( {{{\left\| {{\mathbf{u}}_h^n} \right\|}^2} + {{\left\| {{\mathbf{u}}_h^{n,\star}} \right\|}^2}} \right) + \frac{{2{\tau ^2}}}{3}{\left\| {\nabla p_h^n} \right\|^2}.
\end{aligned}
\end{equation}
\end{theorem}

\begin{proof}
Firstly, let ${\varphi _h} = 2\tau w_h^{n + 1},\ {\psi _h} =  - \left( {3\phi _h^{n + 1} - 4\phi _h^n + \phi _h^{n - 1}} \right),\ {{\mathbf{v}}_h} = 2\tau {\mathbf{\tilde u}}_h^{n + 1}$ in (\ref{se3-18}a)-(\ref{se3-18}c), respectively. Then, it follows
\begin{equation}
\begin{aligned} 
\left( {3\phi _h^{n + 1} - 4\phi _h^n + \phi _h^{n - 1},w_h^{n + 1}} \right) + 2\tau \left( {\nabla  \cdot \left( {{\mathbf{\tilde u}}_h^{n + 1}\phi _h^{n,\star}} \right),w_h^{n + 1}} \right) + 2\tau \gamma a\left( {w_h^{n + 1},w_h^{n + 1}} \right) = 0&, \\
- \left( {w_h^{n + 1},3\phi _h^{n + 1} - 4\phi _h^n + \phi _h^{n - 1}} \right) + \lambda a\left( {\phi _h^{n + 1},3\phi _h^{n + 1} - 4\phi _h^n + \phi _h^{n - 1}} \right)\qquad\\
+ \lambda \left( {H\left( {\phi _h^{n,\star}} \right)U^{n + 1},3\phi _h^{n + 1} - 4\phi _h^n + \phi _h^{n - 1}} \right) = 0,&\\
\left( {3{\mathbf{\tilde u}}_h^{n + 1} - 4{\mathbf{u}}_h^n + {\mathbf{u}}_h^{n - 1},{\mathbf{\tilde u}}_h^{n + 1}} \right) + 2\tau \mu \tilde a\left( {{\mathbf{\tilde u}}_h^{n + 1},{\mathbf{\tilde u}}_h^{n + 1}} \right) - 2\tau \left( {p_h^n,\nabla  \cdot {\mathbf{\tilde u}}_h^{n + 1}} \right)\qquad\\
+ 2\tau \left( {\phi _h^{n,\star}\nabla w_h^{n + 1},{\mathbf{\tilde u}}_h^{n + 1}} \right) = 0.&
\end{aligned} \label{se3-22}
\end{equation}
The summation of \eqref{se3-22} with using (\ref{se3-18}e) upon regrouping gives
\begin{equation}
\begin{aligned}
2\tau \gamma {\left\| {\nabla w_h^{n + 1}} \right\|^2}
+ \lambda a\left( {\phi _h^{n + 1},3\phi _h^{n + 1} - 4\phi _h^n + \phi _h^{n - 1}} \right) \\
+ 2\lambda \left( {U^{n + 1},3U^{n + 1} - 4U_h^n + U_h^{n - 1}} \right)\\
+ \left( {3{\mathbf{u}}_h^{n + 1} - 4{\mathbf{u}}_h^n + {\mathbf{u}}_h^{n - 1},{\mathbf{\tilde u}}_h^{n + 1}} \right) + \left( {3{\mathbf{\tilde u}}_h^{n + 1} - 3{\mathbf{u}}_h^{n + 1},{\mathbf{\tilde u}}_h^{n + 1} + {\mathbf{u}}_h^{n + 1}} \right)\\
+ 2\tau \mu \left\| {\nabla {\mathbf{\tilde u}}_h^{n + 1}} \right\| - 2\tau \left( {p_h^n,\nabla  \cdot {\mathbf{\tilde u}}_h^{n + 1}} \right) = 0,
\end{aligned} \label{se3-23} 
\end{equation}
where we have used the identity $ \left( {3{\mathbf{\tilde u}}_h^{n + 1} - 3{\mathbf{u}}_h^{n + 1}, {\mathbf{u}}_h^{n + 1}} \right) = 0$, which follows by choosing $\chi_h = 2\tau \mathbf{u}_h^{n+1} $ and $q_h = 2\tau \nabla (p_h^{n+1}-p_h^n)$ in equation \eqref{se3-19}. Then according to Lemma \ref{lem2e2}, and 
\begin{equation}
\left( {3{\mathbf{u}}_h^{n + 1} - 4{\mathbf{u}}_h^n + {\mathbf{u}}_h^{n - 1},{\mathbf{\tilde u}}_h^{n + 1}} \right) = \left( {3{\mathbf{u}}_h^{n + 1} - 4{\mathbf{u}}_h^n + {\mathbf{u}}_h^{n - 1},{\mathbf{ u}}_h^{n + 1}} \right),
\end{equation}
which follows by choosing $\chi_h = 2\tau (3\mathbf{u}_h^{n+1} - 4\mathbf{u}_h^{n} + \mathbf{u}_h^{n-1}) $ in equation \eqref{se3-19}, it holds
\begin{equation}
\begin{aligned}
2\tau \gamma {\left\| {\nabla w_h^{n + 1}} \right\|^2} 
+ 2\tau \mu \left\| {\nabla {\mathbf{\tilde u}}_h^{n + 1}} \right\|^2 - 2\tau \left( {p_h^n,\nabla  \cdot {\mathbf{\tilde u}}_h^{n + 1}} \right)\\+\frac{\lambda }{2}\left( {{{\left\| {\nabla \phi _h^{n + 1}} \right\|}^2} 
	+ {{\left\| {\nabla \phi _h^{n + 1,\star}} \right\|}^2} - {{\left\| {\nabla \phi _h^n} \right\|}^2} - {{\left\| {\nabla \phi _h^{n,\star}} \right\|}^2} + {{\left\| {\nabla \left( {\phi _h^{n + 1} - \phi _h^{n,\star}} \right)} \right\|}^2}} \right)\\
+ \lambda \left( {{{\left\| {U^{n + 1}} \right\|}^2} + {{\left\| {{2U^{n+1}-U_h^n}} \right\|}^2} - {{\left\| {U_h^n} \right\|}^2} - {{\left\| {U_h^{n,\star}} \right\|}^2} + {{\left\| {U^{n + 1} - U_h^{n,\star}} \right\|}^2}} \right)\\
+ \frac{1}{2}\left( {{{\left\| {{\mathbf{u}}_h^{n + 1}} \right\|}^2} + {{\left\| {{\mathbf{u}}_h^{n + 1,\star}} \right\|}^2} - {{\left\| {{\mathbf{u}}_h^n} \right\|}^2} - {{\left\| {{\mathbf{u}}_h^{n,\star}} \right\|}^2} + {{\left\| {{\mathbf{u}}_h^{n + 1} - {\mathbf{u}}_h^{n,\star}} \right\|}^2}} \right)\\
+ 3\left( {{{\left\| {{\mathbf{\tilde u}}_h^{n + 1}} \right\|}^2} - {{\left\| {{\mathbf{u}}_h^{n + 1}} \right\|}^2}} \right)  = 0,
\end{aligned} \label{se3-24}
\end{equation}
here the properties of the numerical solution $\nabla \cdot \mathbf{u}^{m} = 0,\;m=1,2,\dots,n+1$ is applied. 
Secondly, taking ${{\mathbf{\chi }}_h},\  {q_h}$ as $ 2\tau \nabla p_h^n,\;  p_h^n$ in \eqref{se3-19}, respectively, we get
\begin{equation}
\begin{aligned}
\left( {3{\mathbf{u}}_h^{n + 1} - 3{\mathbf{\tilde u}}_h^{n + 1},\nabla p_h^n} \right) + 2\tau \left( {\nabla \left( {p_h^{n + 1} - p_h^n} \right),\nabla p_h^n} \right) = 0,\\
\left( {{\mathbf{u}}_h^{n + 1},\nabla p_h^n} \right) = 0,
\end{aligned} \label{se3-25}
\end{equation}
which can be simplified as
\begin{equation}
- \left( {3{\mathbf{\tilde u}}_h^{n + 1},\nabla p_h^n} \right) + \tau \left( {{{\left\| {\nabla p_h^{n + 1}} \right\|}^2} - {{\left\| {\nabla \left( {p_h^{n + 1} - p_h^n} \right)} \right\|}^2} - {{\left\| {\nabla p_h^n} \right\|}^2}} \right) = 0. \label{se3-26}
\end{equation}
Plugging \eqref{se3-26} into \eqref{se3-24} to replace  $- 2\tau \left( {p_h^n,\nabla  \cdot {\mathbf{\tilde u}}_h^{n + 1}} \right)$ yields
\begin{equation}
\begin{aligned}
\frac{\lambda }{2}\left( {{{\left\| {\nabla \phi _h^{n + 1}} \right\|}^2} + {{\left\| {\nabla \phi _h^{n + 1,\star}} \right\|}^2} - {{\left\| {\nabla \phi _h^n} \right\|}^2} - {{\left\| {\nabla \phi _h^{n,\star}} \right\|}^2} + {{\left\| {\nabla \left( {\phi _h^{n + 1} - \phi _h^{n,\star}} \right)} \right\|}^2}} \right)\\
+ \lambda \left( {{{\left\| {U^{n + 1}} \right\|}^2} + {{\left\| {{2U^{n+1}-U_h^n}} \right\|}^2} - {{\left\| {U_h^n} \right\|}^2} - {{\left\| {U_h^{n,\star}} \right\|}^2} + {{\left\| {U^{n + 1} - U_h^{n,\star}} \right\|}^2}} \right)\\
+ \frac{1}{2}\left( {{{\left\| {{\mathbf{u}}_h^{n + 1}} \right\|}^2} + {{\left\| {{\mathbf{u}}_h^{n + 1,\star}} \right\|}^2} - {{\left\| {{\mathbf{u}}_h^n} \right\|}^2} - {{\left\| {{\mathbf{u}}_h^{n,\star}} \right\|}^2} + {{\left\| {{\mathbf{u}}_h^{n + 1} - {\mathbf{u}}_h^{n,\star}} \right\|}^2}} \right)\\
+ 3\left( {{{\left\| {{\mathbf{\tilde u}}_h^{n + 1}} \right\|}^2} - {{\left\| {{\mathbf{u}}_h^{n + 1}} \right\|}^2}} \right) + 2\tau \mu \left\| {\nabla {\mathbf{\tilde u}}_h^{n + 1}} \right\|^2 + 2\tau \gamma {\left\| {\nabla w_h^{n + 1}} \right\|^2}\\
+ \frac{{2{\tau ^2}}}{3}\left( {{{\left\| {\nabla p_h^{n + 1}} \right\|}^2} - {{\left\| {\nabla \left( {p_h^{n + 1} - p_h^n} \right)} \right\|}^2} - {{\left\| {\nabla p_h^n} \right\|}^2}} \right) = 0.
\end{aligned} \label{se3-27}
\end{equation}

Finally, let ${{\mathbf{\chi }}_h} = 2\tau \left( {{\mathbf{u}}_h^{n + 1} + {\mathbf{\tilde u}}_h^{n + 1}} \right)$, ${q_h} = p_h^{n + 1} - p_h^n$ in \eqref{se3-19}. Then,
\begin{equation}
\begin{aligned}
\left( {3{\mathbf{u}}_h^{n + 1} - 3{\mathbf{\tilde u}}_h^{n + 1},{\mathbf{u}}_h^{n + 1} + {\mathbf{\tilde u}}_h^{n + 1}} \right) + 2\tau \left( {\nabla \left( {p_h^{n + 1} - p_h^n} \right),{\mathbf{u}}_h^{n + 1} + {\mathbf{\tilde u}}_h^{n + 1}} \right) = 0,\\
\left( {{\mathbf{u}}_h^{n + 1},\nabla \left( {p_h^{n + 1} - p_h^n} \right)} \right) = 0,
\end{aligned} \label{se3-28}
\end{equation}
leading to
\begin{equation}\label{se3-29}
3{\left\| {{\mathbf{u}}_h^{n + 1}} \right\|^2} - 3{\left\| {{\mathbf{\tilde u}}_h^{n + 1}} \right\|^2} + 2\tau \left( {\nabla \left( {p_h^{n + 1} - p_h^n} \right),{\mathbf{\tilde u}}_h^{n + 1}} \right) = 0. 
\end{equation}

In addition, taking ${{\mathbf{\chi }}_h} = \tau \nabla \left( {p_h^{n + 1} - p_h^n} \right)$ in \eqref{se3-19} gives
\begin{equation} \label{se3-30}
\left( {3{\mathbf{u}}_h^{n + 1} - 3{\mathbf{\tilde u}}_h^{n + 1},\nabla \left( {p_h^{n + 1} - p_h^n} \right)} \right) + 2\tau \left( {\nabla \left( {p_h^{n + 1} - p_h^n} \right),\nabla \left( {p_h^{n + 1} - p_h^n} \right)} \right) = 0.
\end{equation}
Upon simplification, it follows
\begin{equation}\label{se3-31}
\left( {{\mathbf{\tilde u}}_h^{n + 1},\nabla \left( {p_h^{n + 1} - p_h^n} \right)} \right) = \frac{2}{3}\tau \left( {\nabla \left( {p_h^{n + 1} - p_h^n} \right),\nabla \left( {p_h^{n + 1} - p_h^n} \right)} \right).
\end{equation}
Plugging \eqref{se3-31} into \eqref{se3-29} gives
\begin{equation}\label{se3-32}
{\left\| {{\mathbf{\tilde u}}_h^{n + 1}} \right\|^2} - {\left\| {{\mathbf{u}}_h^{n + 1}} \right\|^2} = \frac{4}{9}{\tau ^2}\left( {\nabla \left( {p_h^{n + 1} - p_h^n} \right),\nabla \left( {p_h^{n + 1} - p_h^n} \right)} \right) .
\end{equation}
By summing \eqref{se3-32} with  \eqref{se3-27}, 
\begin{equation}
\begin{aligned}
\frac{\lambda }{2}\left( {{{\left\| {\nabla \phi _h^{n + 1}} \right\|}^2} + {{\left\| {\nabla \phi _h^{n + 1,\star}} \right\|}^2}} \right) + \lambda \left( {{{\left\| {U^{n + 1}} \right\|}^2} + {{\left\| {{2U^{n+1}-U_h^n}} \right\|}^2}} \right) \\
+ \frac{1}{2}\left( {{{\left\| {{\mathbf{u}}_h^{n + 1}} \right\|}^2} + {{\left\| {{\mathbf{u}}_h^{n + 1,\star}} \right\|}^2}} \right) + \frac{{2{\tau ^2}}}{3}{\left\| {\nabla p_h^{n + 1}} \right\|^2}\\
+ 2\tau \mu \left\| {\nabla {\mathbf{\tilde u}}_h^{n + 1}} \right\|^2 + 2\tau \gamma {\left\| {\nabla w_h^{n + 1}} \right\|^2} + \frac{\lambda }{2}{\left\| {\nabla \left( {\phi _h^{n + 1} - \phi _h^{n,\star}} \right)} \right\|^2} \\+ \lambda {\left\| {U^{n + 1} - U_h^{n,\star}} \right\|^2}
+ \frac{1}{2}{\left\| {{\mathbf{u}}_h^{n + 1} - {\mathbf{u}}_h^{n,\star}} \right\|^2} + \frac{{2{\tau ^2}}}{3}{\left\| {\nabla \left( {p_h^{n + 1} - p_h^n} \right)} \right\|^2}\\
= \frac{\lambda }{2}\left( {{{\left\| {\nabla \phi _h^n} \right\|}^2} + {{\left\| {\nabla \phi _h^{n,\star}} \right\|}^2}} \right) + \lambda \left( {{{\left\| {U_h^n} \right\|}^2} + {{\left\| {U_h^{n,\star}} \right\|}^2}} \right) \\+ \frac{1}{2}\left( {{{\left\| {{\mathbf{u}}_h^n} \right\|}^2} + {{\left\| {{\mathbf{u}}_h^{n,\star}} \right\|}^2}} \right) + \frac{{2{\tau ^2}}}{3}{\left\| {\nabla p_h^n} \right\|^2}.
\end{aligned} \label{se3-33}
\end{equation}
{Based on the \textcolor{black}{stability} property of the $L^2$ projection, \textcolor{black}{$$\|U_h^{n+1}\| \le \|U^{n+1}\|,$$} and the following fact 
$$ (2U^{n+1}_h-U^{n}_h,\mu_h) = (2U^{n+1}-U^{n}_h,\mu_h), \qquad \forall \mu_h \in Y_h
$$}
implies the desired result \eqref{se3-6-2}.
\end{proof}

\subsubsection{\textcolor{black}{C-BDF2} scheme} 
We also process with another second-order
fully discrete IEQ-FEM scheme (\textcolor{black}{C-BDF2} scheme) for CHNS equations by approximating the intermediate function $U^{n+1}\in C^{0}(\Omega)$ and functions $(\phi^{n+1}_h, w_h^{n+1}, {{\mathbf{\tilde u}}_h^{n + 1}, {\mathbf{u}}_h^{n+1}}, p_h^{n+1}) \in Y_h  \times Y_h\times {\mathbf{X}}_h
\times {\mathbf{V}}_h\times M_h$, such that
\begin{subequations}\label{apdx3-18-1}
\begin{align}
\left( {\dfrac{{3\phi _h^{n + 1} - 4\phi _h^n + \phi _h^{n - 1}}}{{2\tau }},{\varphi _h}} \right) - \left( {\left( {{\mathbf{\tilde u}}_h^{n + 1}\phi _h^{n,\star}} \right),{\nabla \varphi _h}} \right)&\nonumber\\ + \gamma a\left( {w_h^{n + 1},{\varphi _h}} \right) = 0,\qquad\forall {\varphi _h} \in {Y_h},\\
\left( {w_h^{n + 1},{\psi _h}} \right) - \lambda a\left( {\phi _h^{n + 1},{\psi _h}} \right) - \lambda\left( H{{\left( {\phi _h^{n,\star}} \right)}U ^{n + 1}},{\psi _h}\right)=0,&\qquad\forall {\psi _h} \in {Y_h},\\
\left( {\dfrac{{3{\mathbf{\tilde u}}_h^{n + 1} - 4{\mathbf{u}}_h^n + {\mathbf{u}}_h^{n - 1}}}{{2\tau }},{{\mathbf{v}}_h}} \right) + \mu \tilde a\left( {{\mathbf{\tilde u}}_h^{n + 1},{{\mathbf{v}}_h}} \right) + b\left( {{\mathbf{u}}_h^{n,\star},{\mathbf{\tilde u}}_h^{n + 1},{{\mathbf{v}}_h}} \right)\nonumber&\\
- \left( {p_h^n,\nabla  \cdot {{\mathbf{v}}_h}} \right) + \left( {\phi _h^{n,\star}\nabla w_h^{n + 1},{{\mathbf{v}}_h}} \right) = 0
,&\qquad\forall {{\mathbf{v}}_h} \in {{\mathbf{X}}_h},\\
U^{n + 1} = \dfrac{{4U^n - U^{n - 1}}}{3} + \dfrac{1}{2}H\left( {\phi _h^{n,\star}} \right)\left( {\dfrac{{3\phi _h^{n + 1} - 4\phi _h^n + \phi _h^{n - 1}}}{3}} \right),&
\end{align}
\end{subequations}
and
\begin{equation}\label{apdx3-19-1}
\begin{aligned}
\left( {\dfrac{{3{\mathbf{u}}_h^{n + 1} - 3{\mathbf{\tilde u}}_h^{n + 1}}}{{2\tau }},{{\mathbf{\chi }}_h}} \right) - \left( { \left( {p_h^{n + 1} - p_h^n} \right),{\nabla \cdot {\mathbf{\chi }}_h}} \right) = 0,\qquad&\forall {{\mathbf{\chi }}_h} \in {{\mathbf{V}}_h},\\
\left( {\nabla  \cdot {\mathbf{u}}_h^{n + 1},{q_h}} \right) = 0,\qquad&\forall {q_h} \in {M_h}.
\end{aligned}
\end{equation}
\textcolor{black}{Similar to the P-BDF2 scheme, the initial data set $(\phi^0_h,\mathbf{u}_h^0,U^0)$ is given by \eqref{phiuInitial} and \eqref{CHinit+}, and the first-step values $(\phi^1_h,\mathbf{u}_h^1,U^1)$ are obtained using the C-BDF1 scheme \eqref{apdxA2.1}-\eqref{apdxA2.3}.} 

\textcolor{black}{\begin{remark}
The C-BDF2 scheme \eqref{apdx3-18-1}-\eqref{apdx3-19-1} first solves for $(\phi_h^{n+1},w_h^{n+1},\mathbf{\tilde u}_h^{n+1})$ by computing the coupled system. Compared with the C-BDF1 scheme \eqref{apdxA2.1}-\eqref{apdxA2.3}, the C-BDF2 requires solving three unknown variables simultaneously, whereas the C-BDF1 involves at most two, which increases the CPU time.
\end{remark}}

The following result holds for the \textcolor{black}{C-BDF2} scheme \eqref{apdx3-18-1}-\eqref{apdx3-19-1}. 
\begin{lem}\label{lem-appendix2}
The \textcolor{black}{C-BDF2} scheme \eqref{apdx3-18-1}-\eqref{apdx3-19-1} satisfies the following energy dissipation law
\begin{equation}\label{energy-appendix2}
\begin{aligned}
&\bar E\left( \phi_h^{n+1},\phi_h^{n+1,\star},{\mathbf{u}}_h^{n+1},{\mathbf{u}}_h^{n+1,\star},U^{n+1},U^{n+1,\star},p_h^{n+1} \right)\\ 
=&\bar E\left( \phi_h^{n},\phi_h^{n,\star},{\mathbf{u}}_h^{n},{\mathbf{u}}_h^{n,\star},U^{n},U^{n,\star},p_h^n \right)- 2\tau \mu \left\| {\nabla {\mathbf{\tilde u}}_h^{n + 1}} \right\|^2 - 2\tau \gamma {\left\| {\nabla w_h^{n + 1}} \right\|^2} \\&- \frac{\lambda }{2}{\left\| {\nabla \left( {\phi _h^{n + 1} - \phi _h^{n,\star}} \right)} \right\|^2} - \lambda {\left\| {U^{n + 1} - U^{n,\star}} \right\|^2}
- \frac{1}{2}{\left\| {{\mathbf{u}}_h^{n + 1} - {\mathbf{u}}_h^{n,\star}} \right\|^2} \\&- \frac{{2{\tau ^2}}}{3}{\left\| {\nabla \left( {p_h^{n + 1} - p_h^n} \right)} \right\|^2},
\end{aligned}
\end{equation}
where 
\begin{equation}
\begin{aligned}
&\bar E\left( \phi_h^{n},\phi_h^{n,\star},{\mathbf{u}}_h^{n},{\mathbf{u}}_h^{n,\star},U^{n},U^{n,\star},p_h^n \right)=\frac{\lambda }{2}\left( {{{\left\| {\nabla \phi _h^n} \right\|}^2} + {{\left\| {\nabla \phi _h^{n,\star}} \right\|}^2}} \right) \\&\qquad\qquad\qquad+ \lambda \left( {{{\left\| {U^n} \right\|}^2} + {{\left\| {U^{n,\star}} \right\|}^2}} \right) + \frac{1}{2}\left( {{{\left\| {{\mathbf{u}}_h^n} \right\|}^2} + {{\left\| {{\mathbf{u}}_h^{n,\star}} \right\|}^2}} \right) + \frac{{2{\tau ^2}}}{3}{\left\| {\nabla p_h^n} \right\|^2},
\end{aligned}
\end{equation}
with the term $v^{{n+1},\star}$ being defined as 
\[v^{{n+1},\star}=2v^{n+1}-v^{n}.\]
\end{lem}

\begin{remark}
Similar to the \textcolor{black}{CP-BDF1} scheme proposed in Algorithm 1, we can also present the \textcolor{black}{CP-BDF2} scheme by combining the \textcolor{black}{C-BDF2} scheme and the \textcolor{black}{P-BDF2} scheme. 
\end{remark}

\begin{remark}
The numerical solution from the second-order Crank–Nicolson (CN) IEQ-FEM shows instability, a phenomenon that has also been observed in the CN-IEQ-DG method \cite{liu_unconditionally_2021, yin2019efficient} and the CN-IEQ-FEM \cite{chen_IEQ_FEM} when solving the CH equation with logarithmic potential. Additionally, extending the current work to higher-order time discretizations \cite{KANGUnconditionally} is a potential direction for future research. The investigation into these time discretizations, including adaptive time discretizations \cite{YANGCompatibleL2}, in combination with the IEQ-FEMs, will be addressed in our future work.
\end{remark}

\section{Numerical examples}\label{se:4}

In this section, we present 2D and 3D numerical examples to validate the theoretical results presented in this paper. In implementation, we consider \textcolor{black}{uniform triangular meshes in 2D, and quasi-uniform tetrahedral meshes in 3D for the following numerical examples. 
The 2D triangular meshes are mainly generated as follows: the domain $\Omega = I_x \times I_y$ is first divided into an $N_x \times N_y$ rectangles, each of which is subdivided into two equal triangles. The mesh size of the triangulation is identified by either the partition number $N=N_x=N_y$ in each direction or the mesh size $h=h_x=h_y$, where $h_x= |I_x|/N_x$ and $h_y = |I_y|/N_y$. 3D tetrahedral meshes are generated and identified similarly.}

\subsection{Convergent rates}
In this part, we give an example  to validate the temporal and spatial convergence rates of the proposed numerical schemes for the CHNS equations. 

\begin{example}\label{exam1}
In the first example, we consider the following CHNS equations \begin{equation}\label{se4-1}
\begin{aligned}
{\partial _t}\phi  + \nabla  \cdot \left( {{\mathbf{u}}\phi } \right) - \gamma \Delta w = g(\textcolor{black}{t,x,y}),\qquad & {\it \text{in}}\ \Omega  \times \textcolor{black}{(0,T]},\\
w + \lambda \left( {\Delta \phi  - f\left( \phi  \right)} \right) = 0,\qquad & {\it \text{in}}\ \Omega  \times \textcolor{black}{(0,T]},\\
{\partial _t}{\mathbf{u}} - \mu \Delta {\mathbf{u}} + \left( {{\mathbf{u}} \cdot \nabla } \right){\mathbf{u}} + \nabla p + \phi \nabla w = {\mathbf{h}}(\textcolor{black}{t,x,y}),\qquad & {\it \text{in}}\ \Omega  \times \textcolor{black}{(0,T]},\\
\nabla  \cdot {\mathbf{u}} = 0,\qquad & {\it \text{in}}\ \Omega  \times \textcolor{black}{(0,T]},\\
{\mathbf{u}}\left( { \cdot ,0} \right) = {\mathbf{u}}_0,\ \phi \left( { \cdot ,0} \right) = {\phi _0},\qquad & {\it \text{in}}\ \Omega  \times \left\{ {t = 0} \right\},\\
{\mathbf{u}} = 0,\ \dfrac{{\partial \phi }}{{\partial {\textcolor{black}{\nu}}}} = 0,\ \dfrac{{\partial w}}{{\partial {\textcolor{black}{\nu}}}} = 0,\qquad & {\it \text{on}}\ \partial \Omega  \times \textcolor{black}{(0,T]},
\end{aligned}
\end{equation}
with $\Omega = \left[ 0 , 4\pi \right]^2$, the exact solution satisfies

\begin{equation}\label{se4-2}
\begin{aligned}
\phi \left( {t,x,y} \right) &= \sin \left( t \right)\cos \left( {\frac{x}{2}} \right)\cos \left( {\frac{y}{2}} \right),\\
\mathbf{u}\left( {t,x,y} \right) &= \left(
e^{ { - \frac{49t}{64}} }{{\sin }^2}\left( {\frac{x}{4}} \right)\sin \left( {\frac{y}{2}} \right),  - e^ { { - \frac{49t}{64}} }\sin \left( {\frac{x}{2}} \right){{\sin }^2}\left( \frac{y}{4} \right) \right)^\top,\\
p\left( {t,x,y} \right) &= \cos \left( \frac{x}{2} \right)\sin \left( \frac{y}{2} \right)\left( {1 - \sin \left( t \right)} \right), 
\end{aligned}
\end{equation}
the parameters $\epsilon = 1,\; \lambda = 1,\; \mu = 1,\; \gamma = 1,\; B=50$, and the corresponding right terms $g(x,t), \; {\mathbf{h}}(x,t)$ can be obtained by taking the the exact solution \eqref{se4-2} into \eqref{se4-1}. 
\end{example}
\begin{table}[H]
\begin{center}
\begin{tabular}{ c|c c c c c c c c }
\hline
& $\left\| \phi - \phi_h \right\| _ {L^2}$ & & $\left\| \mathbf{u} - \mathbf{u}_h \right\| _ {L^2}$ & & $\left| \phi - \phi_h \right| _ {H^1}$ & & $\left| \mathbf{u} - \mathbf{u}_h \right| _ {H^1}$ \\
\hline
$N = 4$ & 1.733e-05 &  $--$ & 2.742e-01 &$--$ & 9.544e-05 &$--$ & 7.499e-01 &$--$ \\
$N = 8$ & 3.442e-06 & 2.33 & 3.672e-02 & 2.90 & 2.643e-05 & 1.85& 1.999e-01 &1.91 \\
$N = 16$ & 5.338e-07 & 2.69 & 4.671e-03 &2.97 & 6.710e-06 & 1.98& 5.078e-02 &1.98 \\
$N = 32$ & 7.275e-08 & 2.88 & 5.864e-04 &2.99 & 1.678e-06 & 2.00& 1.275e-02 &1.99 \\
$N = 64$ & 9.286e-09 & 2.97 & 7.340e-05 &3.00 & 4.211e-07 & 1.99& 3.191e-03 &2.00 \\
\hline
\end{tabular}
\end{center}
\caption{$\mathbf{Example\ \ref{exam1}}$,  $L^2$, $H^1$ error and convergent rate of spatial discretization for the \textcolor{black}{P-BDF1} scheme \eqref{se3-1}-\eqref{se3-3} with $\tau=10^{-7},\, T=10^{-5}$.}\label{BDF1 SPACE ERROR}
\end{table}

\begin{table}[H]
\begin{center}
\begin{tabular}{ c|c c c c c c c c }
\hline
& $\left\| \phi - \phi_h \right\| _ {L^2}$ & & $\left\| \mathbf{u} - \mathbf{u}_h \right\| _ {L^2}$ & &$\left| \phi - \phi_h \right| _ {H^1}$ & &$\left| \mathbf{u} - \mathbf{u}_h \right| _ {H^1}$ &\\
\hline
$\tau = 0.04$ & 1.783e-01 & $--$ & 5.235e-02 & $--$ & 1.478e-01 & $--$ & 3.089e-02 & $--$ \\
$\tau = 0.02$ & 8.978e-02 & 0.99& 2.633e-02 & 0.99& 7.504e-02 & 0.98&1.530e-02 &1.01\\
$\tau = 0.01$ & 4.504e-02 & 1.00&1.320e-02 & 1.00&3.781e-02 & 0.99&7.646e-03 &1.00\\
$\tau = 0.005$ & 2.256e-02 & 1.00&6.610e-03 & 1.00&1.899e-02 & 0.99&3.839e-03 &0.99\\
\hline
\end{tabular}
\end{center}
\caption{$\mathbf{Example\ \ref{exam1}}$, $L^2$, $H^1$ error and convergent rate of temporal discretization for the \textcolor{black}{P-BDF1} scheme \eqref{se3-1}-\eqref{se3-3} with $N=160$, \textcolor{black}{$T=2$}.} \label{BDF1 TIME ERROR}
\end{table}


\begin{table}[H]
\centering
\begin{tabular}{ c|c c c c c c c c }
\hline
& $\left\| \phi - \phi_h \right\| _ {L^2}$ & & $\left\| \mathbf{u} - \mathbf{u}_h \right\| _ {L^2}$ & & $\left| \phi - \phi_h \right| _ {H^1}$ & &$\left| \mathbf{u} - \mathbf{u}_h \right| _ {H^1}$ & \\
\hline
$N = 4$ & 1.733e-06 & $--$ & 2.743e-01 & $--$ & 9.545e-06 & $--$& 7.499e-01 & $--$ \\
$N = 8$ & 3.442e-07 & 2.33 & 3.673e-02 & 2.90& 2.648e-06 &1.85 & 1.999e-01 &1.91\\
$N = 16$ & 5.313e-08 & 2.70& 4.672e-03 &2.97 & 6.782e-07 &1.97 & 5.079e-02 & 1.98\\
$N = 32$ & 7.196e-09 &2.88 & 5.866e-04 &2.99 & 1.684e-07 &2.01 & 1.275e-02 &1.99 \\
$N = 64$ & 9.267e-10 & 2.96& 7.341e-05 & 3.00& 4.211e-08 &2.00 & 3.191e-03 &2.00 \\
\hline
\end{tabular}
\caption{$\mathbf{Example\ \ref{exam1} }$, $L^2$, $H^1$ error and convergent rate of spatial discretization for the \textcolor{black}{P-BDF2} scheme \eqref{se3-18}-\eqref{se3-19} with $\tau=10^{-7},\, T=10^{-5}$.}\label{BDF2 SPACE ERROR}
\end{table}

\begin{table}[H]
\centering
\begin{tabular}{ c|c c c c c c c c }
\hline
& $\left\| \phi - \phi_h \right\| _ {L^2}$ & &  $\left\| \mathbf{u} - \mathbf{u}_h \right\| _ {L^2}$ & &  $\left| \phi - \phi_h \right| _ {H^1}$ & &  $\left| \mathbf{u} - \mathbf{u}_h \right| _ {H^1}$ &  \\
\hline
$\tau = 0.4$ & 7.754e-01 &$--$ &  1.585e-01 &$--$ &  5.547e-01 & $--$ & 1.204e-01 &$--$  \\
$\tau = 0.2$ & 2.022e-01 & 1.94 &  3.969e-02 & 2.00&  1.455e-01 &1.93 &  3.467e-02 & 1.80 \\
$\tau = 0.1$ & 5.051e-02 & 2.00 & 1.008e-02 & 1.98&  3.648e-02 & 2.00&  1.103e-02  &1.65 \\
$\tau = 0.05$ &  1.253e-02 & 2.01&  2.551e-03 & 1.98&  9.075e-03 &2.01 &  3.341e-03 &1.72  \\
\hline
\end{tabular}
\caption{$\mathbf{Example\ \ref{exam1} }$, $L^2$, $H^1$ error and convergent rate of temporal discretization for the \textcolor{black}{P-BDF2} scheme \eqref{se3-18}-\eqref{se3-19} with $N=160$, \textcolor{black}{$T=2$}.}\label{BDF2 TIME ERROR}
\end{table}

The errors and convergent rates for spatial discretization and temporal discretization between the numerical solution and exact solution based on the \textcolor{black}{P-BDF1} scheme \eqref{se3-1}-\eqref{se3-3} and \textcolor{black}{P-BDF2} scheme \eqref{se3-18}-\eqref{se3-19} are shown in Tables \ref{BDF1 SPACE ERROR}-\ref{BDF2 TIME ERROR}, respectively. From the Tables, we numerically verify that for the \textcolor{black}{P-BDF1} scheme 
\begin{equation}
\begin{aligned}
&\left\|\phi(t_{n+1})-\phi_{h}^{n+1}\right\|_{L^{2}(\Omega)} \approx O(h^3 + \tau),\qquad
&\left\|\mathbf{u}(t_{n+1})-\mathbf{u}_{h}^{n+1}\right\|_{(L^{2}(\Omega))^{2}} \approx O(h^3 + \tau), \\
&\left|\phi(t_{n+1})-\phi_{h}^{n+1}\right|_{H^{1}(\Omega)} \approx O(h^2 + \tau),\qquad
&\left|\mathbf{u}(t_{n+1})-\mathbf{u}_{h}^{n+1}\right|_{(H^{1}(\Omega))^{2}} \approx O(h^2 + \tau),
\end{aligned}    
\end{equation}
and for the \textcolor{black}{P-BDF2} scheme 
\begin{equation}
\begin{aligned}
&\left\|\phi(t_{n+1})-\phi_{h}^{n+1}\right\|_{L^{2}(\Omega)} \approx O(h^3 + \tau^2),\qquad
&\left\|\mathbf{u}(t_{n+1})-\mathbf{u}_{h}^{n+1}\right\|_{(L^{2}(\Omega))^{2}} \approx O(h^3 + \tau^2), \\
&\left|\phi(t_{n+1})-\phi_{h}^{n+1}\right|_{H^{1}(\Omega)} \approx O(h^2 + \tau^2),\qquad
&\left|\mathbf{u}(t_{n+1})-\mathbf{u}_{h}^{n+1}\right|_{(H^{1}(\Omega))^{2}} \approx O(h^2 + \tau^2),
\end{aligned}    
\end{equation}
which are consistent with the expectations. 

\subsection{Numerical solutions for the CHNS equations}

In the following four examples, we numerically investigate the performance of the proposed IEQ-FEM schemes for solving the CHNS equations. We conduct the computation by using schemes \eqref{se3-1}-\eqref{se3-3} and \eqref{se3-18}-\eqref{se3-19} separately and present the results as phase field\textcolor{black}{s} and velocity field\textcolor{black}{s}. Besides, the energy dissipation and mass conservation phenomena are also displayed.

\begin{example} \cite{zhao_second-order_2021} \label{exam4}
In this example, we consider the CHNS equations \eqref{eq1.3} with the domain $\Omega = \left[0, 1\right]^2 $ and the initial condition
\begin{equation}
\begin{aligned}
\phi_0 &= 1 - \tanh\left({\frac{-r+\sqrt{(x-x_a)^2+(y-y_a)^2}}{2\epsilon}}\right) - \tanh\left({\frac{-r+\sqrt{(x-x_b)^2+(y-y_b)^2}}{2\epsilon}}\right) ,\\
\mathbf{u}_0 &= \left[ 0, \; 0 \right]^\top , 
\end{aligned}
\end{equation}
where $x_a=0.5-\frac{r}{\sqrt{2}},\;y_a=0.5+\frac{r}{\sqrt{2}},\;x_b=0.5+\frac{r}{\sqrt{2}},\;y_b=0.5-\frac{r}{\sqrt{2}},\;r=0.15$. The parameters are chosen as $\gamma = \mu = \epsilon = 0.01, \; \lambda = 0.01\epsilon,\; B=100$. We set the \textcolor{black}{partition number $N \times N = 128 \times 128$}, 
and time step $\tau = 5 \times 10^{-4}$ with total time $T=3.2$. 
\end{example}

\begin{figure}[!htbp]
$\begin{array}{c} 
\includegraphics[width=5.5cm,height=4.5cm]{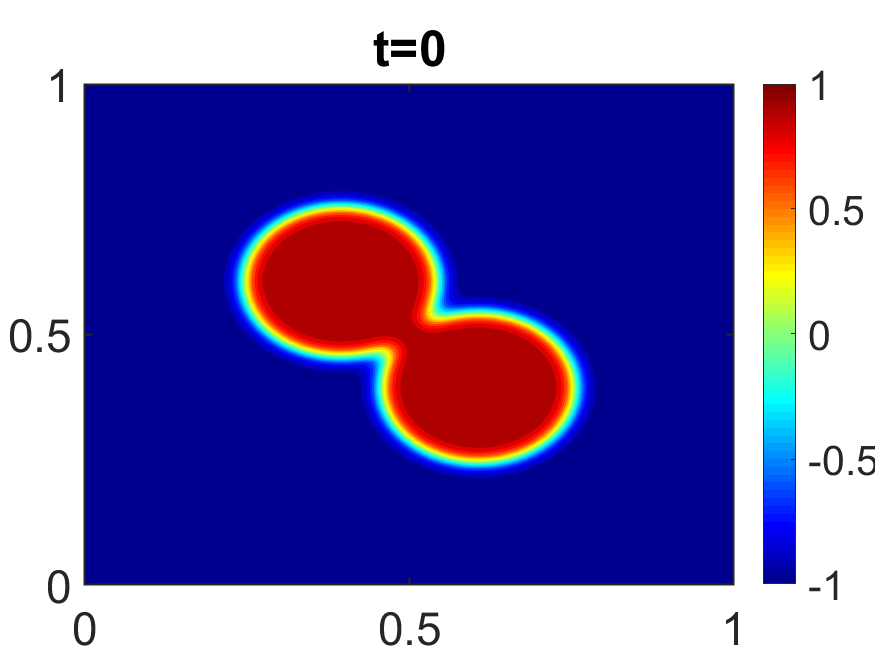}
\includegraphics[width=5.5cm,height=4.5cm]{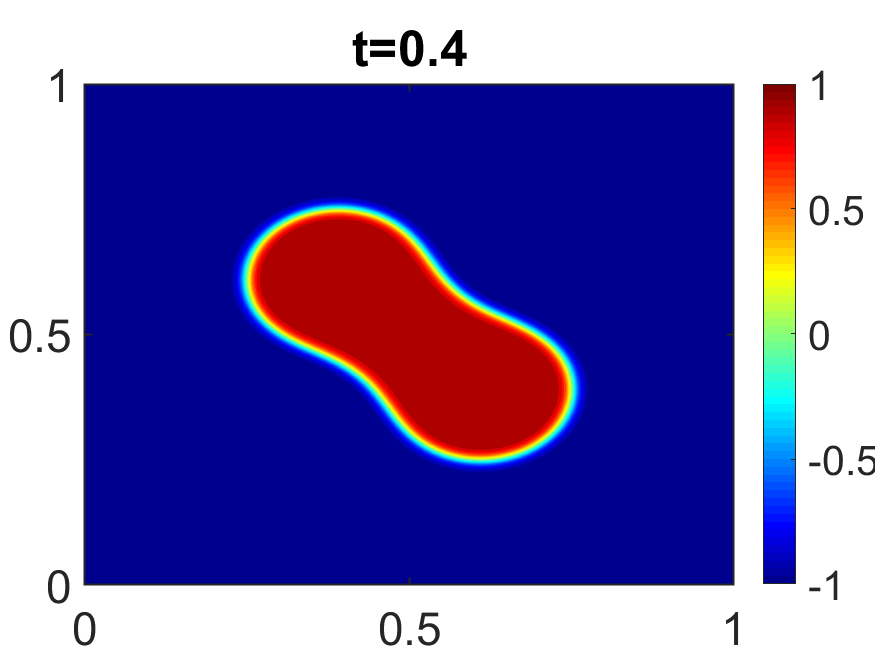}
\includegraphics[width=5.5cm,height=4.5cm]{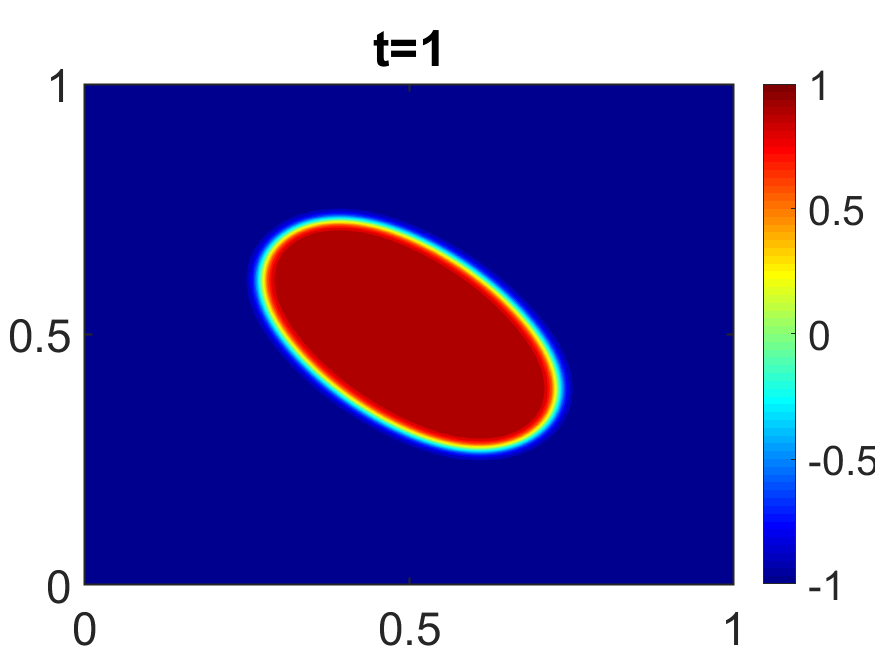}\\ 
\includegraphics[width=5.5cm,height=4.5cm]{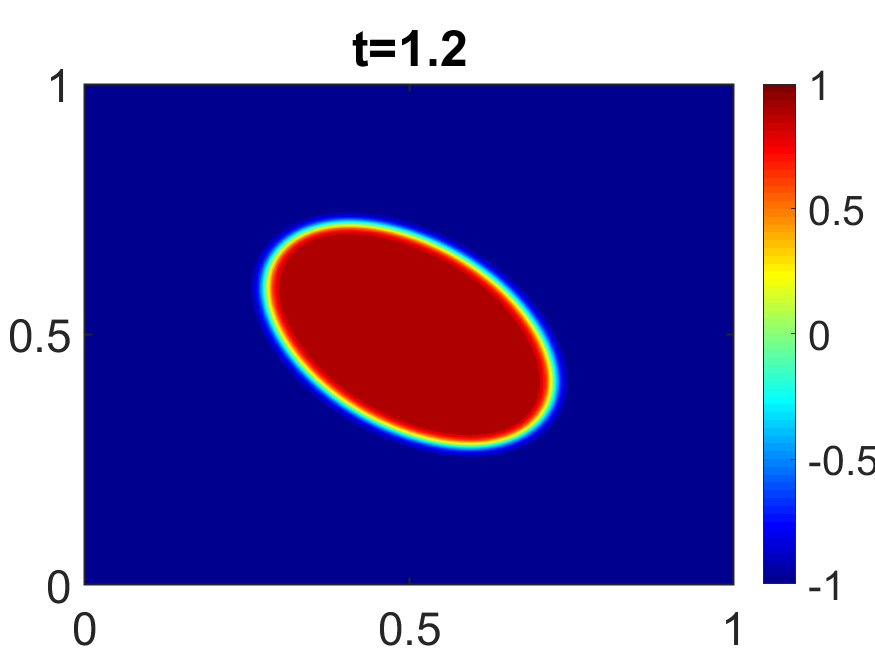}
\includegraphics[width=5.5cm,height=4.5cm]{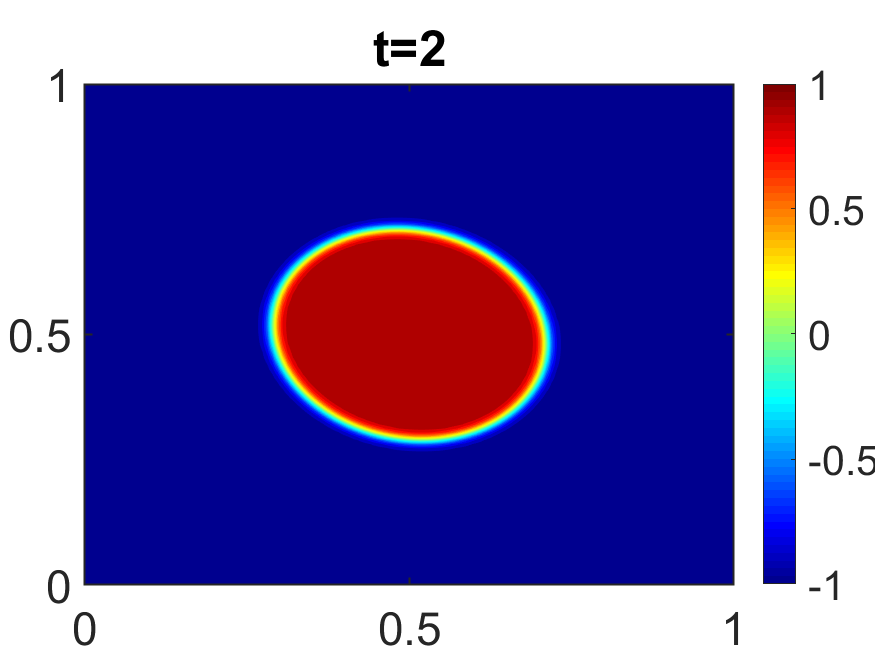}
\includegraphics[width=5.5cm,height=4.5cm]{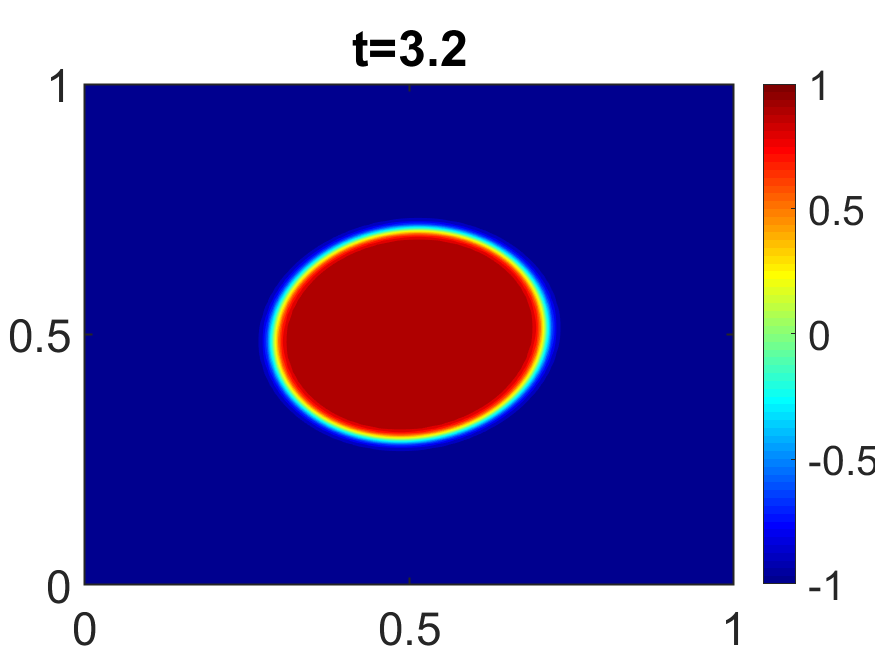}\\
\includegraphics[width=5.5cm,height=4.5cm]{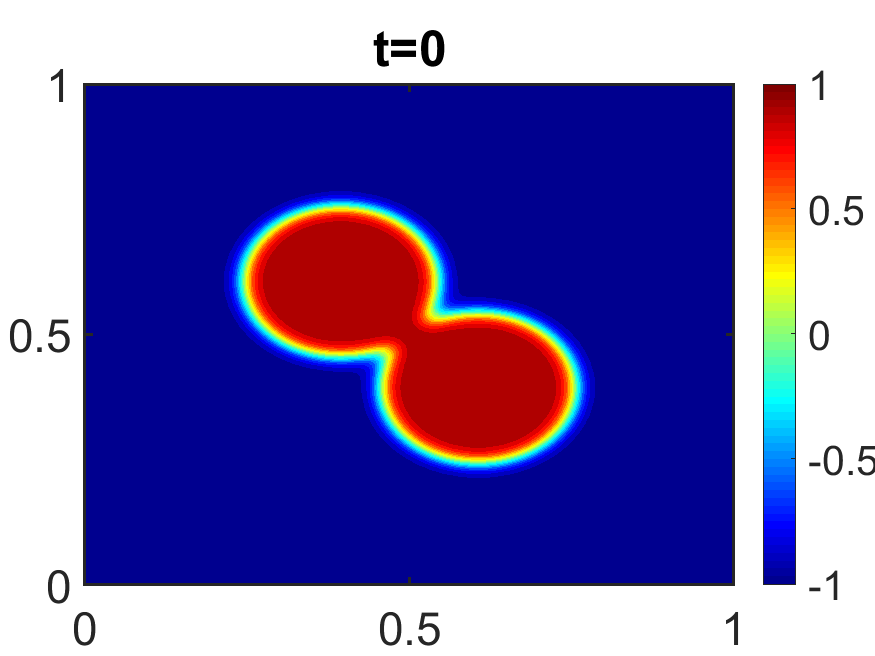}
\includegraphics[width=5.5cm,height=4.5cm]{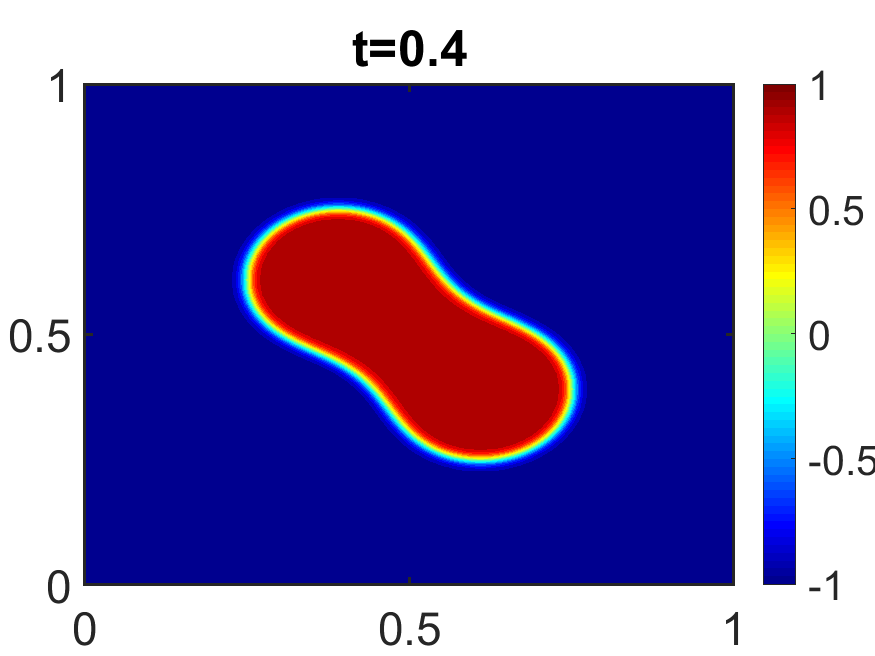}
\includegraphics[width=5.5cm,height=4.5cm]{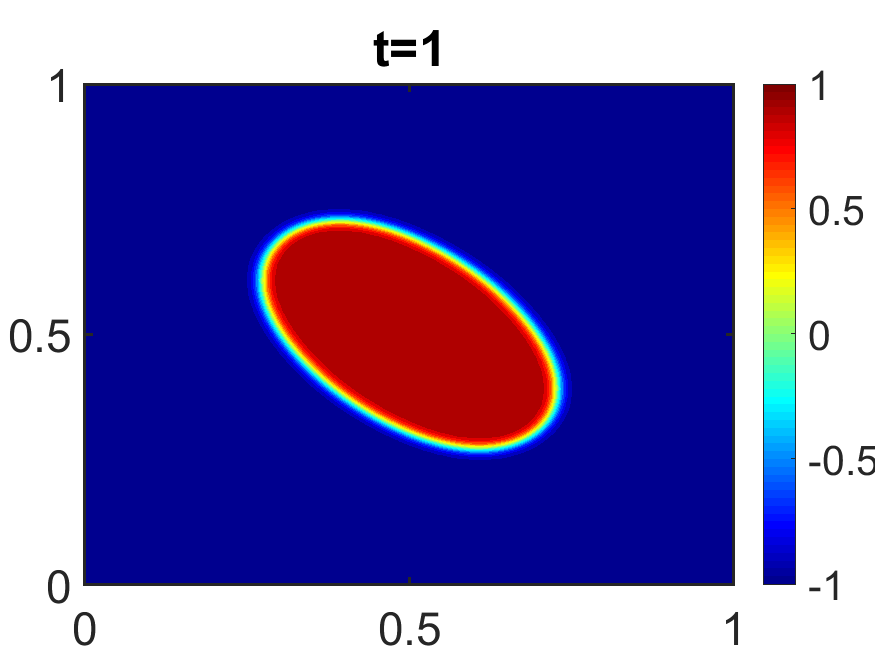}\\ 
\includegraphics[width=5.5cm,height=4.5cm]{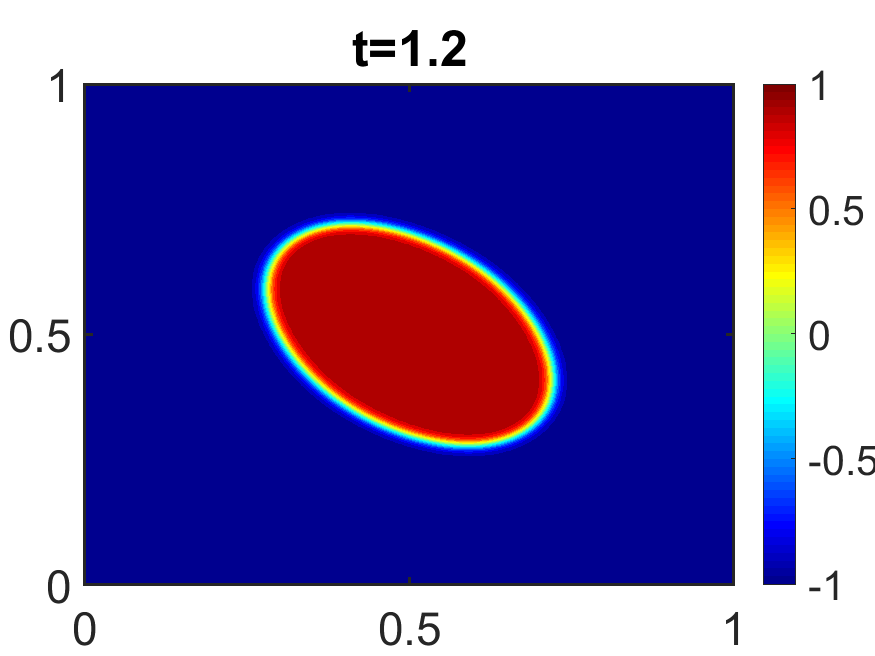}
\includegraphics[width=5.5cm,height=4.5cm]{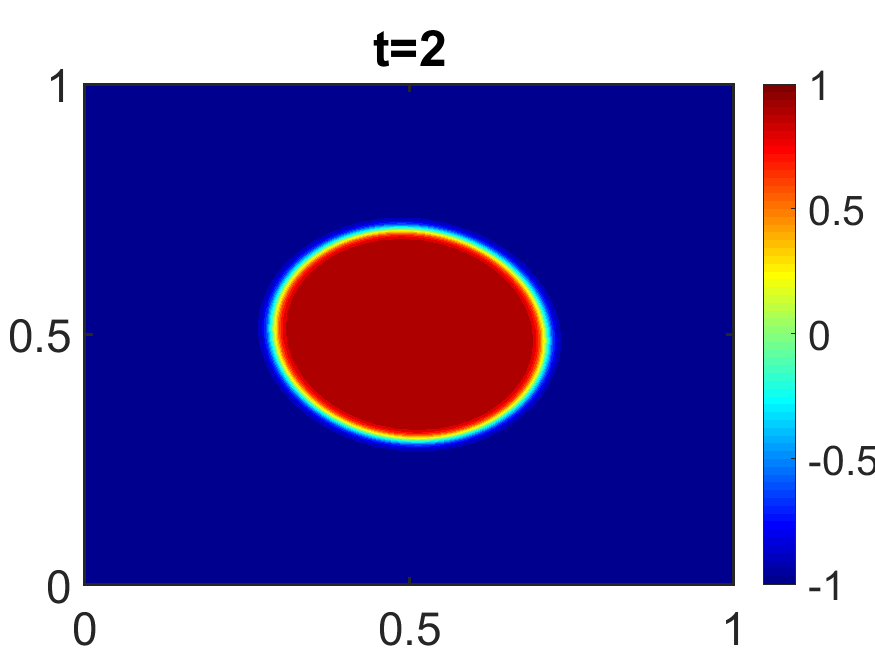}
\includegraphics[width=5.5cm,height=4.5cm]{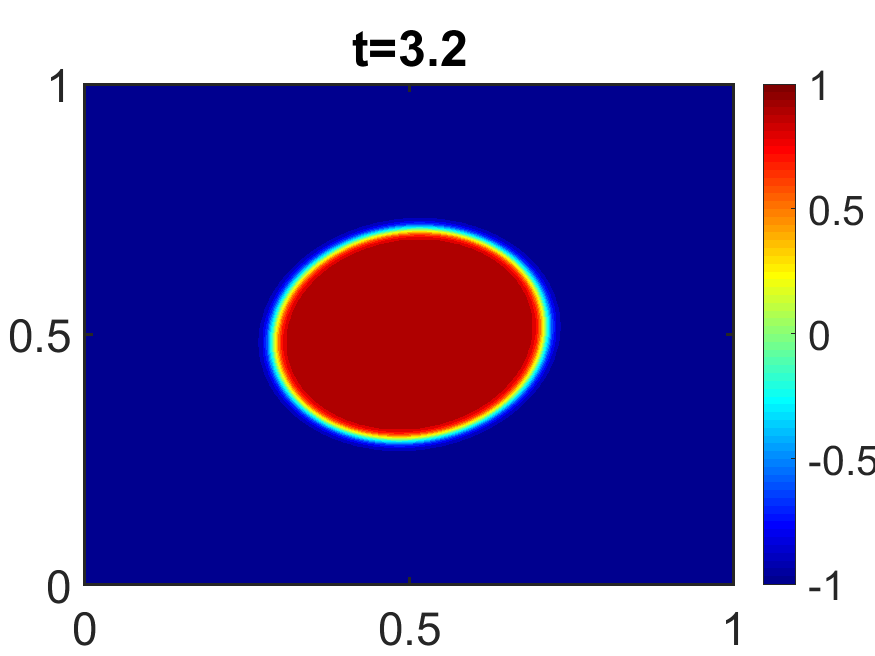}
\end{array}$\vspace{-0.2cm}
\caption{$\mathbf{Example\ \ref{exam4}}$, snapshots of numerical solutions for phase field\textcolor{black}{s} function, First and second lines: \textcolor{black}{P-BDF1} scheme \eqref{se3-1}-\eqref{se3-3}; Third and fourth lines: \textcolor{black}{P-BDF2} scheme \eqref{se3-18}-\eqref{se3-19}. } \label{ex2 BDF2 phi plot}
\end{figure}

\begin{figure}[!htbp]
$\begin{array}{c} 
\includegraphics[width=5.5cm,height=4.5cm]{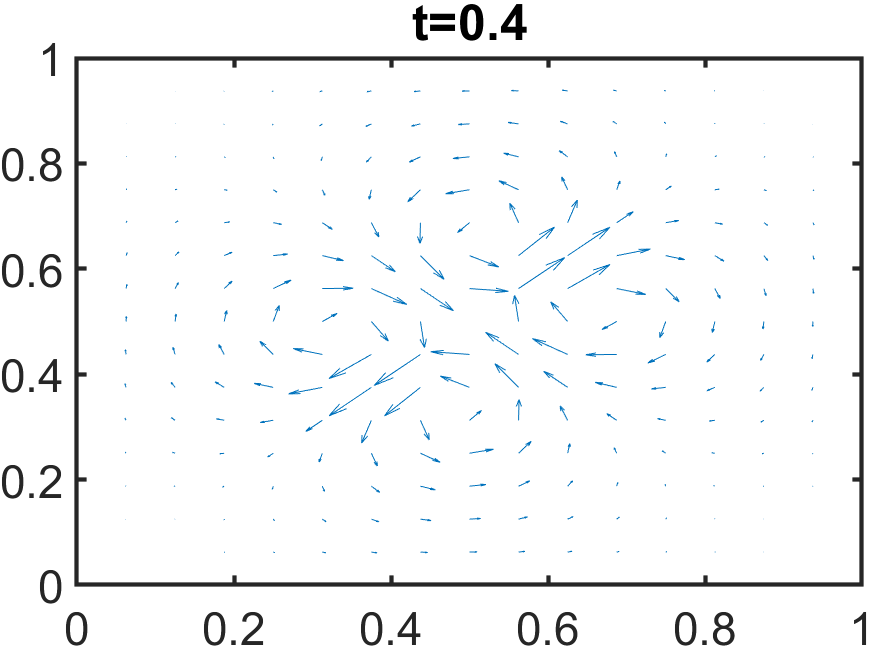}\;
\includegraphics[width=5.5cm,height=4.5cm]{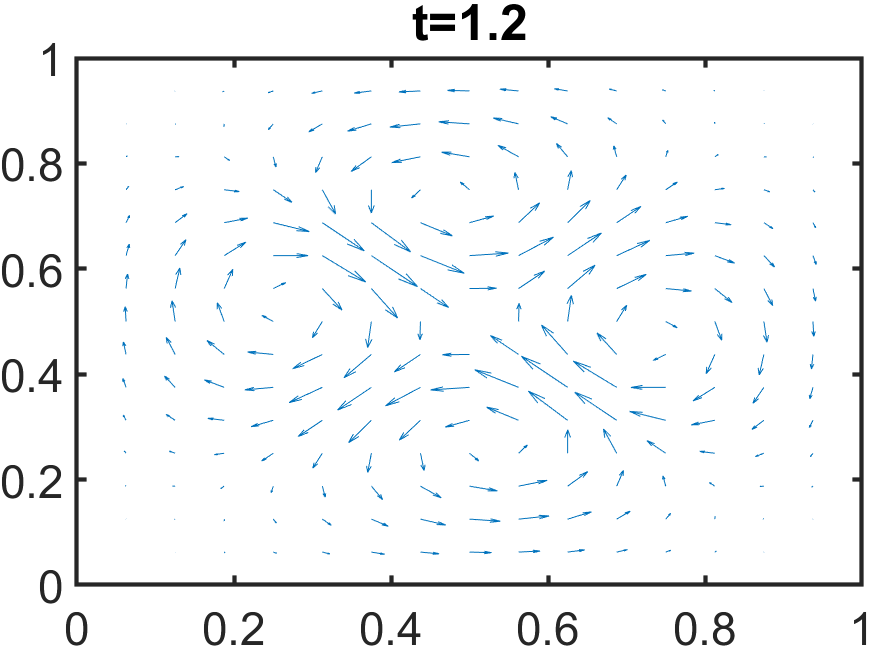}\;
\includegraphics[width=5.5cm,height=4.5cm]{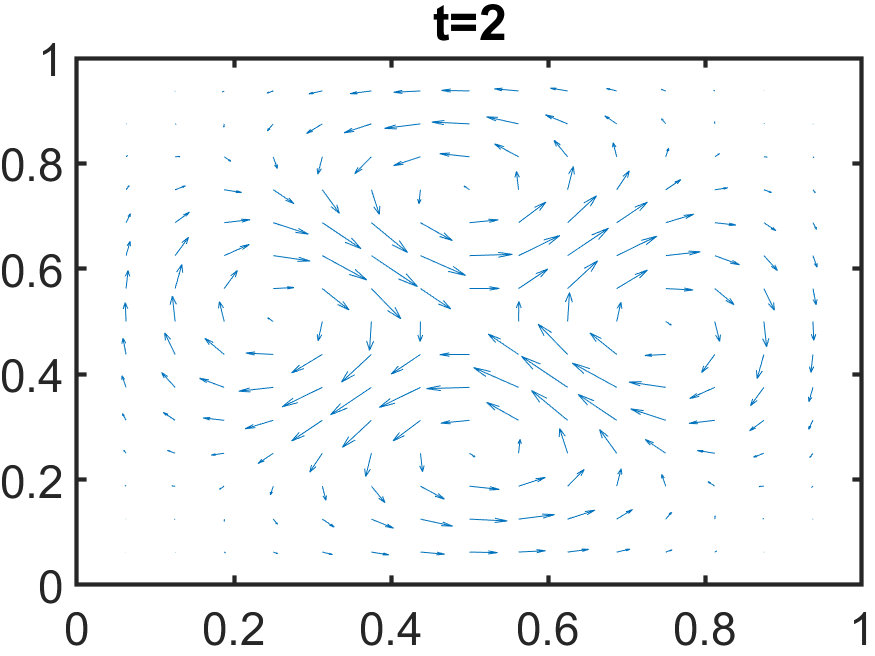}\\

\includegraphics[width=5.5cm,height=4.5cm]{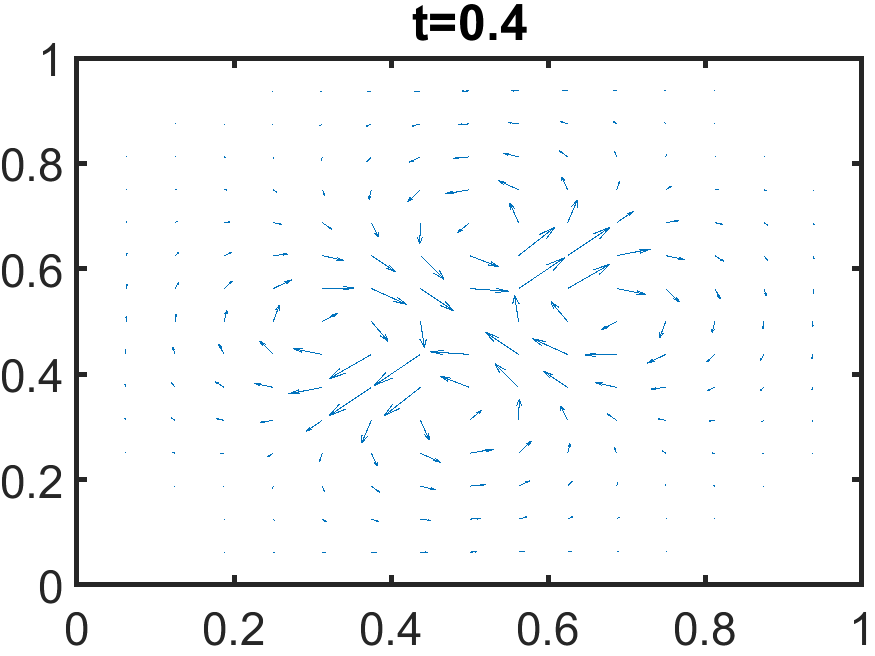}\;
\includegraphics[width=5.5cm,height=4.5cm]{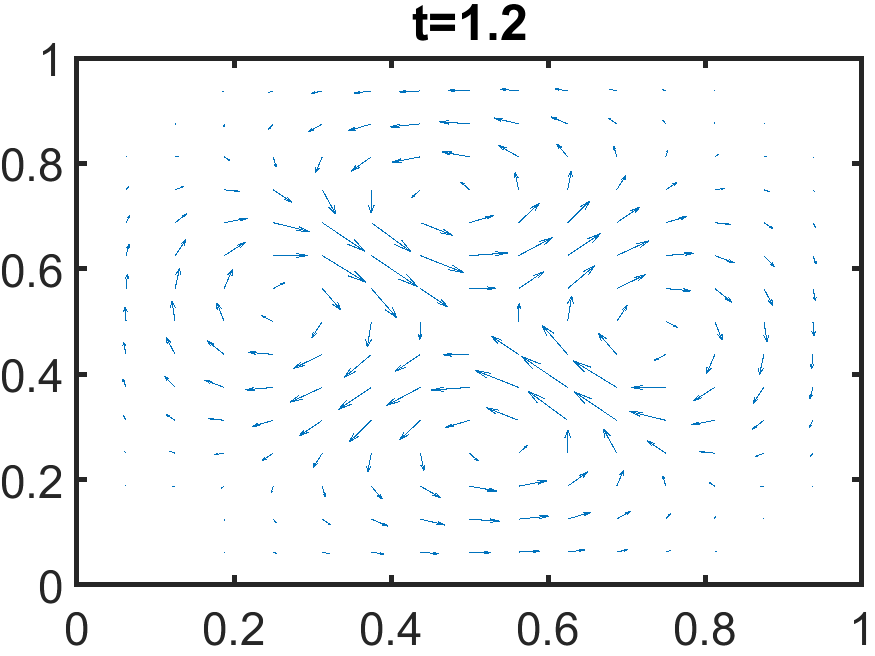}\;
\includegraphics[width=5.5cm,height=4.5cm]{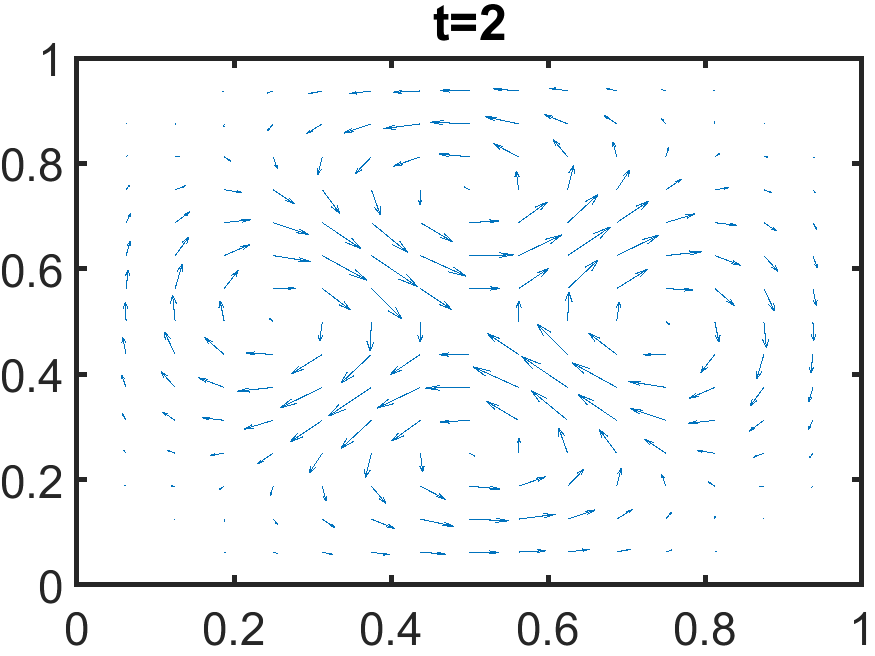}
\end{array}$\vspace{-0.2cm}
\caption{$\mathbf{Example\ \ref{exam4}}$, snapshots of numerical solutions for velocity field\textcolor{black}{s} function, First line: \textcolor{black}{P-BDF1} scheme \eqref{se3-1}-\eqref{se3-3}; Second line: \textcolor{black}{P-BDF2} scheme \eqref{se3-18}-\eqref{se3-19}. 
} \label{ex2 BDF2 u}
\end{figure}

\begin{figure}[!htbp]
$\begin{array}{c}
\includegraphics[width=8cm,height=6.5cm]{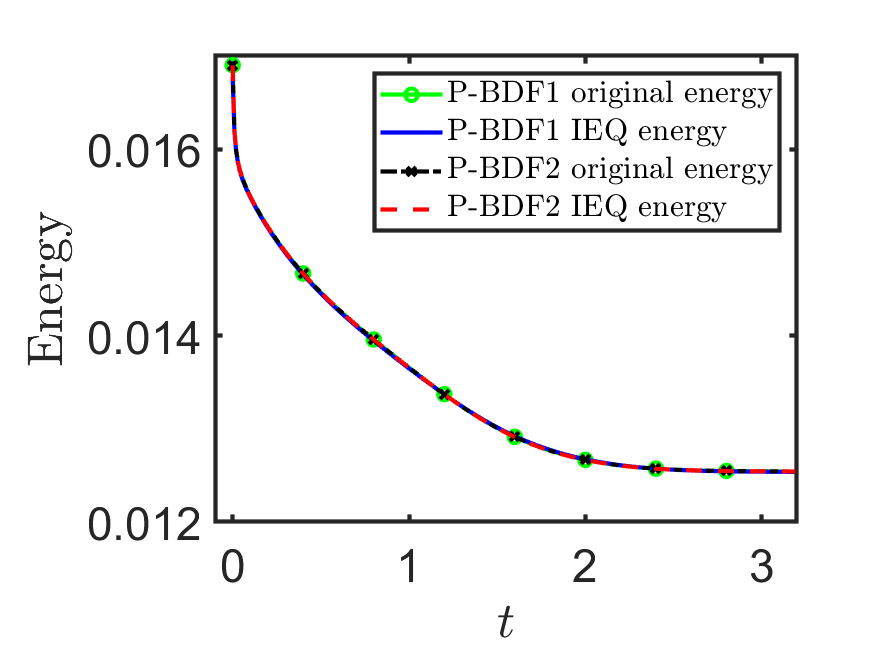}
\includegraphics[width=8cm,height=6.5cm]{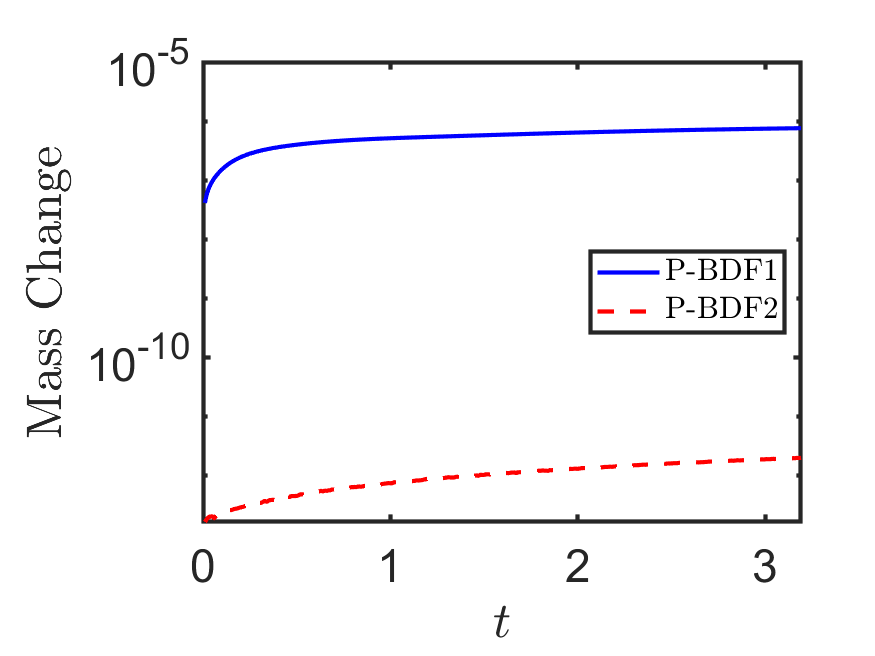}
\end{array}$
\caption{$\mathbf{Example\ \ref{exam4}}$, the modified discrete energy history and discrete mass \textcolor{black}{variation}.} \label{ex2 Energy and mass}
\end{figure}

Figures \ref{ex2 BDF2 phi plot}-\ref{ex2 BDF2 u} show the numerical solution of phase field\textcolor{black}{s} $\phi_h^{n+1}$ and velocity field\textcolor{black}{s} $\mathbf{u}_h^{n+1}$ by using the \textcolor{black}{P-BDF1} scheme \eqref{se3-1}-\eqref{se3-3} and \textcolor{black}{P-BDF2} scheme \eqref{se3-18}-\eqref{se3-19}. From the pictures, we can see that the numerical solutions of phase field\textcolor{black}{s} and velocity field\textcolor{black}{s} obtained by the two schemes are almost identical. We can see that the patterns are comparable to those in \cite{zhao_second-order_2021}. The evolution of the discrete energy and  total mass are shown in Figure \ref{ex2 Energy and mass}. The results indicate that both two schemes preserve  the total mass and satisfy the energy dissipation law.
\begin{table}
\centering
\begin{tabular}{c|cc}
\hline
{Example\ \ref{exam4} } & P scheme & C scheme  \\
\hline
BDF1 & 535.0966s & 496.5006s  \\
BDF2 & 2011.748s & 1972.774s  \\
\hline
\end{tabular}
\caption{CPU time for 100 steps calculated using four schemes: (1) scheme \eqref{se3-1}-\eqref{se3-3}; (2) scheme \eqref{se3-18}-\eqref{se3-19}; (3) scheme \eqref{apdxA2.1}-\eqref{apdxA2.3}; (4) scheme \eqref{apdx3-18-1}-\eqref{apdx3-19-1}. }
\label{tab:my_label}
\end{table}

Besides, we also compare the CPU time of four classes of schemes: (1) \textcolor{black}{P-BDF1} scheme; (2) \textcolor{black}{P-BDF2} scheme; (3) \textcolor{black}{C-BDF1} scheme; (4) \textcolor{black}{C-BDF2} scheme, and the corresponding results are shown in Table \ref{tab:my_label}.  The results show that the additional projection step $\eqref{se3-1 projection}$ for BDF1 scheme and $\eqref{se3-18 projection}$ for BDF2 scheme in the CHNS equations slightly increases the computation time, which is acceptable.


In the following, we further present more examples based on the \textcolor{black}{P-BDF2} scheme to investigate the performance of the proposed methods.

\begin{example}\label{exam2}
In this example, we consider the CHNS equations \eqref{eq1.3} with the domain $\Omega = \left[-1, 1\right]^2 $ and the initial condition 
\begin{equation}
\begin{aligned}
\phi_0 &= \tanh \left( \left( \left( x-0.3 \right)^2 +y^2 -0.2^2  \right)/ \epsilon^2 \right) \times \tanh \left( \left( \left( x+0.3 \right)^2 +y^2 -0.2^2  \right)/ \epsilon^2 \right) \times \\
&\quad\tanh \left( \left( x^2 + \left( y-0.3 \right)^2 -0.2^2  \right)/ \epsilon^2 \right) \times \tanh \left( \left( x^2 + \left( y+0.3 \right)^2 -0.2^2  \right)/ \epsilon^2 \right), \\
\mathbf{u}_0 &= \left[ \sin \left(\pi x \right)^2 \sin \left(2\pi y \right),
\quad \sin \left(2\pi x \right) \sin \left(\pi y \right)^2 \right]^\top ,
\end{aligned}
\end{equation}
where the parameters $\epsilon = \lambda = 0.04,\;\mu = \gamma = 1$, $B=50$, the mesh $N \times N = 80 \times 80$, time step $\tau = 10^{-6}$ and the final time $T=10^{-1}$.  
\end{example}

\begin{figure}[!htbp]
$\begin{array}{c}
\includegraphics[width=5.5cm,height=4.5cm]{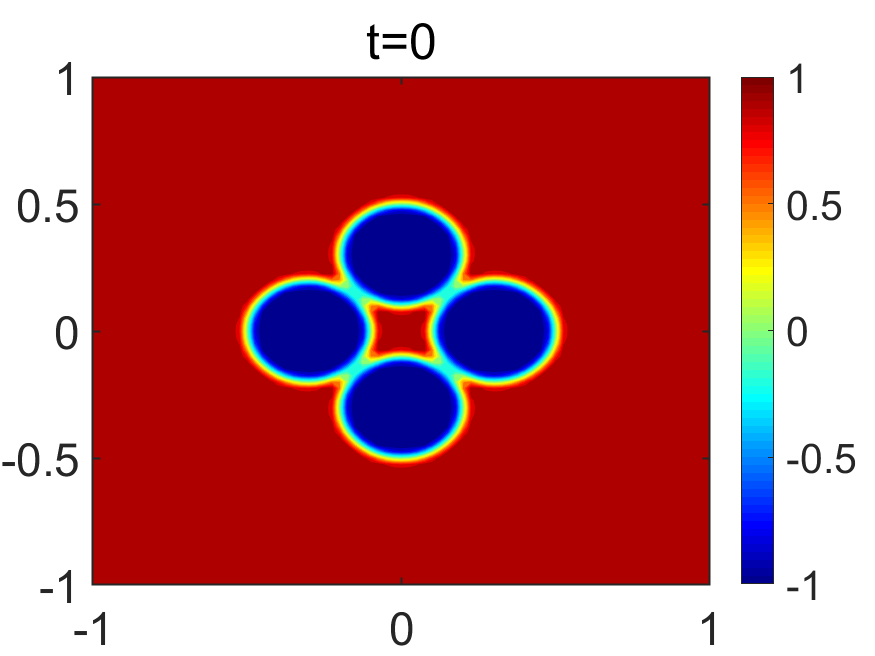} \;
\includegraphics[width=5.5cm,height=4.5cm]{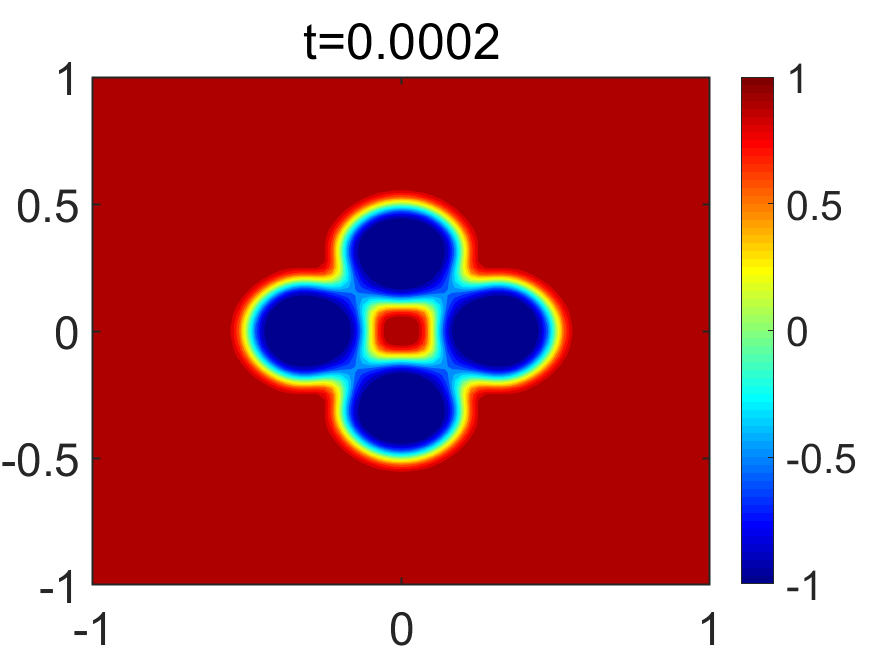} \;
\includegraphics[width=5.5cm,height=4.5cm]{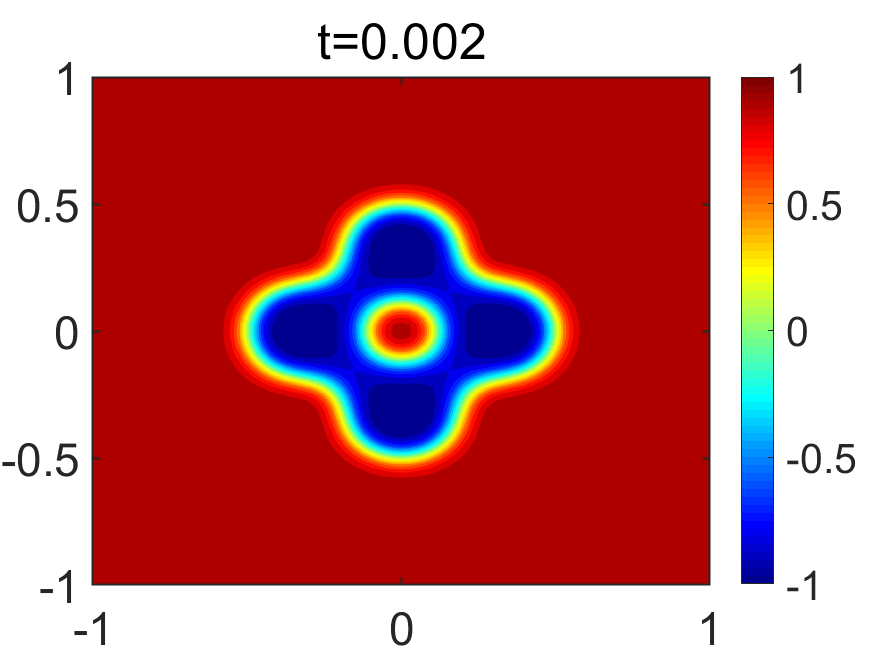} \\
\includegraphics[width=5.5cm,height=4.5cm]{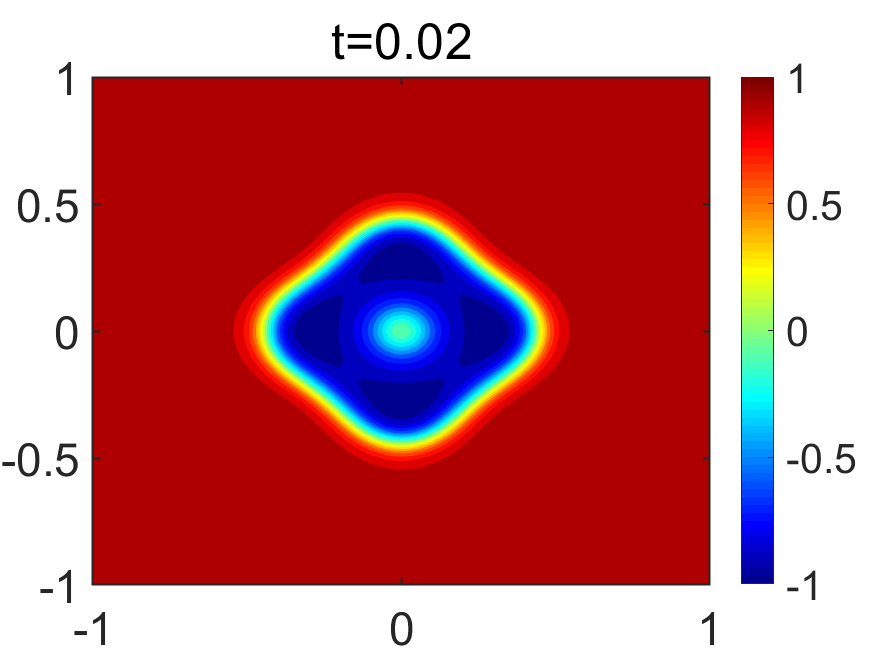}\; 
\includegraphics[width=5.5cm,height=4.5cm]{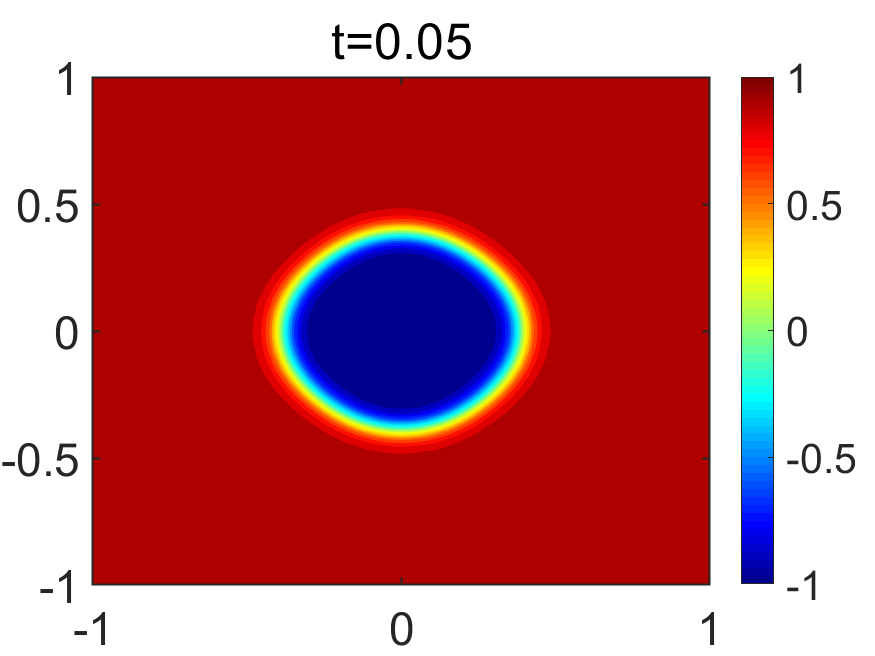}\;
\includegraphics[width=5.5cm,height=4.5cm]{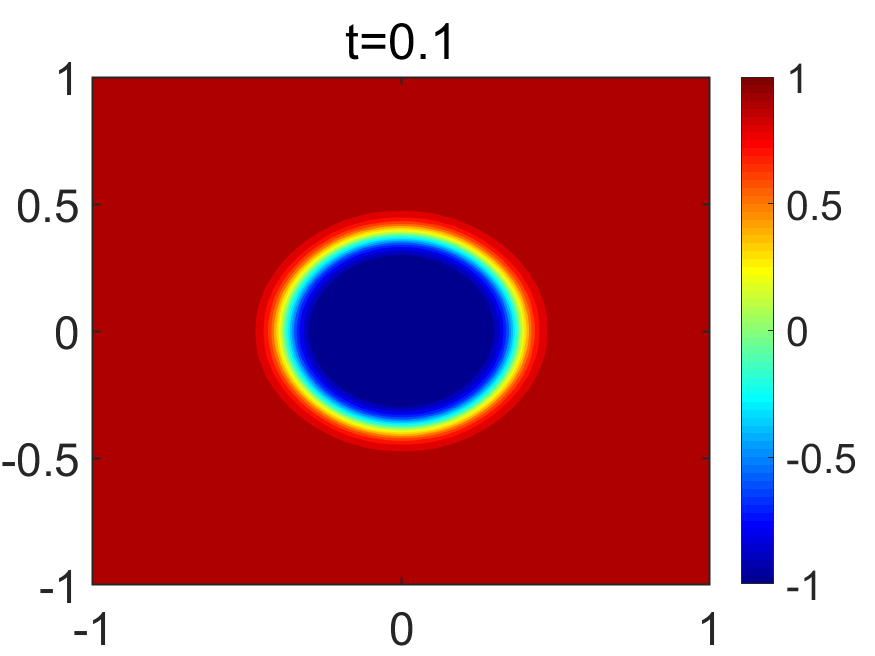}
\end{array}$\vspace{-0.2cm}
\caption{$\mathbf{Example\ \ref{exam2}}$, \textcolor{black}{P-BDF2} scheme \eqref{se3-18}-\eqref{se3-19}, snapshots of numerical solutions for \textcolor{black}{the} phase field function.}\label{Cexp3phi}
\end{figure}


\begin{figure}[!htbp]
$\begin{array}{c}
\includegraphics[width=5.5cm,height=4.5cm]{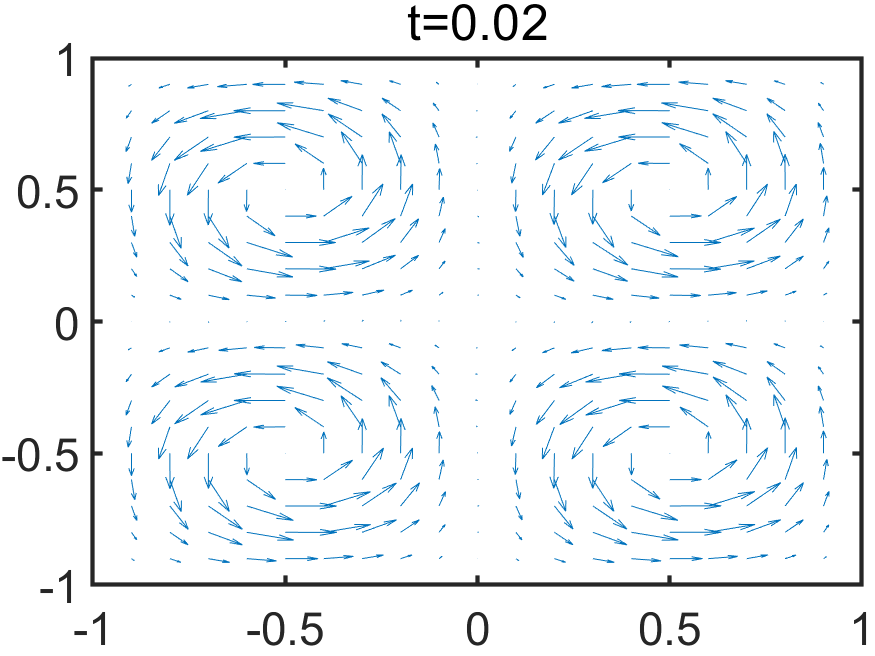}\; 
\includegraphics[width=5.5cm,height=4.5cm]{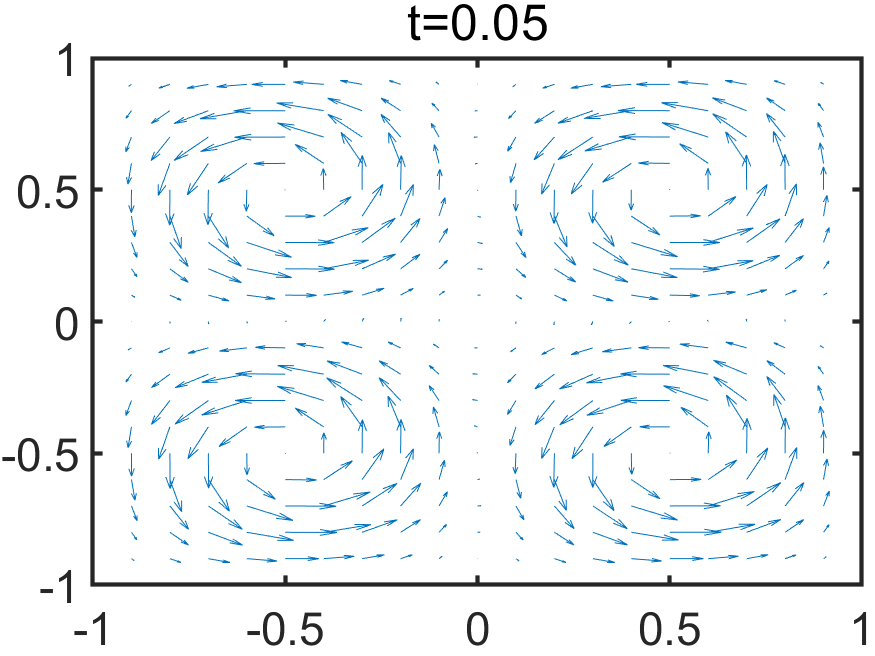}\;
\includegraphics[width=5.5cm,height=4.5cm]{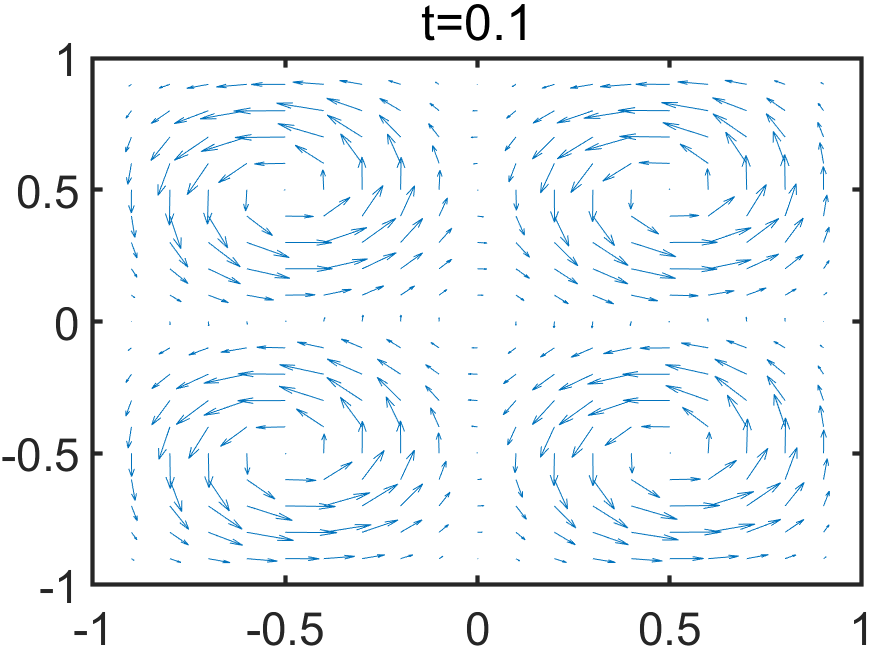}
\end{array}$\vspace{-0.2cm}
\caption{$\mathbf{Example\ \ref{exam2}}$, \textcolor{black}{P-BDF2} scheme \eqref{se3-18}-\eqref{se3-19}, snapshots of numerical solutions for \textcolor{black}{the} velocity field function.}\label{ex3_u}
\end{figure}

\begin{figure}[!htbp]
$\begin{array}{c} \includegraphics[width=8cm,height=6.5cm]{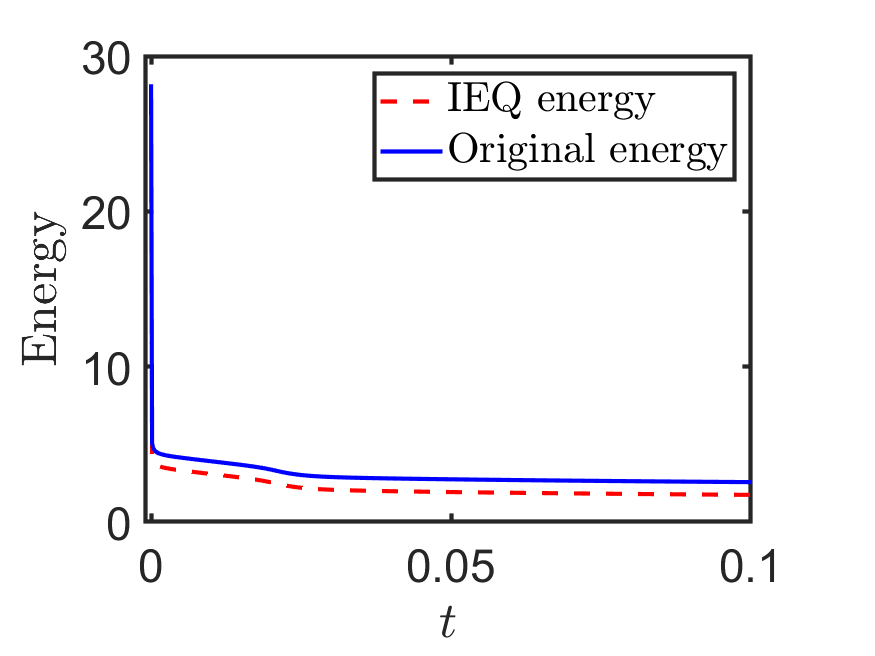} \;
\includegraphics[width=8cm,height=6.2cm]{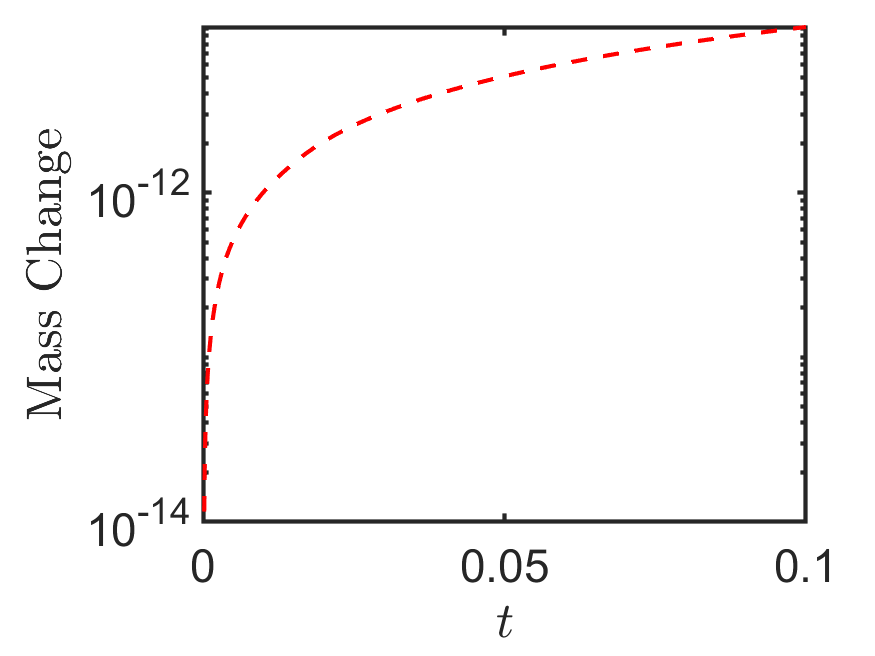}
\end{array}$\vspace{-0.2cm}
\caption{$\mathbf{Example\ \ref{exam2}}$, the modified discrete energy history and the change of total mass.}\label{Cexp3 energy and mass}
\end{figure}

A sequence snapshots of the approximate solutions for \textcolor{black}{the} phase field function and the velocity field function are produced in Figures \ref{Cexp3phi}-\ref{ex3_u} by the proposed \textcolor{black}{P-BDF2} scheme. It is easy to see the numerical solutions satisfy the expectation. 
The graphs depicting the evolution of discrete energy and change of total mass for Example \ref{exam2} are presented in
Figure \ref{Cexp3 energy and mass}. As illustrated, the discrete energy decreases over time and the variation in mass approaches machine precision.

\begin{example}\label{exam3}
Let domain $\Omega = \left[-2, 2\right]^2 $, define $m_1 = \left[0,2 \right],\; m_2 = \left[0,0 \right],\; m_3 = \left[0,-2 \right]$. For given $\epsilon = \frac{1}{16}$, let $ \; r_1 = r_3 = 2-\frac{3\epsilon}{2},\; r_2 =1$ and set $d\left(x \right) = \max\left\{ -d_1\left(x \right), d_2\left(x \right), d_3\left(x \right)\right\}$, $d_j\left(x \right) = \left|x-m_j\right|-r_j$ for $j=1,2,3$, we consider the CHNS equations \eqref{eq1.3} with the following initial condition 
\begin{equation}
\begin{aligned}
\phi_0 &= - \tanh\left(\frac{d\left( x\right)}{\sqrt{2}\epsilon} \right),  \\
\mathbf{u}_0 &= C\left[ \sin \left(\pi x \right)^2 \sin \left(2\pi y \right) \quad \sin \left(2\pi x \right) \sin \left(\pi y \right)^2 \right]^\top, 
\end{aligned}   
\end{equation}
here the parameters $\mu = \gamma =  1,\;\lambda = \frac{1}{16}, \;C=100$. We set time step $\tau = 10 ^{-5}$, 
the mesh $N \times N = 80 \times 80$, and $B=1$.
\end{example}

\begin{figure}[!htbp]
$\begin{array}{c} \hspace{-0.5cm}
\includegraphics[width=5.5cm,height=4.5cm]{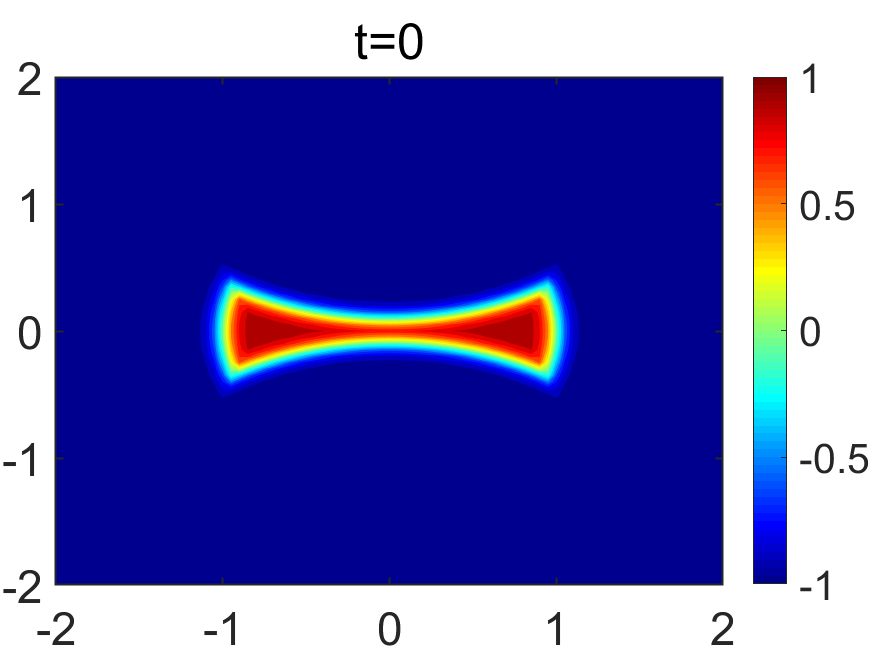}\;
\includegraphics[width=5.5cm,height=4.5cm]{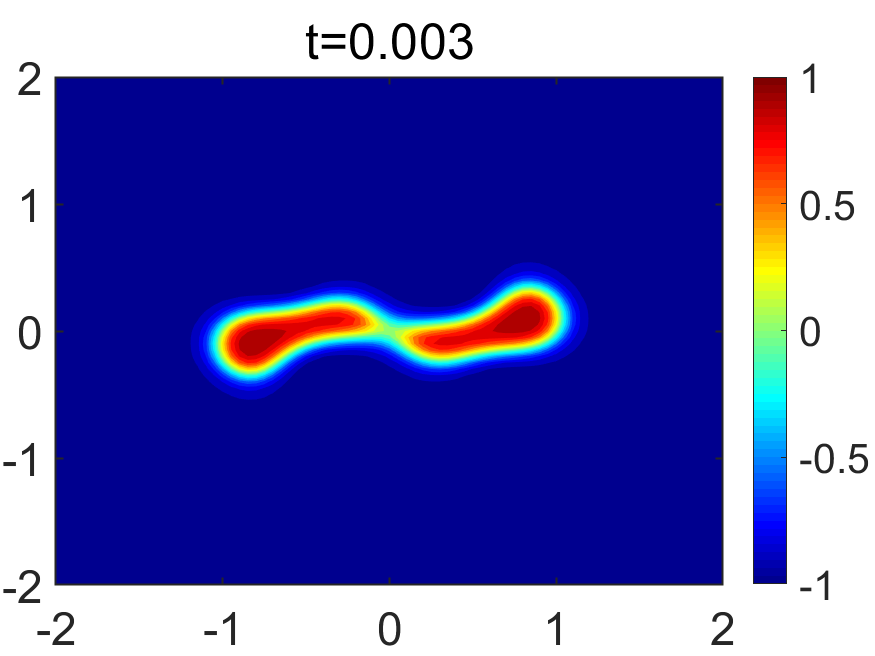}\;
\includegraphics[width=5.5cm,height=4.5cm]{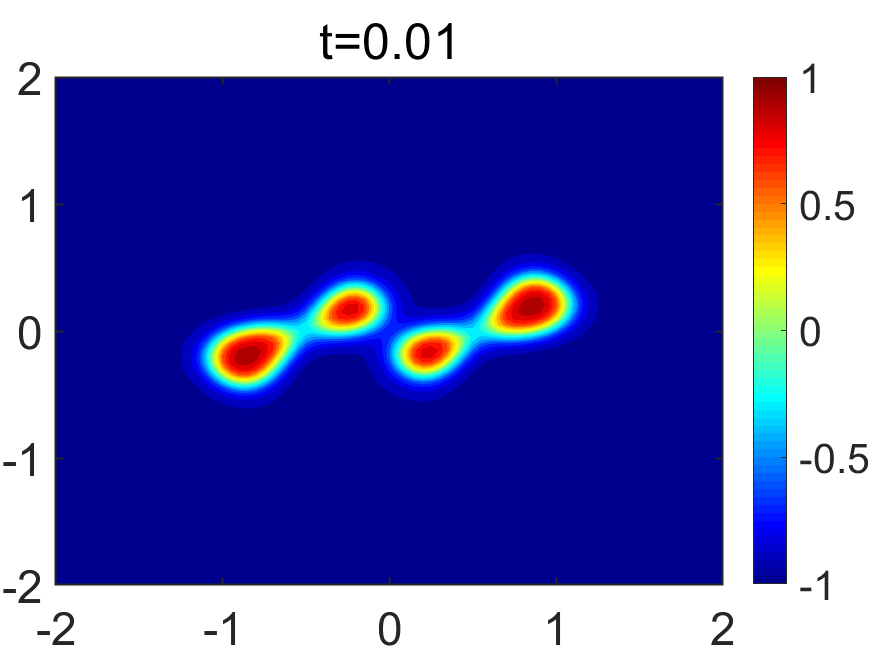}\;
\\ \hspace{-0.5cm}
\includegraphics[width=5.5cm,height=4.5cm]{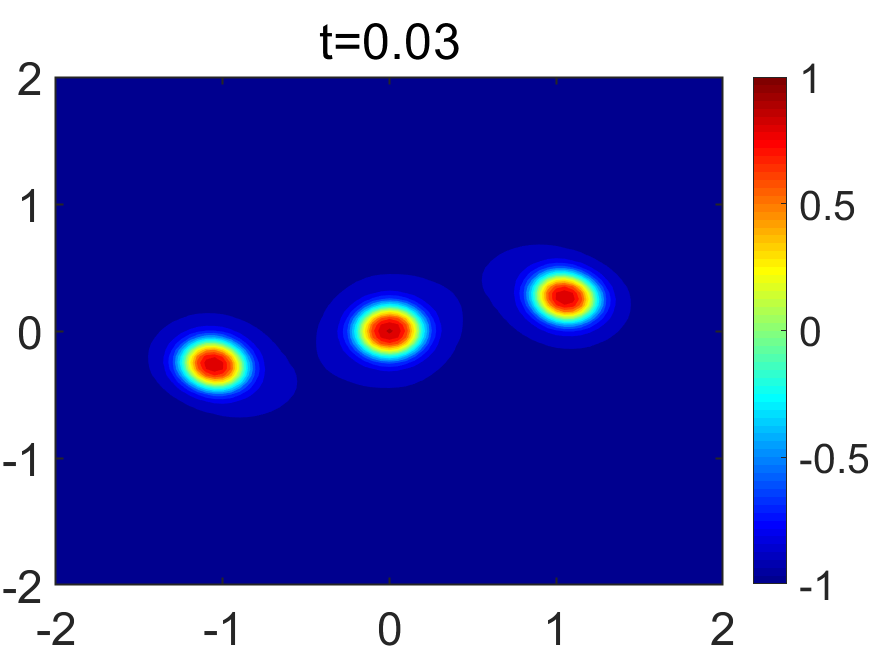}\;
\includegraphics[width=5.5cm,height=4.5cm]{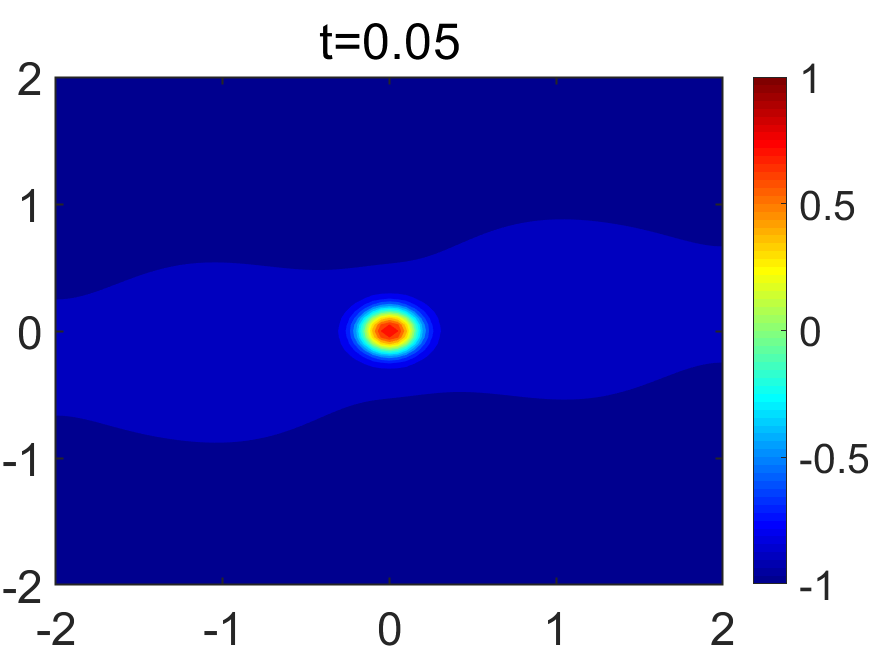}\;
\includegraphics[width=5.5cm,height=4.5cm]{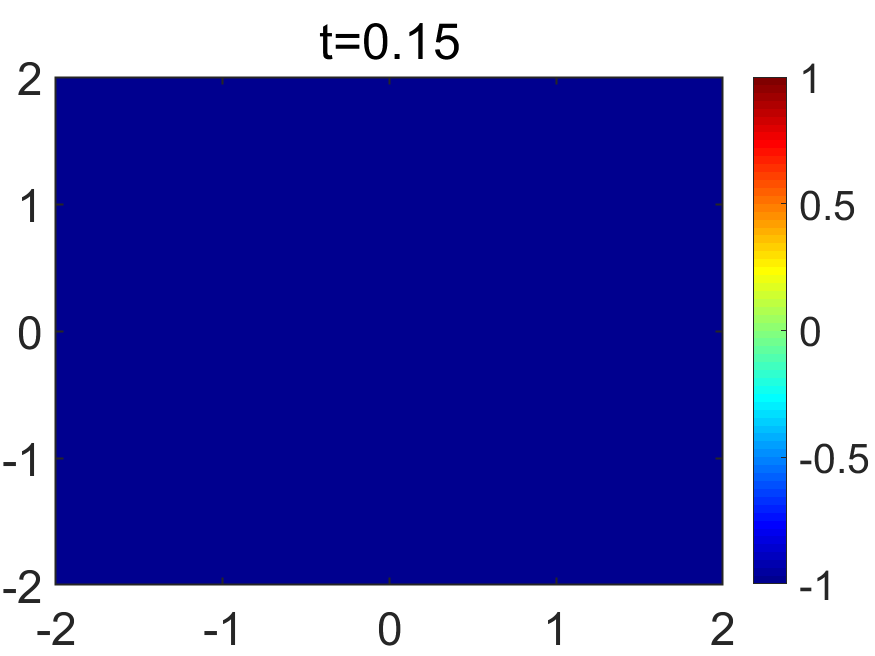}
\end{array}$\vspace{-0.2cm}
\caption{$\mathbf{Example\ \ref{exam3}}$, \textcolor{black}{P-BDF2} scheme \eqref{se3-18}-\eqref{se3-19}, snapshots of numerical solutions for \textcolor{black}{the} phase field function. 
}\label{Cexp4phi}
\end{figure}


\begin{figure}[!htbp]
$\begin{array}{c} \hspace{-1cm}
\includegraphics[width=5.5cm,height=4.5cm]{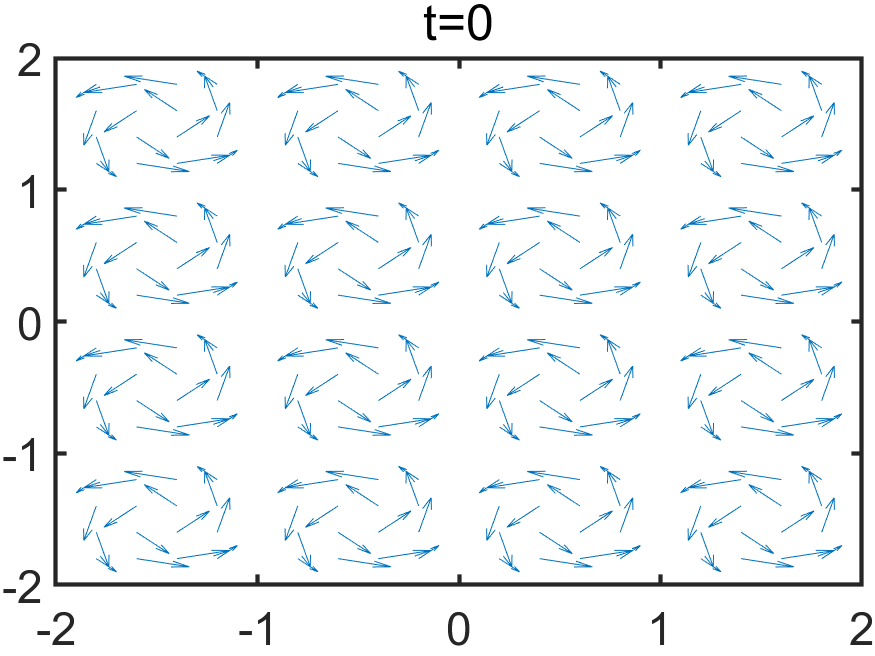}
\includegraphics[width=5.5cm,height=4.5cm]{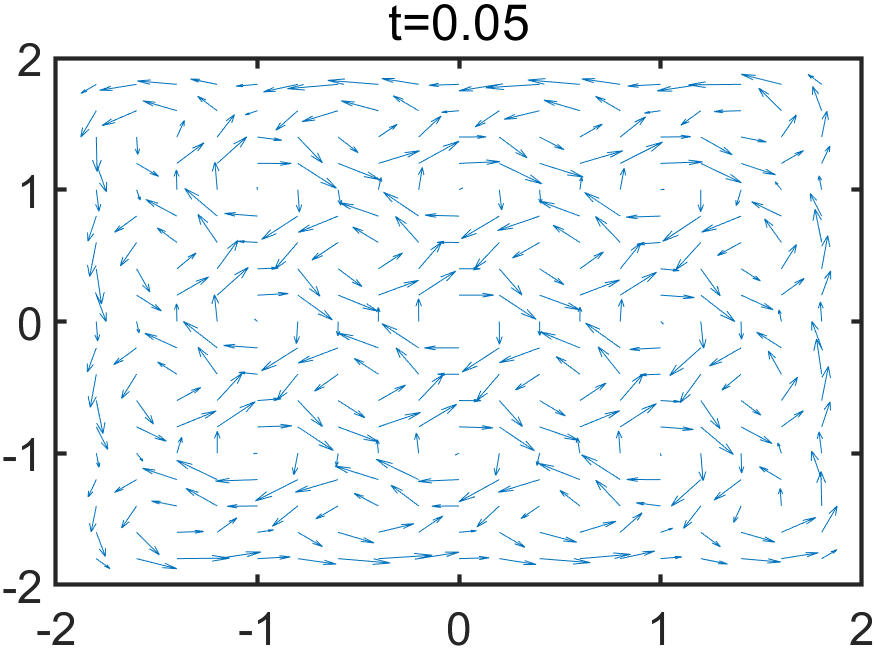}
\includegraphics[width=5.5cm,height=4.5cm]{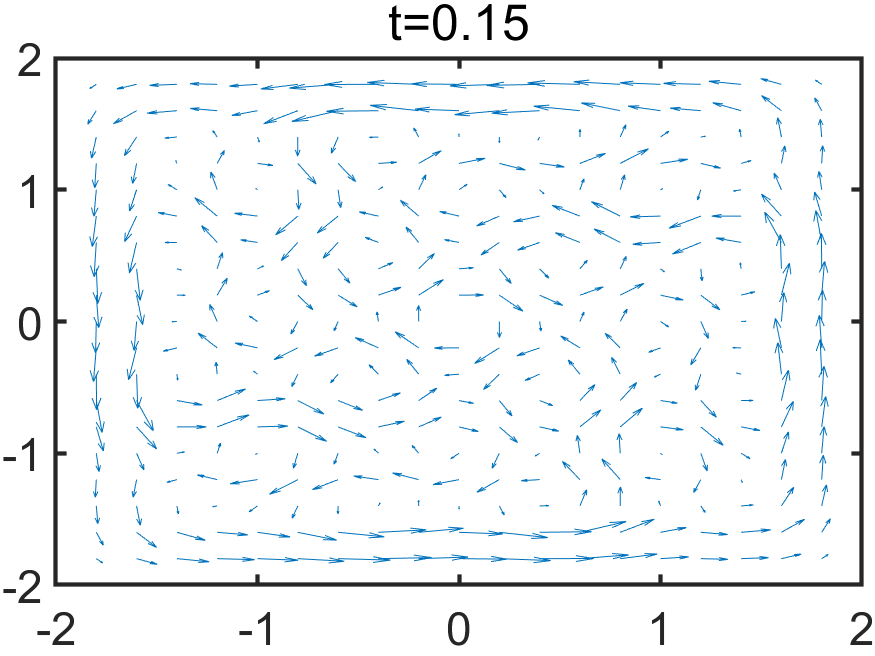}
\end{array}$\vspace{-0.2cm}
\caption{$\mathbf{Example\ \ref{exam3}}$, \textcolor{black}{P-BDF2} scheme \eqref{se3-18}-\eqref{se3-19}, snapshots of numerical solutions for \textcolor{black}{the} velocity field function.}\label{Cexp4u}
\end{figure}

\begin{figure}[!htbp]
$\begin{array}{c}
\includegraphics[width=8cm,height=6.5cm]{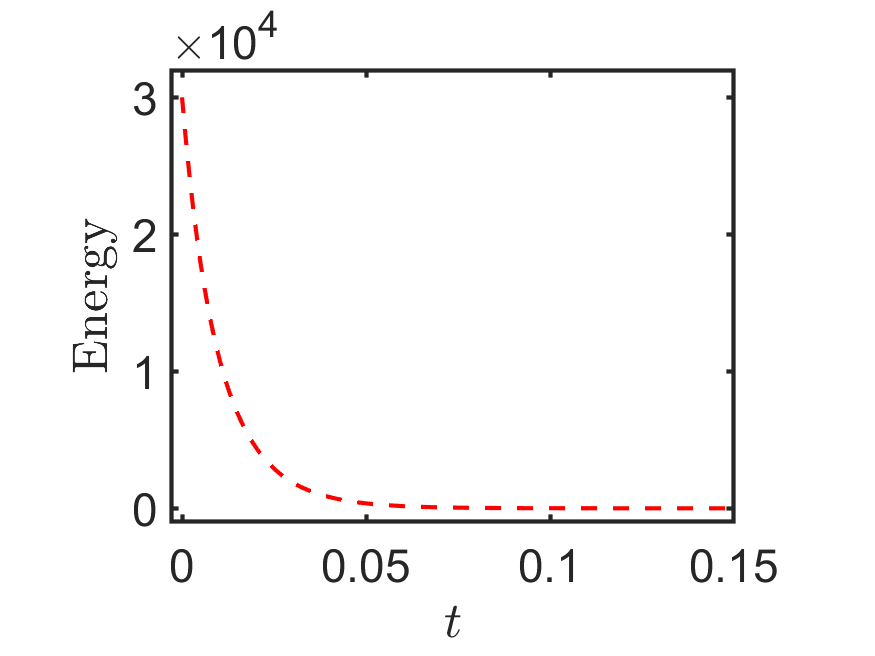}
\includegraphics[width=8cm,height=6.5cm]{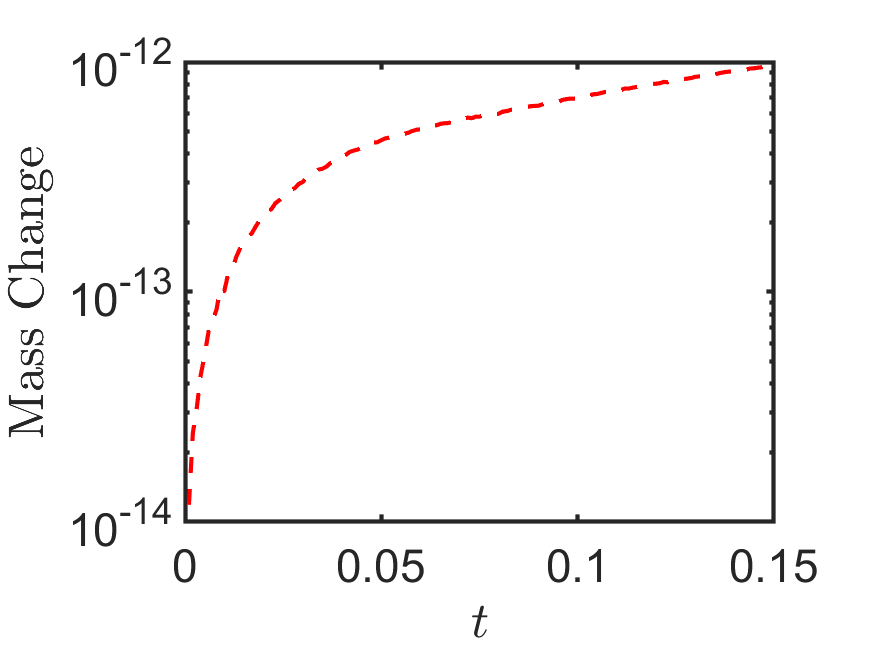}
\end{array}$
\caption{$\mathbf{Example\ \ref{exam3}}$, the modified discrete energy history and discrete mass \textcolor{black}{variation}.}\label{Cexp4 energy and mass}
\end{figure}

As for the $\mathbf{Example\ \ref{exam3}}$, the contour plots of the numerical solutions of $\phi_h^n$ and $\mathbf{u}_h^n$ by using the \textcolor{black}{P-BDF2} scheme are shown in Figures \ref{Cexp4phi}-\ref{Cexp4u}.
And the evolution of discrete energy and total change of discrete mass are shown in Figure \ref{Cexp4 energy and mass}. As it is shown, we can also see that the discrete energy decreases over time and the change of total mass reaches machine precision. 

\begin{example} \label{examlong}
\textcolor{black}{
The purpose of this example is to evaluate the long-time performance of the proposed method. Following \cite{CHEN201640}, we consider the CHNS equations \eqref{eq1.3} with the initial data
\begin{equation} \label{randinit}
\begin{aligned}
\phi_0 &= 2\text{rand}(x,y)-1,  \\
\mathbf{u}_0 &= \mathbf{0},  
\end{aligned}   
\end{equation}
where $\text{rand}$ denotes a two-dimensional random function.
The parameters are set as $\gamma = 0.002,\; \mu = 1\;, \lambda=0.01,\; \epsilon = 0.02$. The computational domain $\Omega=[-1,1] \times [-1,1]$ is uniformly triangulated using a mesh identified by partition number $N \times N = 160 \times 160$. The final simulation time is $T=100$ with a time step size $\tau = 10^{-2}$. 
We solve this problem by the \textcolor{black}{P-BDF1} scheme \eqref{se3-1}-\eqref{se3-3}, with the IEQ auxiliary constant $B=1$.}

\textcolor{black}{
The evolution of the phase function at various times is shown in \Cref{Cexplongphi}, clearly illustrating the coarsening dynamics. The accompanying plots of energy and mass variations in \Cref{CexplongEandM} further confirm that the numerical results preserve the expected physical properties.
}

\begin{figure}[!htbp]
$\begin{array}{c}
\includegraphics[width=5.5cm,height=4.5cm]{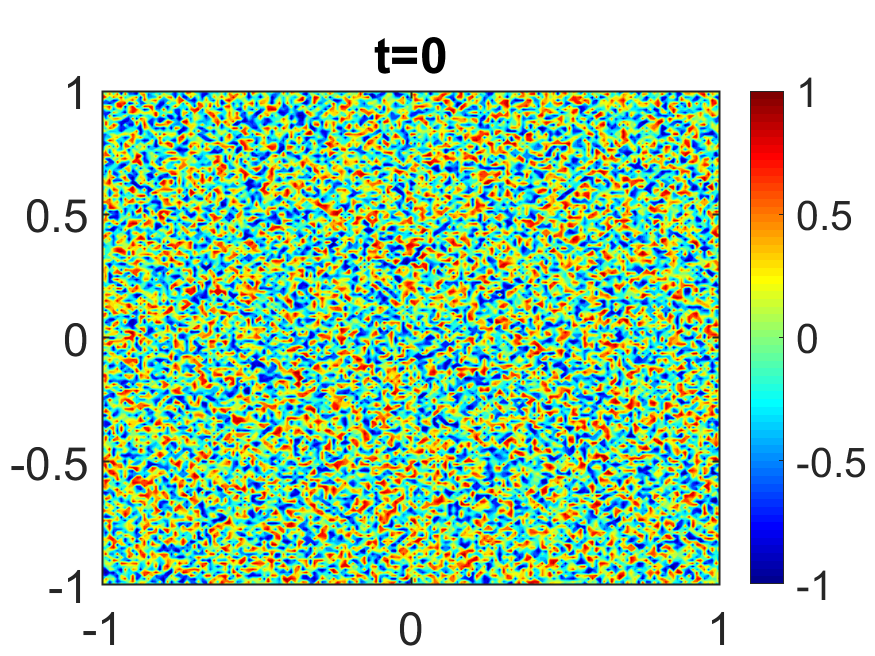} \;
\includegraphics[width=5.5cm,height=4.5cm]{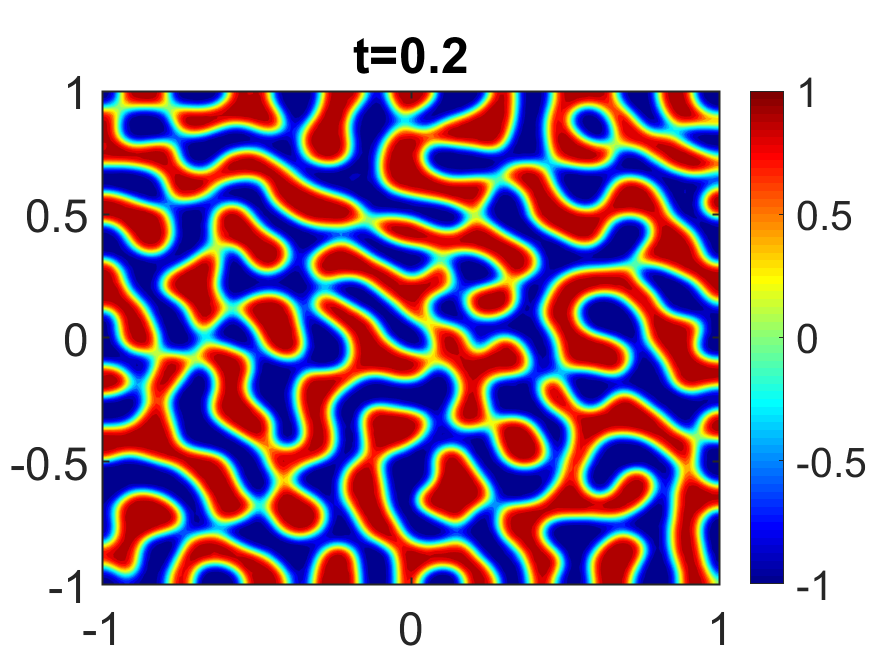} \;
\includegraphics[width=5.5cm,height=4.5cm]{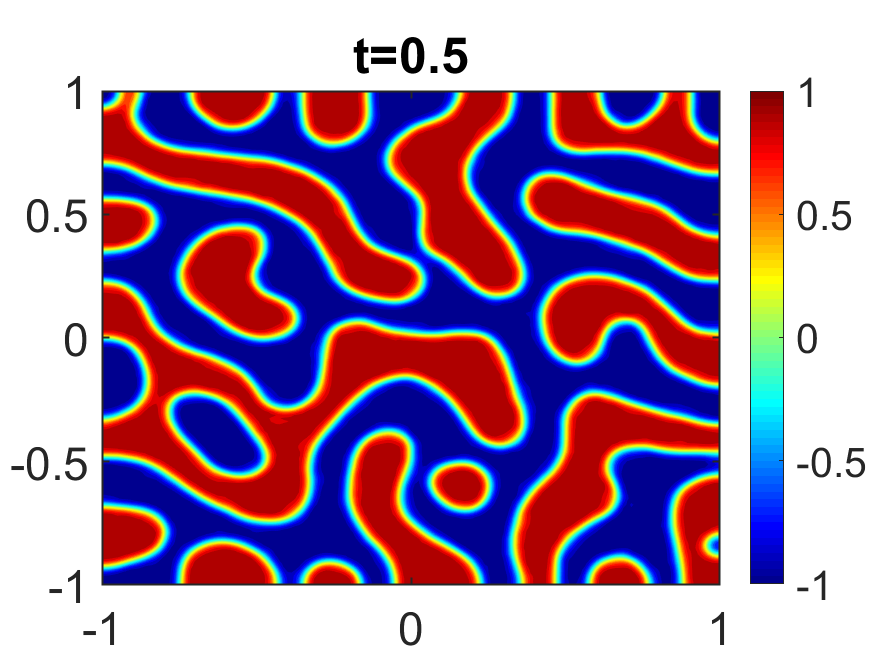} \\
\includegraphics[width=5.5cm,height=4.5cm]{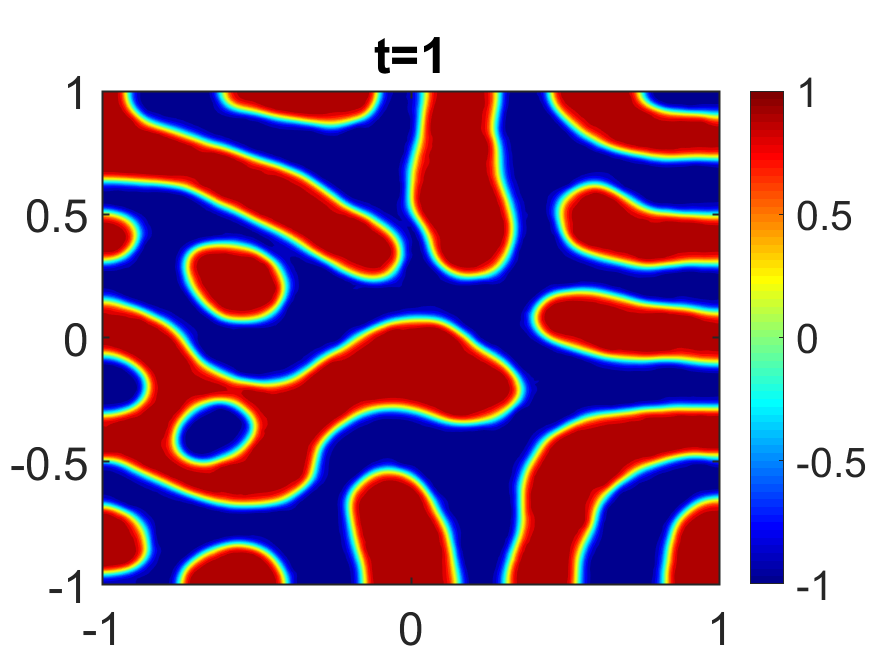}\; 
\includegraphics[width=5.5cm,height=4.5cm]{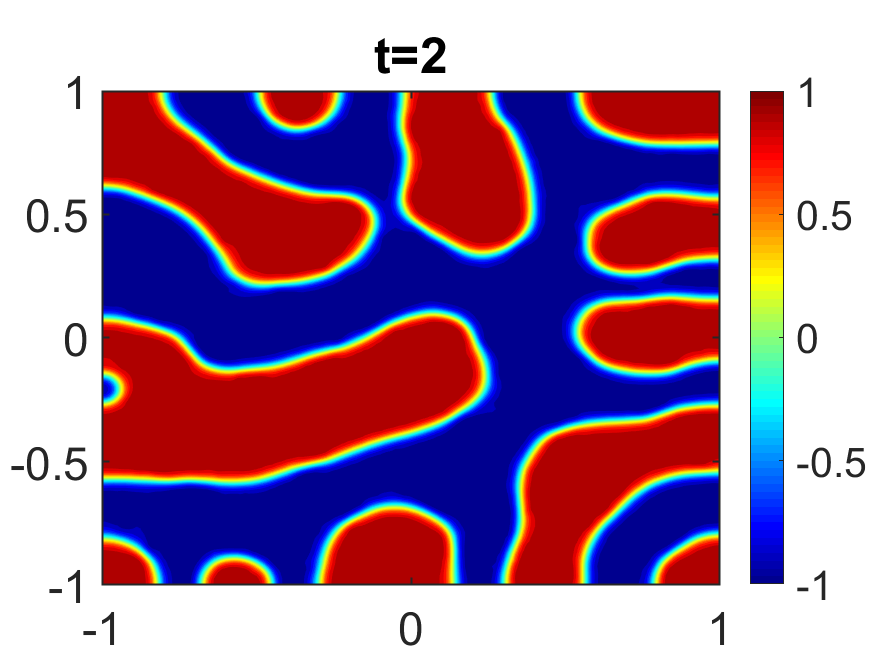}\;
\includegraphics[width=5.5cm,height=4.5cm]{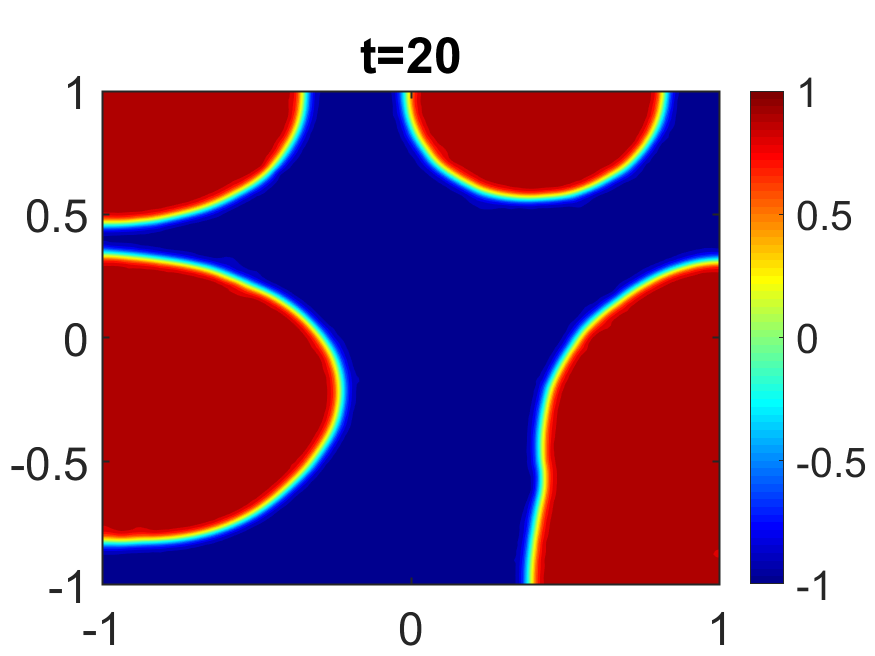} \\
\includegraphics[width=5.5cm,height=4.5cm]{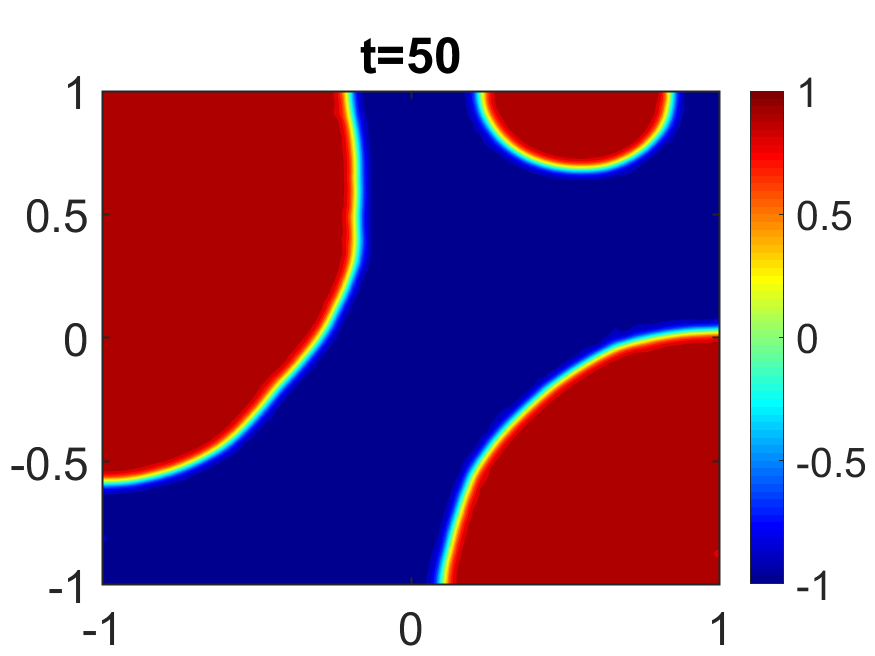}\; 
\includegraphics[width=5.5cm,height=4.5cm]{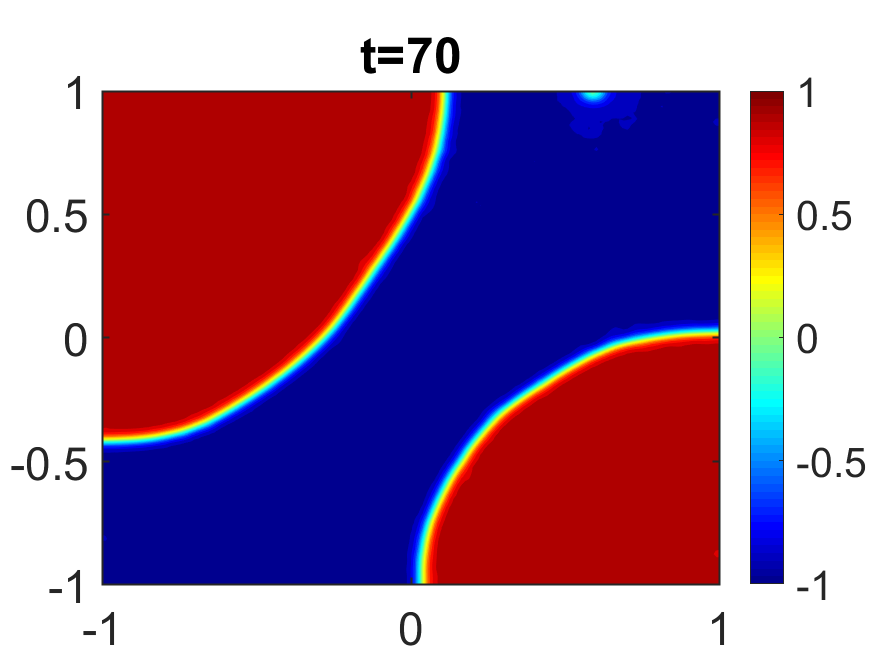}\;
\includegraphics[width=5.5cm,height=4.5cm]{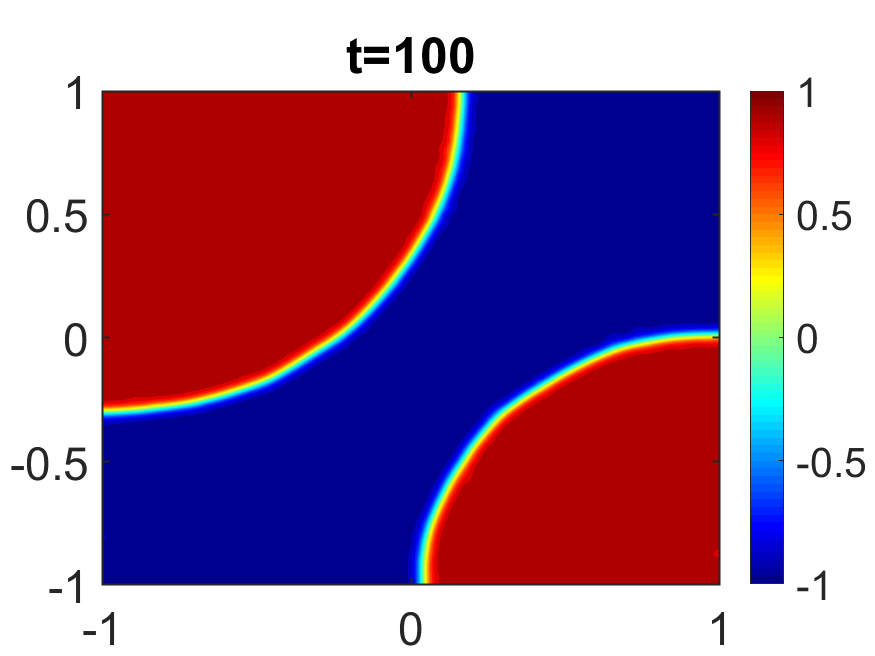}
\end{array}$\vspace{-0.2cm}
\caption{$\mathbf{Example\ \ref{examlong}}$, \textcolor{black}{P-BDF1} scheme \eqref{se3-1}-\eqref{se3-3}, snapshots of numerical solutions for \textcolor{black}{the} phase field function.}\label{Cexplongphi}
\end{figure}

\begin{figure}[!htbp]
$\begin{array}{c}
\includegraphics[width=6.5cm,height=6.5cm]{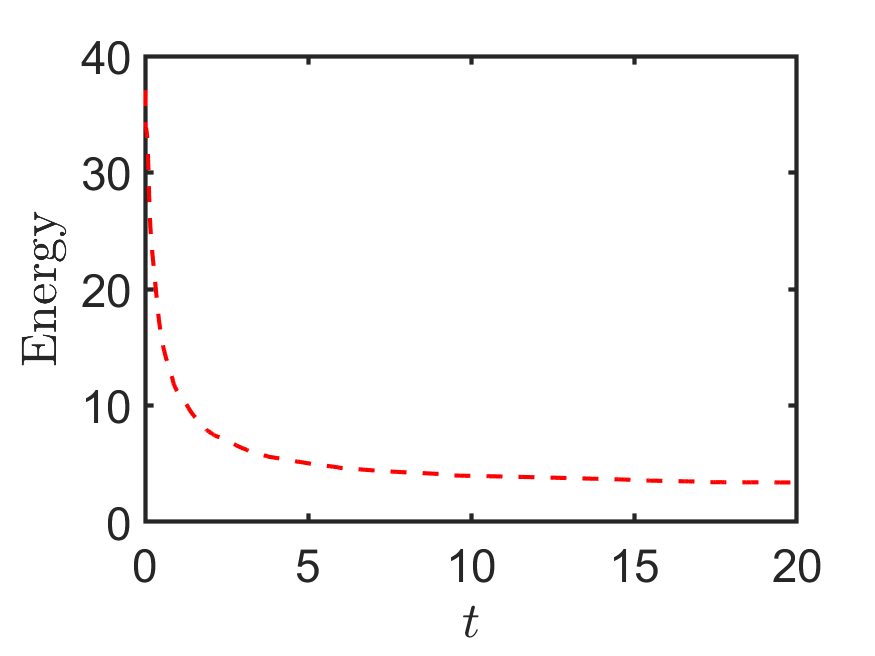} \;
\includegraphics[width=7.3cm,height=6.2cm]{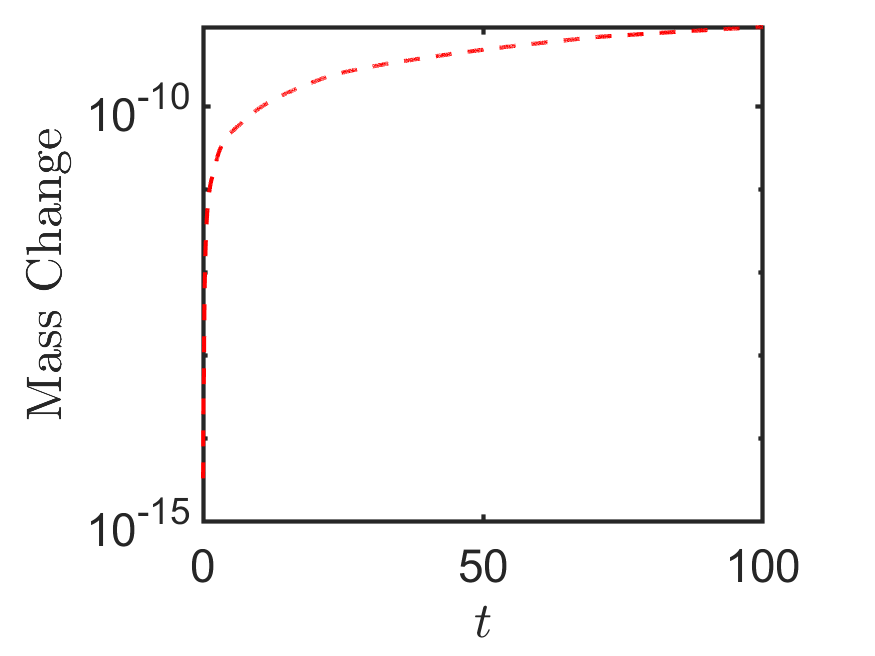} \;
\end{array}$\vspace{-0.2cm}
\caption{$\mathbf{Example\ \ref{examlong}}$, the discrete energy history and discrete mass variation.}\label{CexplongEandM}
\end{figure}

\end{example}
\begin{example} \label{examdisc}
\textcolor{black}{In this example, we demonstrate the flexibility of the proposed method on a curved domain. The same parameters as in $\mathbf{Example\ \ref{examlong}}$ are used for equation \eqref{eq1.3}, along with the same final simulation time, time-step size, and IEQ auxiliary constant $B$. The initial condition is given by \eqref{randinit}, and the computational domain is defined as $\Omega={(x,y):x^2+y^2\le 1}$. A quasi-uniform triangular mesh with $21{,}396$ nodes is employed for discretization.}

\textcolor{black}{Snapshots of the numerical phase function at selected time instances are shown in \Cref{Cexpdiscphi}. The energy and mass evolutions are shown in \Cref{CexpdiscEandM}, demonstrating both energy dissipation and mass conservation. These results confirm that the proposed scheme performs effectively on curved domains as well.}
\begin{figure}[!htbp]
$\begin{array}{c}
\includegraphics[width=5.5cm,height=4.5cm]{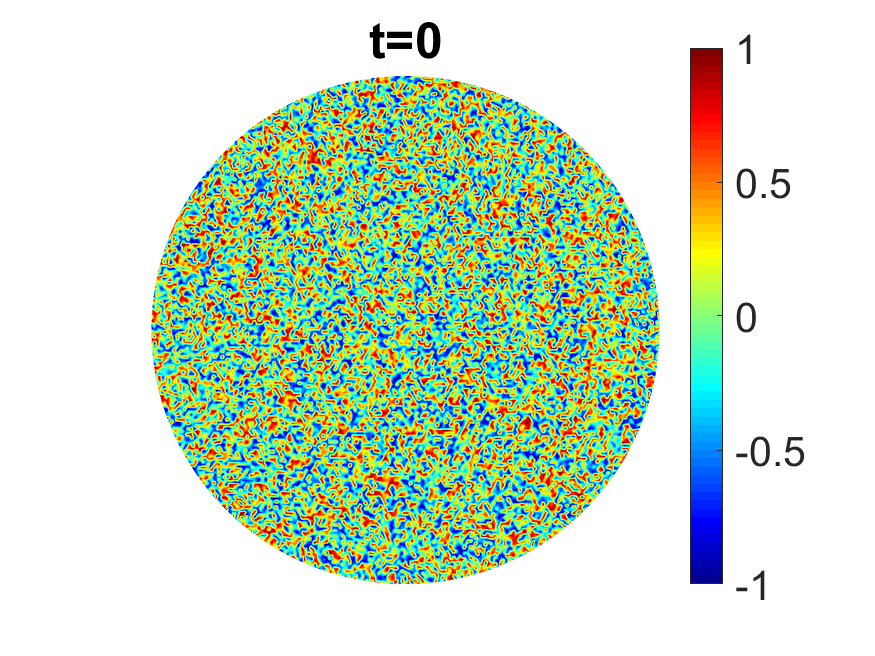} \;
\includegraphics[width=5.5cm,height=4.5cm]{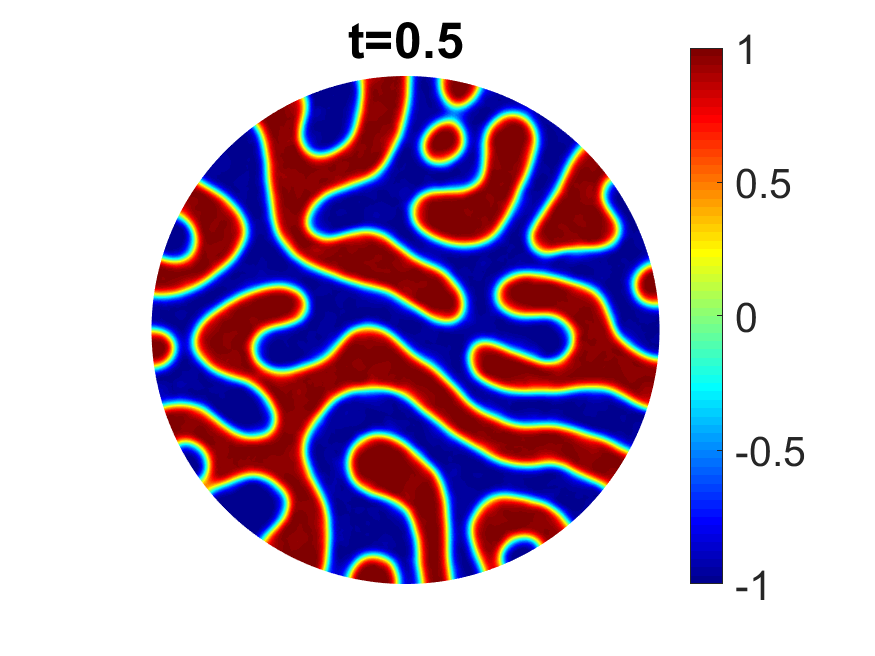} \;
\includegraphics[width=5.5cm,height=4.5cm]{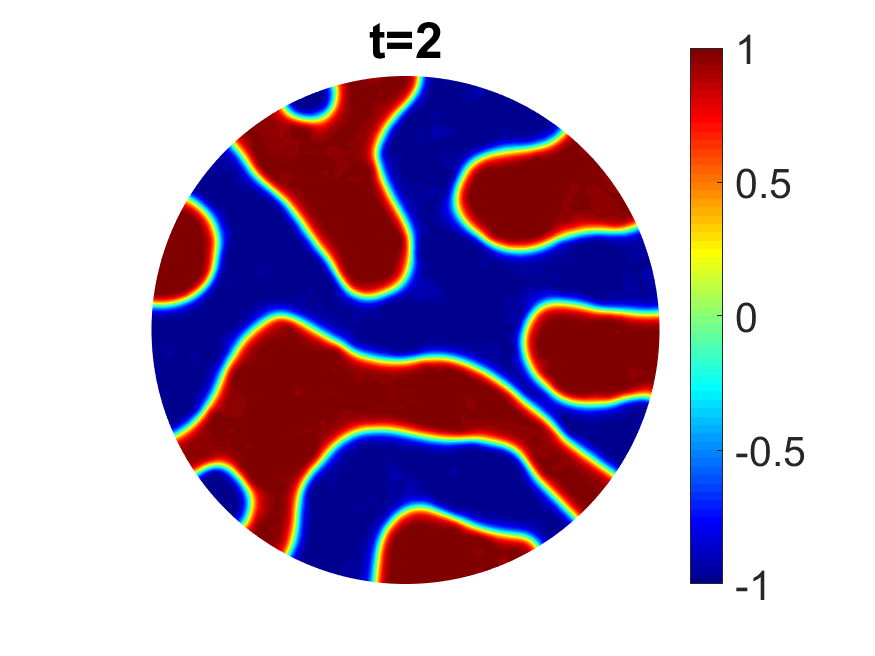} \\
\includegraphics[width=5.5cm,height=4.5cm]{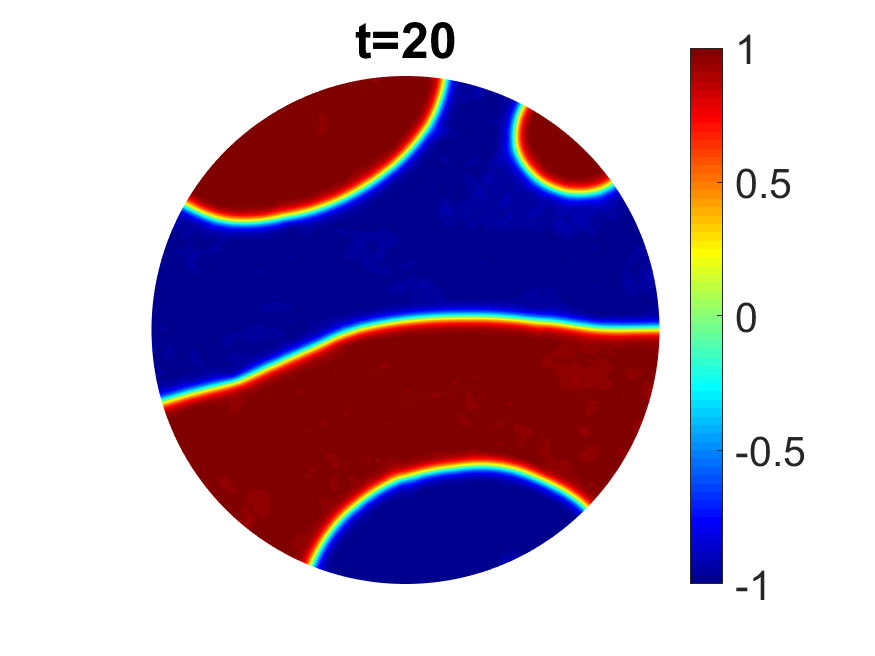}\; 
\includegraphics[width=5.5cm,height=4.5cm]{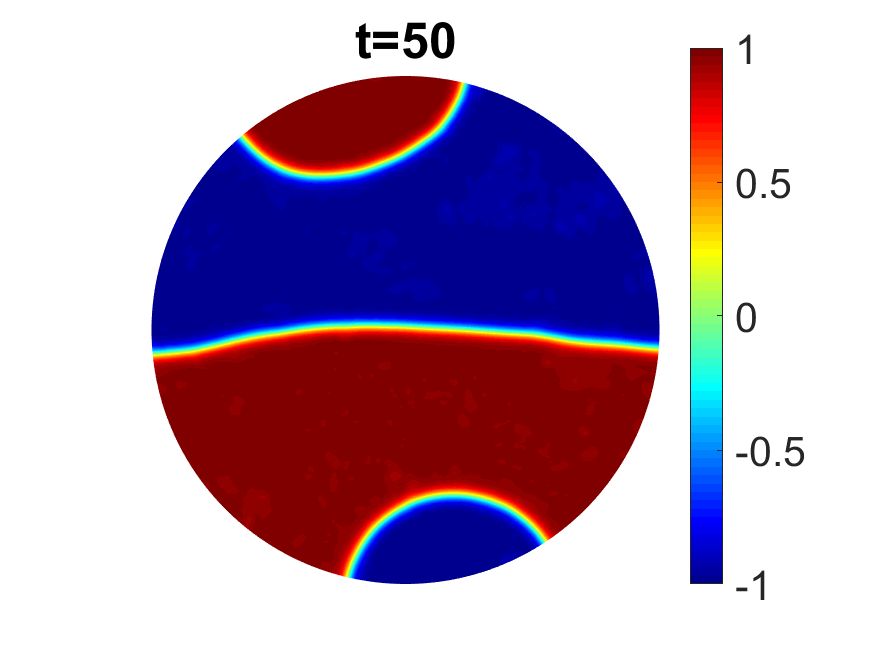}\;
\includegraphics[width=5.5cm,height=4.5cm]{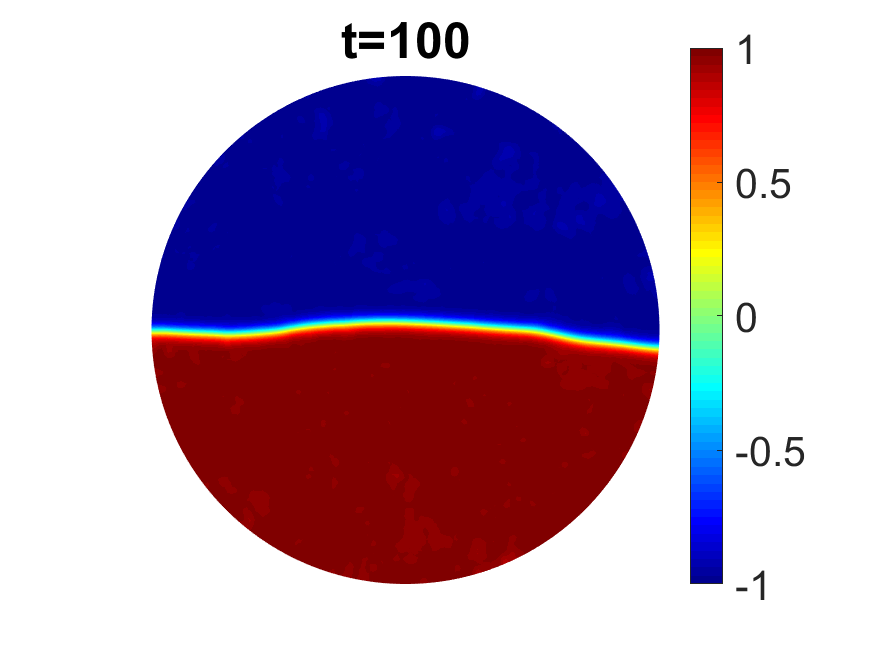}
\end{array}$\vspace{-0.2cm}
\caption{$\mathbf{Example\ \ref{examdisc}}$, \textcolor{black}{P-BDF1} scheme \eqref{se3-1}-\eqref{se3-3}, snapshots of numerical solutions for \textcolor{black}{the} phase field function.}\label{Cexpdiscphi}
\end{figure}
\begin{figure}[!htbp]
$\begin{array}{c}
\includegraphics[width=6.5cm,height=6.5cm]{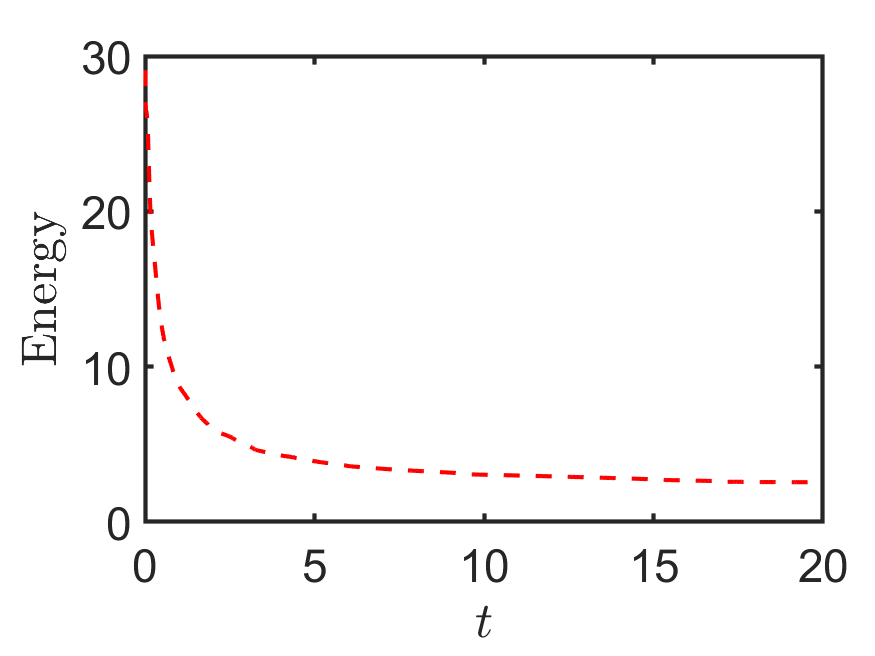} \;
\includegraphics[width=7.3cm,height=6.2cm]{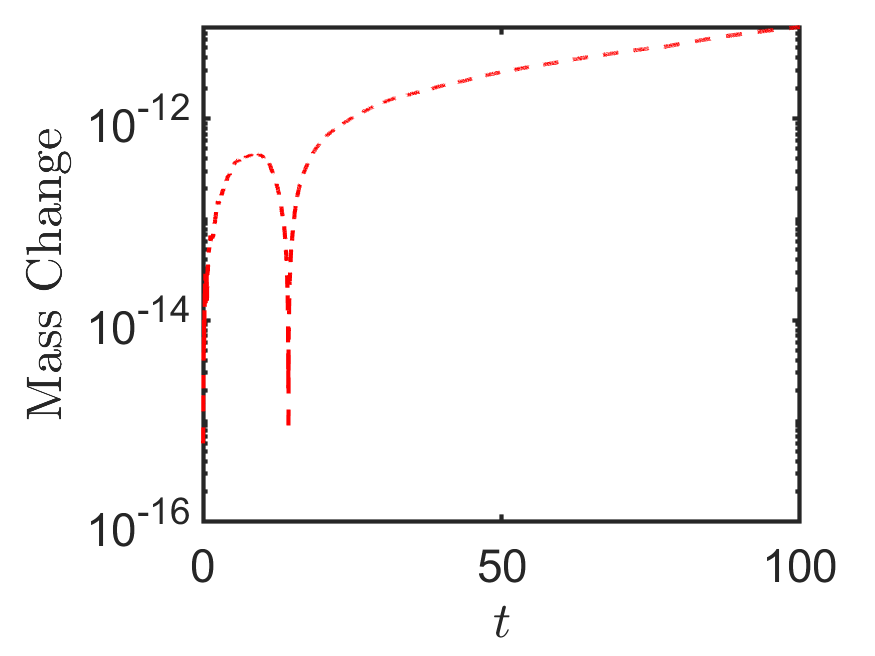} \;
\end{array}$\vspace{-0.2cm}
\caption{$\mathbf{Example\ \ref{examdisc}}$, the discrete energy history and discrete mass variation.}\label{CexpdiscEandM}
\end{figure}
\end{example}
\begin{example} \label{exam5}\cite{chen_IEQ_FEM}
Let $\Omega = [-\frac{1}{2},-\frac{1}{2}] \times [-\frac{1}{5},-\frac{1}{5}]$, we consider the CHNS equations satisfying the following initial condition 
\begin{equation}
\begin{aligned}
& \phi_0 = \left\{ \begin{aligned}
1, & \qquad x<x_0, \\
-1, & \qquad x>x_1, \\
-\sin \left( \frac{\pi x}{2 x_1} \right), & \qquad x_0 \le x \le x_1,
\end{aligned} \right. \\
& \mathbf{u}_0 = [0 \quad 0]^\top,
\end{aligned}
\end{equation}
where $x_1=x_0=\frac{\sqrt{2}}{20}$, $\epsilon=\frac{1}{500\sqrt{10}},\; \lambda=\frac{1}{100},\;\gamma=\frac{1}{10},\;\mu=1$, $\tau = 10^{-7},\; T=10^{-5}$.

We solve this problem using the schemes \textcolor{black}{P-BDF1}, \textcolor{black}{C-BDF1}, and \textcolor{black}{CP-BDF1} with different \textcolor{black} {IEQ constants $B$, which is introduced in the IEQ approach to ensure $F(\phi_h)+B> 0$.}
The evolution of the discrete energy of $E\left( \phi_h^{n+1},{\mathbf{u}}_h^{n+1},U_h^{n+1},p_h^{n+1} \right)$ for \textcolor{black}{P-BDF1} scheme and $E\left( \phi_h^{n+1},{\mathbf{u}}_h^{n+1},U_I^{n+1},p_h^{n+1} \right)$ for \textcolor{black}{C-BDF1} scheme are shown in Figure \ref{energy apdx}. From the results, we observe that solutions of the \textcolor{black}{C-BDF1} scheme satisfy the energy dissipation law only when $B$
is sufficiently large, as explained in Remark \ref{Enerdis}, whereas the \textcolor{black}{P-BDF1} scheme satisfies the energy dissipation law regardless of the choice of $B$. There are noticeable discrepancies between the energy curves produced by the \textcolor{black}{C-BDF1} and \textcolor{black}{P-BDF1} schemes when the mesh size is large ($h = \frac{1}{10}$, as seen in Figure \ref{energy apdx}). However, these differences decrease as $h$ decreases, as illustrated in Figure \ref{energy apdx PC}. 

To demonstrate how the \textcolor{black}{CP-BDF1} scheme improves the \textcolor{black}{C-BDF1} scheme in preserving the energy dissipation law, we present energy plots from all three schemes for comparison in Figure \ref{energy apdx three}. To capture the energy increase phenomenon in this specific example, the initial computation is performed on a coarse grid with $h = \frac{1}{10}$. At this stage, we observe a significant discrepancy between the energy curves obtained from the \textcolor{black}{C-BDF1} and \textcolor{black}{P-BDF1} schemes. When the \textcolor{black}{CP-BDF1} scheme is applied for correction, a noticeable energy decrease occurs at $t = 9.6 \times 10^{-6}$. As the grid refined to $h = \frac{1}{160}$, the magnitude of this change diminished. The energy curve produced by the \textcolor{black}{CP-BDF1} scheme is aligned with that of the \textcolor{black}{P-BDF1} scheme. Moreover, the difference between the energy curves from the \textcolor{black}{P-BDF1} scheme on coarse and fine meshes is smaller than that of the \textcolor{black}{C-BDF1} scheme. 
It implies that the modification introduced by the \textcolor{black}{CP-BDF1} scheme is effective in ensuring energy dissipation.

Pictures in the top row of \Cref{phi apdx three} display the phase field plots at $T = 10^{-5}$, computed by using three distinct schemes. 
slight differences are observed in both the energy curves (Figure \ref{energy apdx three}, left) and the phase field plots (Figure \ref{phi apdx three}, top) on the coarser grid ($h = \frac{1}{10}$). 
These discrepancies are primarily attributed to the significant projection error in the intermediate variable $U^{n}$. As expected, this error diminishes with increased grid resolution ($h = \frac{1}{160}$; Figure \ref{energy apdx three}, right and Figure \ref{phi apdx three}, bottom), but we can still observe increasing energy for the \textcolor{black}{C-BDF1}
scheme. However, the \textcolor{black}{P-BDF1} and \textcolor{black}{CP-BDF1} schemes satisfy the energy dissipation law.

\begin{figure}[!htbp]
$\begin{array}{c} \hspace{-0.5cm}
\includegraphics[width=8cm,height=6.5cm]{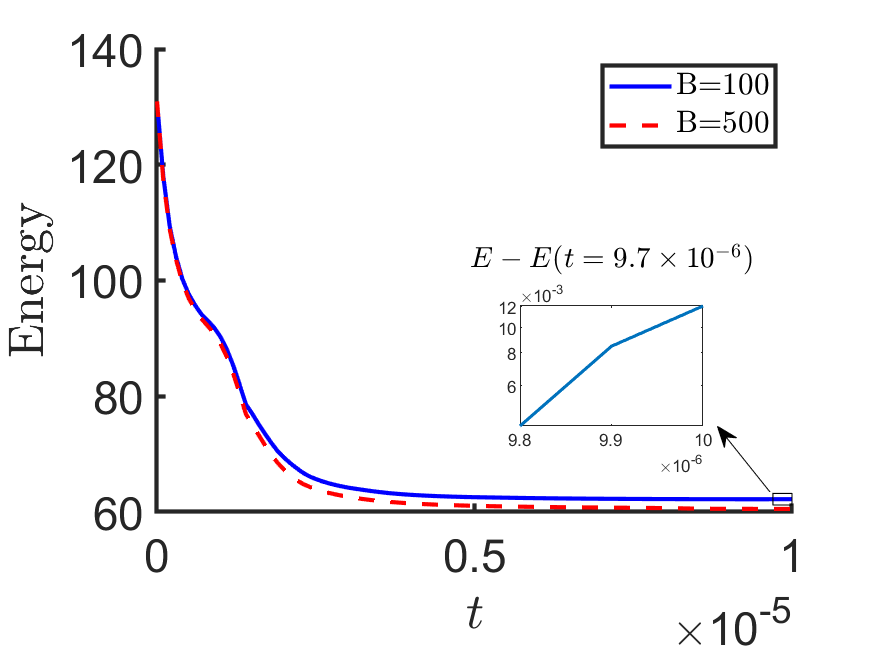}\;
\includegraphics[width=8cm,height=6.5cm]{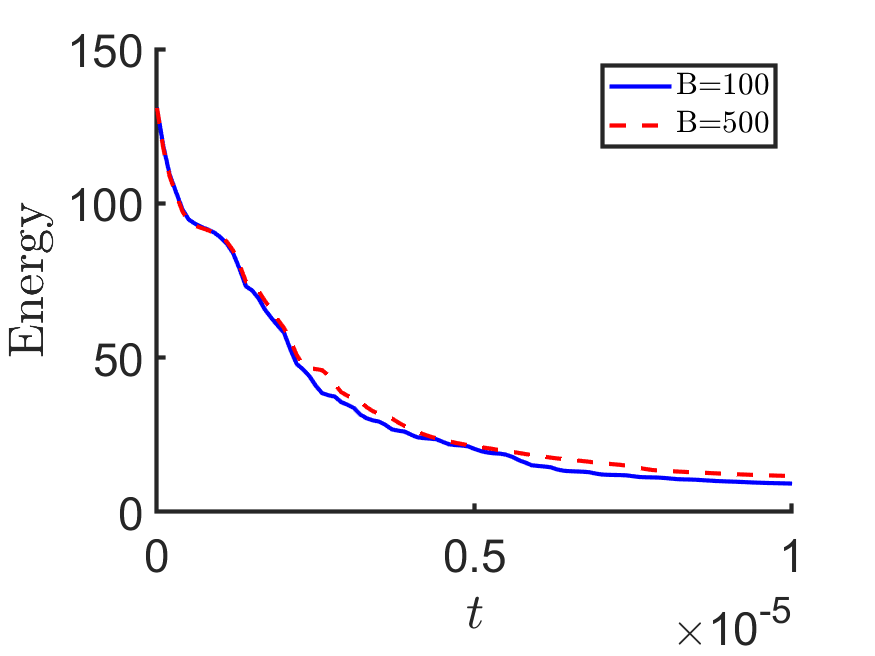}
\end{array}$\vspace{-0.2cm}
\caption{$\mathbf{Example\ \ref{exam5}}$,  Effect of constant B on the discrete energy,  Left: \textcolor{black}{C-BDF1} scheme; Right: \textcolor{black}{P-BDF1} scheme, $h=\frac{1}{10}$.}\label{energy apdx}
\end{figure}

\begin{figure}[!htbp]
$\begin{array}{c} \hspace{-0.5cm}
\includegraphics[width=8cm,height=6.5cm]{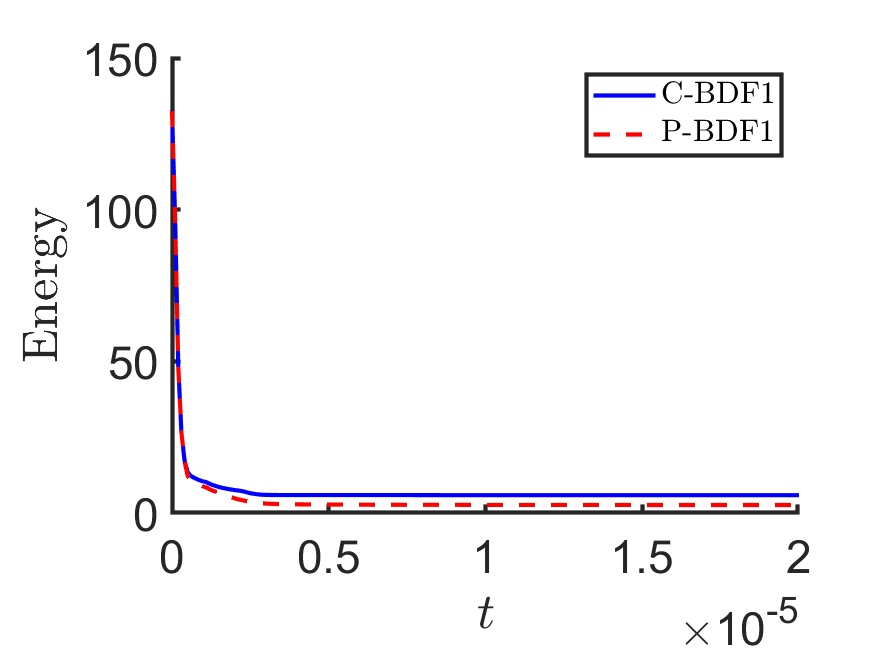}\;
\includegraphics[width=8cm,height=6.5cm]{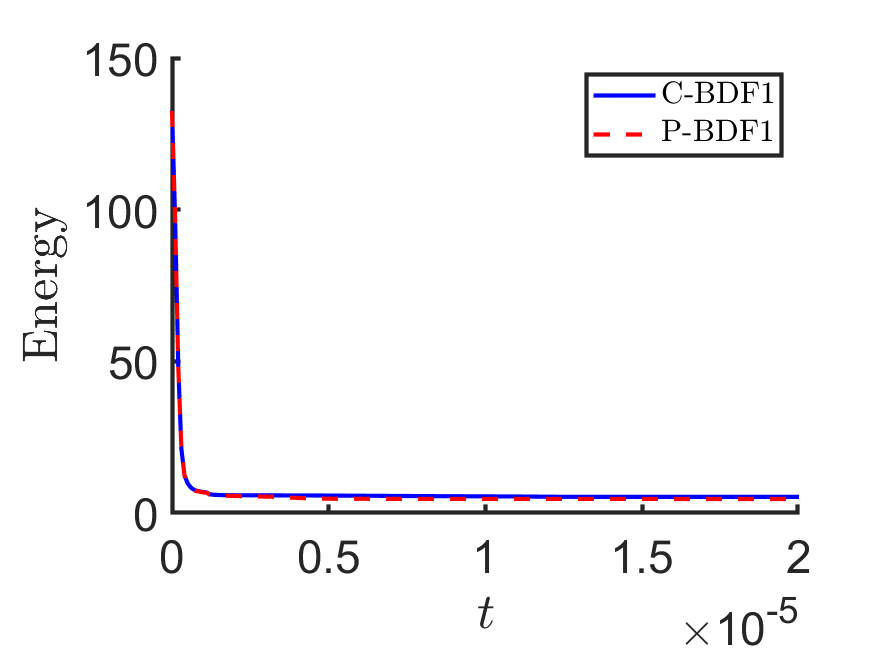}
\end{array}$\vspace{-0.2cm}
\caption{$\mathbf{Example\ \ref{exam5}}$,  Effect of mesh size on the difference of discrete energy, $B=500$. Left: $h = \frac{1}{160}$; Right: $h = \frac{1}{320}$.}\label{energy apdx PC}
\end{figure}

\begin{figure}[!htbp]
$\begin{array}{c} \hspace{-0.5cm} \includegraphics[width=8cm,height=6.5cm]{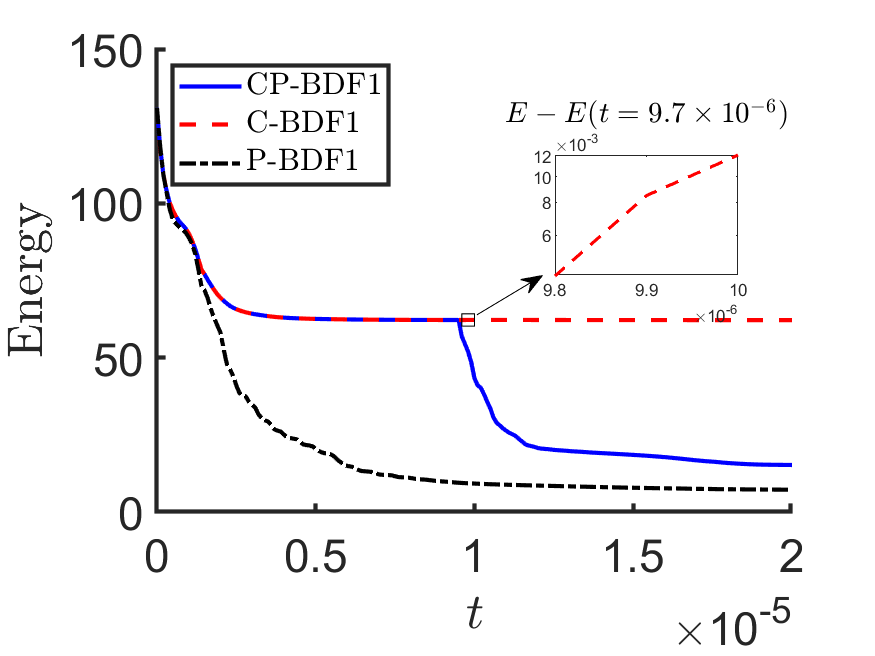}\;
\includegraphics[width=8cm,height=6.5cm]{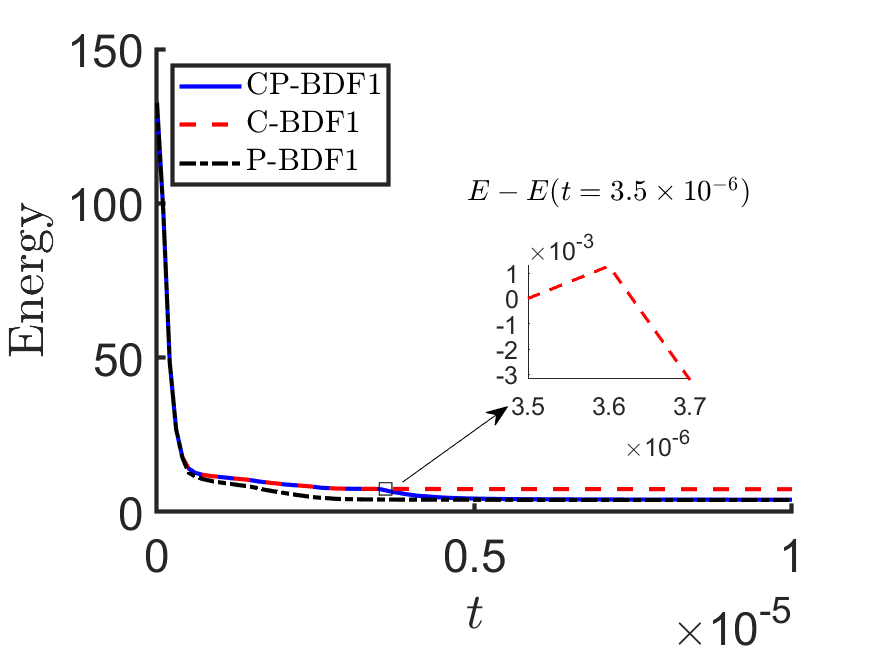}\;
\end{array}$\vspace{-0.2cm}
\caption{$\mathbf{Example\ \ref{exam5}}$, The comparison of energy curves for three methods. $B = 100$. Left: $h=\frac{1}{10}$, right: $h=\frac{1}{160}$}\label{energy apdx three}
\end{figure}

\begin{figure}[!htbp]
$\begin{array}{c} \hspace{-0.5cm}
\includegraphics[width=5.5cm,height=4.5cm]{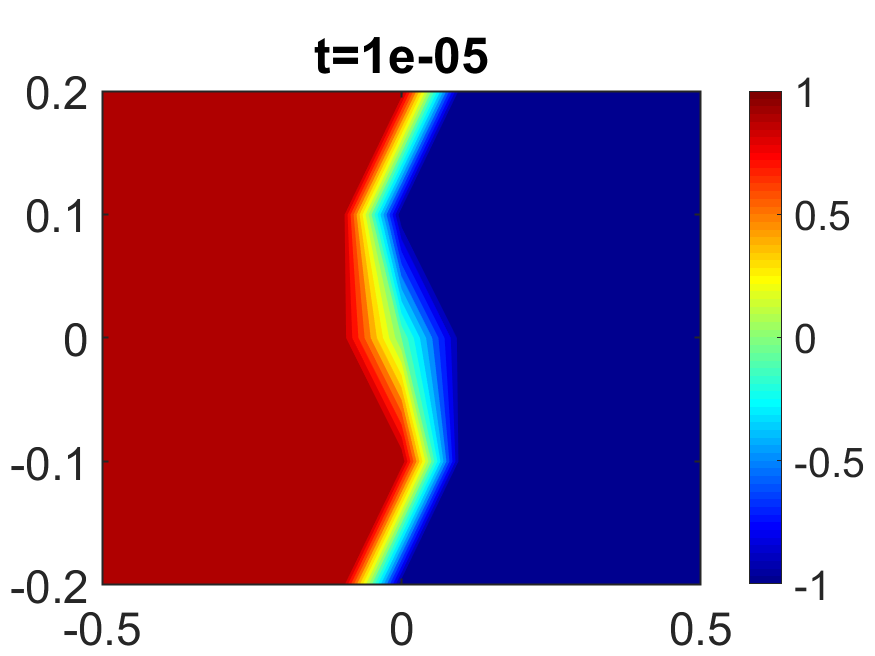}\;\includegraphics[width=5.5cm,height=4.5cm]{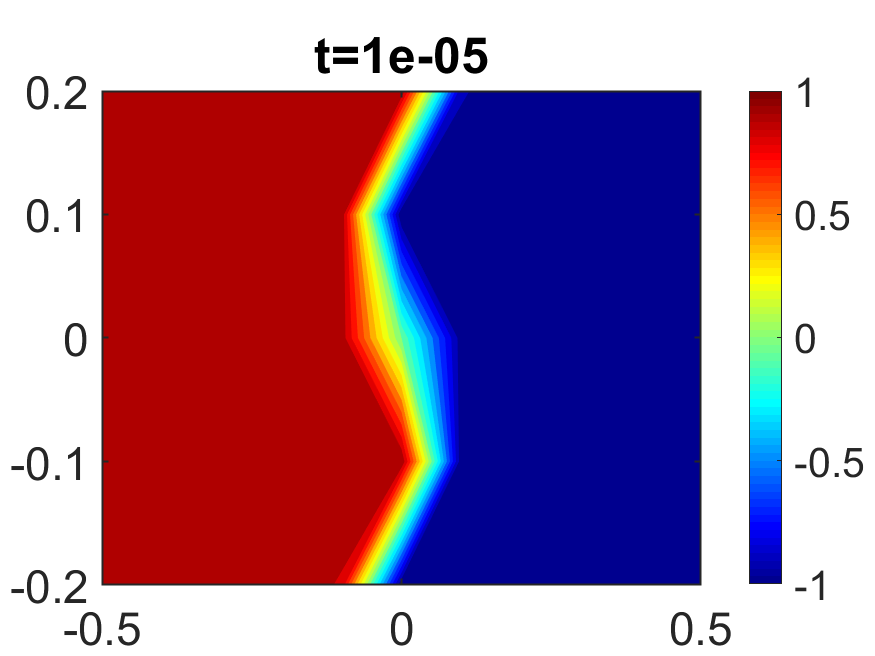}\;\includegraphics[width=5.5cm,height=4.5cm]{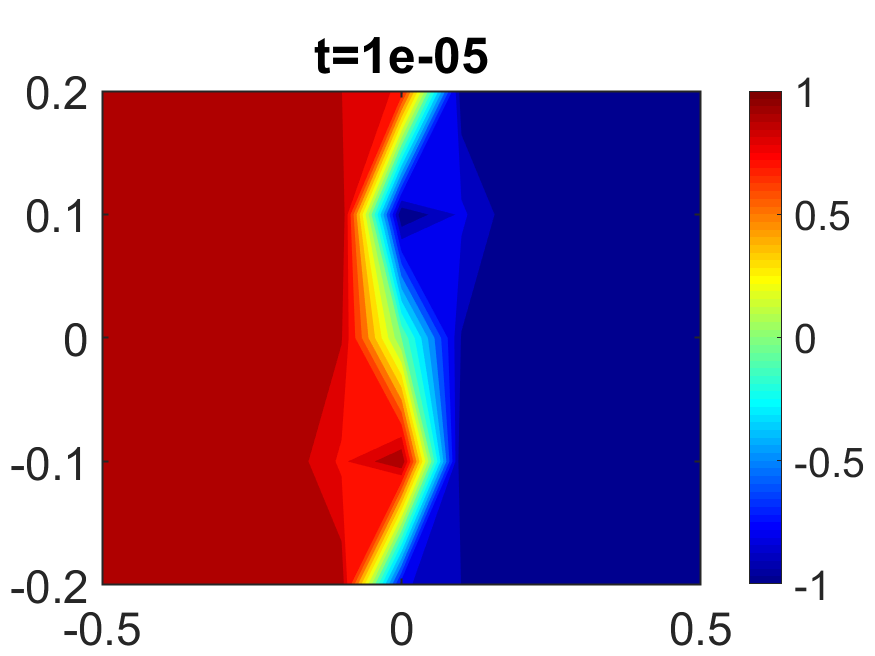} \\ \hspace{-0.5cm}
\includegraphics[width=5.5cm,height=4.5cm]{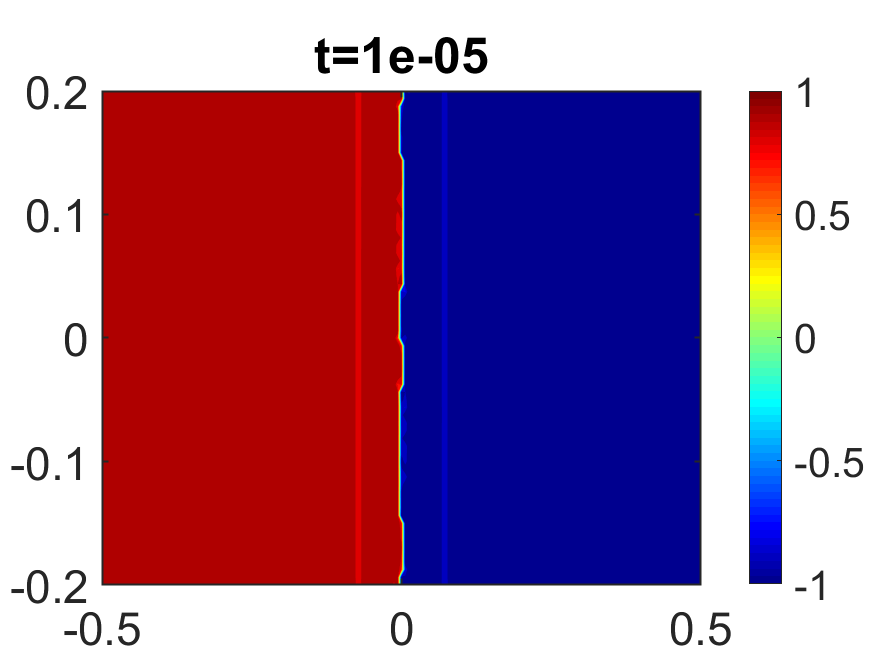}\;\includegraphics[width=5.5cm,height=4.5cm]{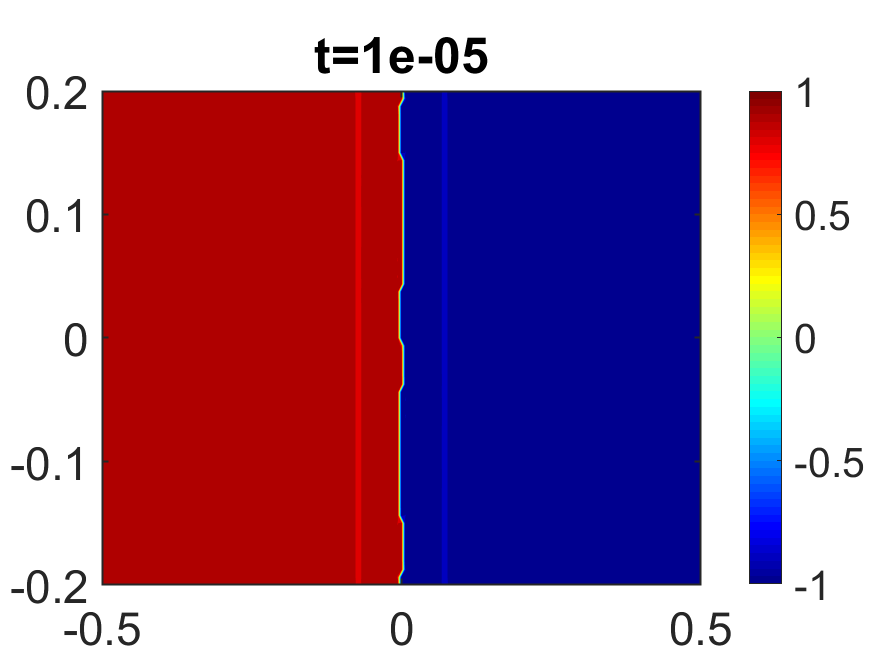}\;\includegraphics[width=5.5cm,height=4.5cm]{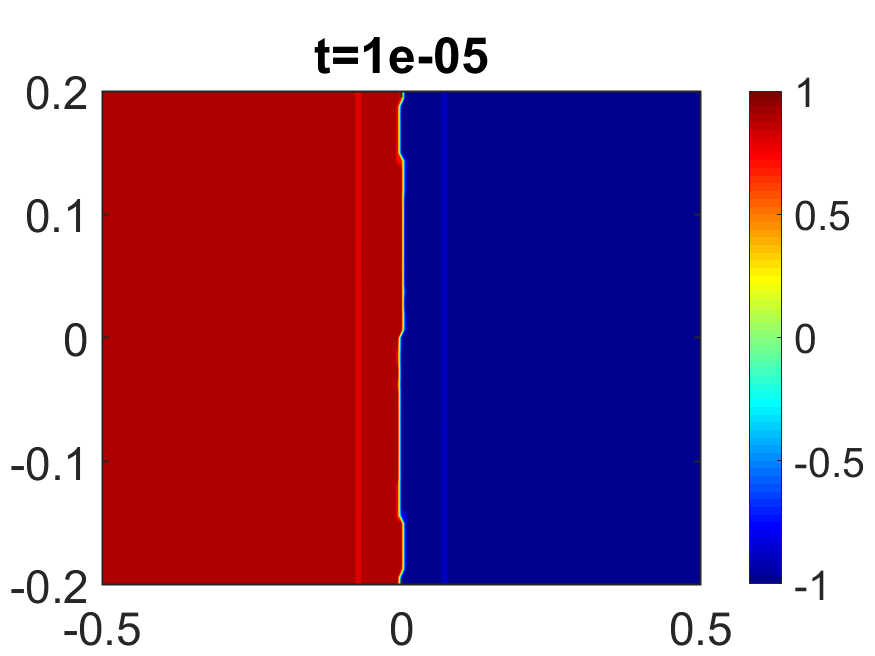}
\end{array}$\vspace{-0.2cm}
\caption{$\mathbf{Example\ \ref{exam5}}$, The comparison of energy curves for three methods. Top: $h=\frac{1}{10}$, 
bottom: $h=\frac{1}{160}$; left: \textcolor{black}{CP-BDF1}, center: \textcolor{black}{C-BDF1 scheme}, right: \textcolor{black}{P-BDF1 scheme}.}\label{phi apdx three}
\end{figure}

\end{example}

\begin{example} \label{examThreeD}
\textcolor{black}{
We design the following three-dimensional test case to assess the capability of the proposed method in handling 3D models. Consider the CHNS equations \eqref{eq1.3} imposed on the domain $\Omega = [0.2, 0.8] \times [0.2, 0.8] \times [0, 1]$, with initial conditions specified as
\begin{equation}
\begin{aligned}
\phi_0 &= 1 - \tanh\left({\frac{-r+\sqrt{(x-x_a)^2+(y-y_a)^2+(z-z_a)^2}}{2\epsilon}}\right) \\ & \quad - \tanh\left({\frac{-r+\sqrt{(x-x_b)^2+(y-y_b)^2+(z-z_b)^2}}{2\epsilon}}\right) ,\\
\mathbf{u}_0 &= \left[ 0, \; 0, \;0 \right]^\top , 
\end{aligned}
\end{equation}
where $x_a=y_a=x_b=y_b=0.5,\;z_a=0.65,\;z_b=0.35$, and $r=0.15$. 
Following \cite{YiShi2013Inverse}, we choose $\gamma = 0.001$, $\mu = 12$, $\epsilon = 0.02$, $\lambda = 0.05$, and $B = 50$. The problem is solved using the P-BDF1 scheme \eqref{se3-1}-\eqref{se3-3} with uniform mesh size $h = 1/50$ in all spatial directions and a time step $\tau = 5 \times 10^{-3}$, up to the final time $T = 1$.}

\textcolor{black}{
Figure~\ref{CThreeDphi} depicts the coalescence of two interacting bubbles in three dimensions. As time evolves, the bubbles gradually merge into a single, larger bubble driven by surface tension effects, in line with physical expectations \cite{YiShi2013Inverse}. The corresponding energy dissipation and mass variation over time are shown in Figure~\ref{Cexp3DEandM}, demonstrating consistency with the theoretical predictions.}

\begin{figure}[!htbp]
$\begin{array}{c} \hspace{-0.5cm}
\includegraphics[width=5cm,height=6cm]{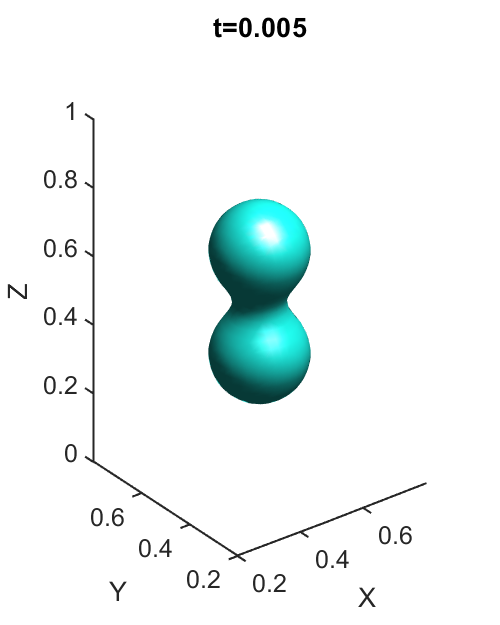}\;
\includegraphics[width=5cm,height=6cm]{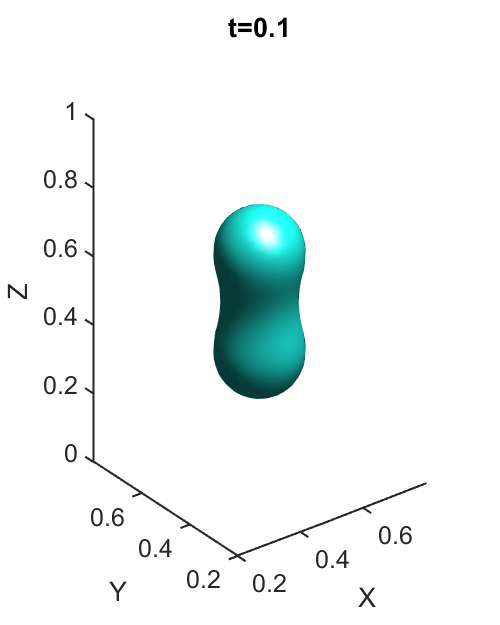}\;
\includegraphics[width=5cm,height=6cm]{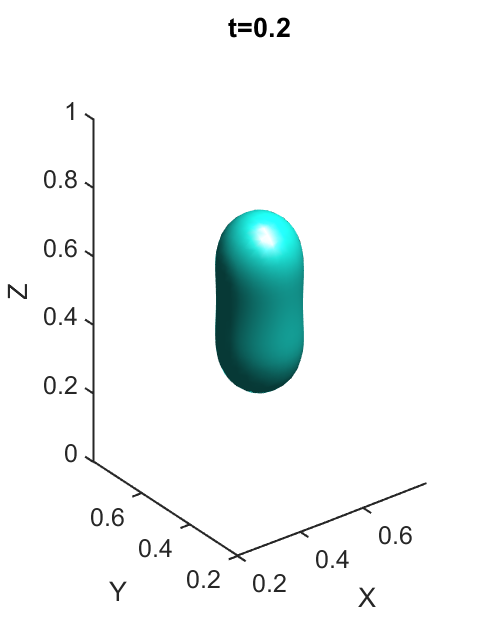}\\
\\ \hspace{-0.5cm}
\includegraphics[width=5cm,height=6cm]{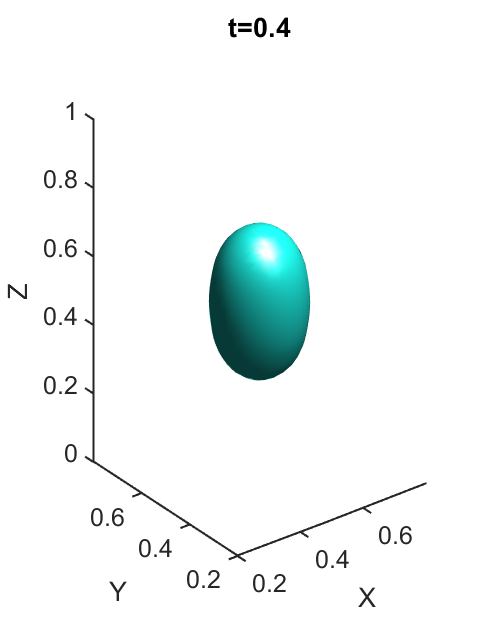}\;
\includegraphics[width=5cm,height=6cm]{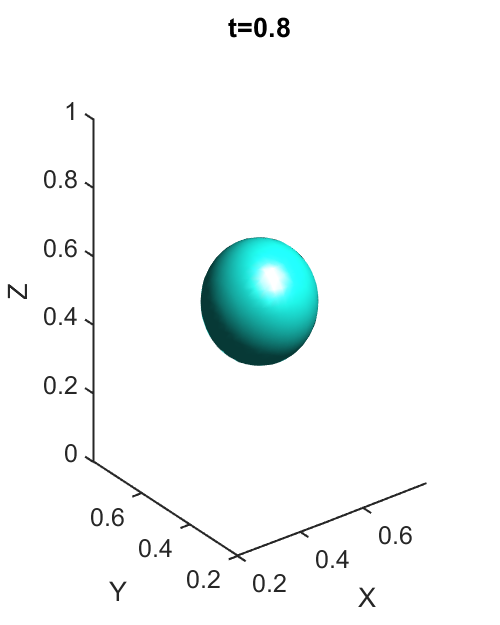}\;
\includegraphics[width=5cm,height=6cm]{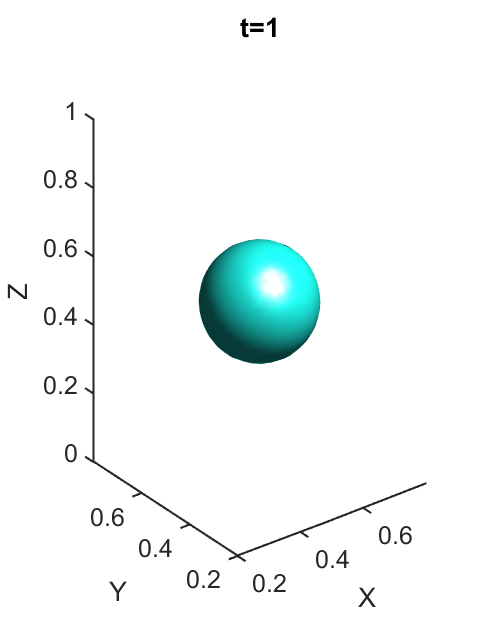}
\end{array}$\vspace{-0.2cm}
\caption{$\mathbf{Example\ \ref{examThreeD}}$, \textcolor{black}{P-BDF1} scheme \eqref{se3-1}-\eqref{se3-3}, snapshots of numerical solutions for \textcolor{black}{the} phase field function.}\label{CThreeDphi}
\end{figure}

\begin{figure}[!htbp]
$\begin{array}{c}
\includegraphics[width=6.7cm,height=6.2cm]{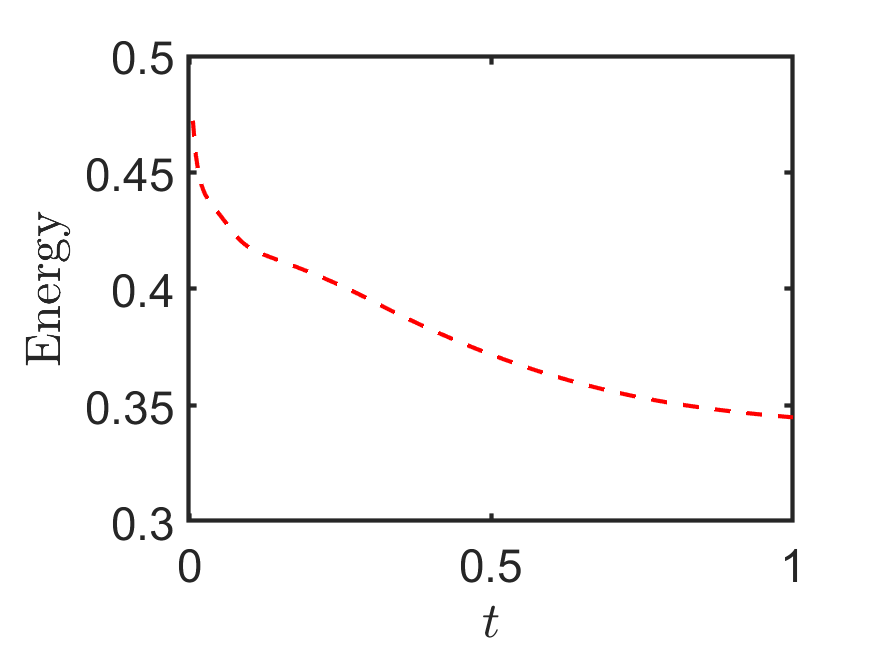} \;
\includegraphics[width=6.7cm,height=6.2cm]{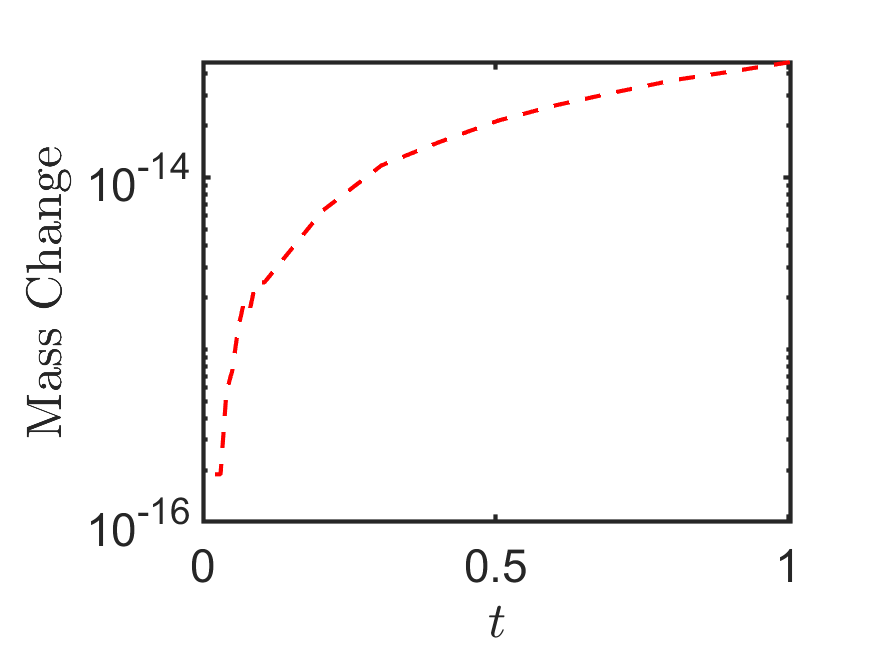} \;
\end{array}$\vspace{-0.2cm}
\caption{$\mathbf{Example\ \ref{examThreeD}}$, the discrete energy history and discrete mass variation.}\label{Cexp3DEandM}
\end{figure}

\end{example}

\section{Conclusion}
In this work, we studied three types of first- and second-order IEQ-FEMs for solving the CHNS equations. The first two types position the intermediate function introduced by the IEQ approach either in the continuous function space or in a combination of the continuous function space and finite element space. These methods each exhibit distinct advantages: the former is computationally fast, though the energy may not be fully stable in the FEM space; the latter ensures unconditional energy stability in the FEM space, but the projection step incurs additional computational costs. The third IEQ-FEM scheme is a hybrid of the two, designed to leverage the benefits of both. It begins with the former IEQ-FEM for computational efficiency but switches to the latter if the energy fails to satisfy the energy dissipation property. Numerical examples were presented to validate the theoretical findings and demonstrate the effectiveness of the proposed methods.

{\color{black}A rigorous error analysis for CHNS equations of the proposed methods is desirable and is deferred to future work. Besides, extending the proposed methods to higher-order time discretizations is a potential direction for future research. }

	\section*{Acknowledgments}
This work was supported by the National Natural Science Foundation of China Project (No. 12571469), Scientific Research Innovation Capability Support Project for Young Faculty of China (No. SRICSPYF-BS2025132), the Project of Scientific Research Fund of the Hunan Provincial Science and Technology Department (No. 2024JJ1008), the 111 Project (No. D23017), and Program for Science and Technology Innovative Research Team in Higher Educational Institutions of Hunan Province of China.
 
 \subsection*{Author Contributions}
 Authors are listed in alphabetical order by last name. All authors contribute equally to this manuscript.
	
	\subsection*{Data availability}
	Not Applicable.
	
	\section*{Declarations}
	\subsection*{Conflict of interests}
	The authors declare no competing interests.
	
%
%
%
%
%

\bibliographystyle{plain}
\bibliography{myref}
	
\end{document}